\documentclass[12pt,a4paper]{amsart}
\usepackage{amscd,amsmath,amssymb,amsfonts,times}
\usepackage{mathrsfs}
\usepackage[all]{xy}
   \newtheorem{thm}{Theorem}[section]
   \newtheorem{lem}[thm]{Lemma}

   \newtheorem{set}[thm]{Setting}
   \newtheorem{prop}[thm]{Proposition}
   \newtheorem{cor}[thm]{Corollary}
   
   \newtheorem{conj}[thm]{Conjecture}
   \newtheorem{defn}[thm]{Definition}
   
   \newtheorem{rem}[thm]{Remark}
   \newtheorem{athm}{Theorem}[subsection]
   
   \newtheorem{aprop}[athm]{Proposition}
   \newtheorem{acor}[athm]{Corollary}
   \newtheorem{acond}[athm]{Condition}
   \newtheorem{adefn}[athm]{Definition}
   \newtheorem{arem}[athm]{Remark}
\numberwithin{equation}{thm}
\newenvironment{pf}{\par\smallskip\noindent\emph{Proof.\;}}{\qed\par\medskip}
\newenvironment{pf*}[1]{\par\smallskip\noindent\emph{#1.\;}}{\qed\par\medskip}

\begin{document}
\title[Cycle class and $p$-adic regulator]{Cycle classes for $\boldsymbol{p}$-adic \'etale Tate twists \endgraf \vspace{2pt} and the image of $\boldsymbol{p}$-adic regulators}
\author[K. Sato]{Kanetomo Sato}
\address{Department of Mathematics, Chuo University \endgraf 
 1-13-27 Kasuga, Tokyo 112-8551, JAPAN \endgraf
Tel: +81-3-3817-1745 \qquad Fax: +81-3-3817-1746}
\email{kanetomo@math.chuo-u.ac.jp}
\date{March 12, 2012}
\thanks{ \endgraf
supported by Grant-in-Aid for Scientific Research B-20740008 and JSPS Core-to-Core Program 18005}
\begin{abstract}
In this paper, we construct Chern class maps and cycle class maps with values in $p$-adic \'etale Tate twists \cite{S2}. We also relate the $p$-adic \'etale Tate twists with the finite part of Bloch-Kato. As an application, we prove that the integral part of $p$-adic regulator maps has values in the finite part of Galois cohomology under certain assumptions.
\end{abstract}
\maketitle
%
%
%
%
%
%
%
%
%
%
\def\Ab{\mathbf{A\hspace{-0.8pt}b}}
\def\abs{{\sf abs}}
\def\bfmod{\text{-}\mathbf{mod}}
\def\BGL{{\sf B_{\star}GL}}
\def\can{{\sf can}}
\def\cd{{\sf cd}}
\def\codim{{\sf codim}}
\def\ch{{\sf ch}}
\def\CH{{\sf CH}}
\def\cl{{\sf cl}}
\def\Coker{{\sf Coker}}
\def\cone{{\sf Cone}}
\def\cont{{\sf cont}}
\def\crys{{\sf crys}}
\def\cocrys{{\sf cont-cr}}
\def\dccrys{{\sf cG-cr}}
\def\div{{\sf div}}
\def\dlog{d{{\sf log}}}
\def\dR{{\sf dR}}
\def\dX{d}
\def\dY{N}
\def\Et{{\sf \acute{E}t}}
\def\et{{\sf \acute{e}t}}
\def\coet{{\sf cont-\acute{e}t}}
\def\F{{\sf F}}
\def\Fil{{\sf Fil}}
\def\Frac{{\sf Frac}}
\def\Gal{{\sf Gal}}
\def\GL{{\sf GL}}
\def\gp{{\sf gp}}
\def\gr{{\sf gr}}
\def\gys{{\sf gys}}
\def\h{{\sf h}}
\def\hom{{\sf hom}}
\def\Hom{{\sf Hom}}
\def\sHom{{\mathscr H}\hspace{-2pt}om}
\def\id{{\sf id}}
\def\Image{{\sf Im}}
\def\ker{{\sf Ker}}
\def\logcrys{{\sf log\text{-}crys}}
\def\Ob{{\sf Ob}}
\def\op{{\sf op}}
\def\ord{{\sf ord}}
\def\pr{{\sf pr}}
\def\Proj{{\sf Proj}\,}
\def\qc{{\sf qc}}
\def\qis{{\sf qis}}
\def\rank{{\sf rank}}
\def\reg{{\sf reg}}
\def\sh{\,\text{\rm sh}}
\def\Spec{{\sf Spec}}
\def\st{{\sf st}}
\def\syn{{\sf syn}}
\def\cosyn{{\sf cont-syn}}
\def\dcsyn{{\sf cG-syn}}
\def\tr{{\sf tr}}
\def\val{{\sf val}}
\def\zar{{\sf Zar}}
\def\bA{\mathbb A}
\def\bC{\mathbb C}
\def\bE{\mathbb E}
\def\bF{\mathbb F}
\def\bH{\mathbb H}
\def\bJ{\mathbb J}
\def\bK{\mathbb K}
\def\bL{\mathbb L}
\def\bN{\mathbb N}
\def\bP{\mathbb P}
\def\bQ{\mathbb Q}
\def\bR{\mathbb R}
\def\bT{\mathbb T}
\def\bZ{\mathbb Z}

\def\cB{\mathscr B}
\def\cD{\mathscr D}
\def\cE{\mathscr E}
\def\cF{\mathscr F}
\def\cH{\mathscr H}
\def\cI{\mathscr I}
\def\cK{\mathscr K}
\def\cL{\mathscr L}
\def\cM{\mathcal M}
\def\cN{\mathcal N}
\def\cO{\mathscr O}
\def\ccS{\mathcal V}
\def\cS{\mathscr S}
\def\cU{\mathscr U}
\def\cV{\mathscr V}
\def\cX{\mathscr X}
\def\cY{\mathscr Y}
\def\cZ{\mathscr Z}

\def\fD{\mathfrak D}
\def\fE{\mathfrak E}
\def\fF{\mathfrak F}
\def\fL{\mathfrak L}
\def\fM{\mathfrak M}
\def\fO{\mathfrak O}
\def\fT{\mathfrak T}
\def\fY{\mathfrak Y}
\def\fb{\mathfrak b}
\def\fd{\mathfrak d}
\def\fe{\mathfrak e}
\def\fk{\mathfrak k}
\def\fl{\mathfrak l}
\def\fm{\mathfrak m}
\def\fn{\mathfrak n}
\def\fo{\mathfrak o}
\def\fp{\mathfrak p}

\def\sfc{{\sf c}}
\def\sfK{{\sf K}}

%
\def\ep{\epsilon}
\def\ka{\kappa}
\def\lam{\lambda}
\def\ze{\zeta}
\def\UK#1{\cU^{\hspace{-.8pt}#1}\hspace{-2.0pt}\cK^M_r}
\def\UKo#1{\cU^{\hspace{-.8pt}#1}\hspace{-2.0pt}\cK^M_{r-1}}
\def\UM#1{\cU^{\hspace{-.8pt}#1}\hspace{-1.5pt}M_n^r}
\def\VM#1{\cV^{\hspace{-.8pt}#1}\hspace{-1.5pt}M_n^r}
\def\UMo#1{\cU^{\hspace{-.8pt}#1}\hspace{-1.5pt}M_1^r}
\def\VMo#1{\cV^{\hspace{-.8pt}#1}\hspace{-1.5pt}M_1^r}
\def\ijcO{i^*\hspace{-.7pt}j_*\cO}
\def\ihcO{i^*\hspace{-.7pt}h_*\cO}

\def\ve{\varepsilon}
\def\vL{\varLambda}
\def\vD{\varDelta}
\def\vG{\varGamma}
\def\vS{\varSigma}
\def\vT{\varTheta}
\def\vt{\vartheta}
\def\nvD#1{\vD^{#1}}
\def\vH#1{D^{#1}}
\def\bvD#1{\ol{\vD^{#1}}}
\def\Ri#1{R^{#1}\hspace{-1pt}\varprojlim}

%
%
%
%
\def\ra{\rightarrow}
\def\lra{\longrightarrow}
\def\Lra{\Longrightarrow}
\def\la{\leftarrow}
\def\lla{\longleftarrow}
\def\Lla{\Longleftarrow}
\def\da{\downarrow}
\def\hra{\hookrightarrow}
\def\lmt{\longmapsto}
\def\sm{\setminus}
\def\bs#1{\boldsymbol{#1}}
\def\wt#1{\widetilde{#1}}
\def\wh#1{\widehat{#1}}
\def\spt{\sptilde}
\def\ol#1{\overline{#1}}
\def\ul#1{\underline{#1}}
\def\us#1#2{\underset{#1}{#2}}
\def\os#1#2{\overset{#1}{#2}}
\def\lim#1{\us{#1}{\varinjlim}}
\def\tcB{\wt{\cB}}
\def\tcZ{\wt{\cZ}}
\def\tD{\wt{D}}
\def\zp{{\bZ_p}}
\def\qp{{\bQ_p}}
\def\zpn{{\bZ/p^n\bZ}}

\def\Xggf{X_{\ol K}}
\def\Xgf{X_K}
\def\Gm{{\mathbb G}_{\hspace{-1pt}{\sf m}}}
\def\isom{\hspace{9pt}{}^\sim\hspace{-16.5pt}\lra}
\def\lisom{\hspace{10pt}{}^\sim\hspace{-17.5pt}\lla}
%
%
%
\def\Wn{W_{\hspace{-2pt}n}}
\def\witt#1#2#3{W_{\hspace{-2pt}#2}{\hspace{1pt}}\Omega_{#1}^{#3}}
\def\logwitt#1#2#3{W_{\hspace{-2pt}#2}{\hspace{1pt}}\Omega_{{#1},{\log}}^{#3}}
\def\cohwitt#1#2#3{W_{\hspace{-2pt}#2}{\hspace{1pt}}\vL_{#1}^{#3}}
\def\homwitt#1#2#3{W_{\hspace{-2pt}#2}{\hspace{1pt}}\varXi_{#1}^{#3}}
\def\llwitt#1#2#3#4{W_{\hspace{-2pt}#3}{\hspace{1pt}}\Omega_{#1}^{#4}(\log #2)_{\log}}
\def\mwitt#1#2#3{W_{\hspace{-2pt}#2}{\hspace{1pt}}\omega_{#1}^{#3}}
\def\mlogwitt#1#2#3{W_{\hspace{-2pt}#2}{\hspace{1pt}}\omega_{{#1},{\log}}^{#3}}
\def\mtwitt#1#2#3{W_{\hspace{-2pt}#2}{\hspace{1pt}}\wt{\omega}{}_{#1}^{#3}}
\def\loglogwitt#1#2{W_{\hspace{-2pt}n}{\hspace{1pt}}\Omega^{#2}_{Y^{(#1)}\hspace{-1.5pt},\hspace{1pt}E^{[#1]}\hspace{-1pt},\hspace{1pt}\log}}
%
%
\vspace{-5pt}
\section{\bf Introduction}
\medskip
Let $p$ be a rational prime number. Let $A$ be a Dedekind ring whose fraction field has characteristic zero and which has a residue field of characteristic $p$. Let $X$ be a noetherian regular scheme of pure-dimension which is flat of finite type over $S:=\Spec(A)$ and which is a smooth or semistable family around its fibers over $S$ of characteristic $p$. Extending the idea of Schneider \cite{Sch}, the author defined in \cite{S2} the objects $\fT_n(r)_X$ ($r,n\ge 0$) of the derived category of \'etale $\bZ/p^n$-sheaves on $X$ playing the role of the $r$-th Tate twist with $\bZ/p^n$-coefficients, which are endowed with a natural product structure with respect to $r$ and both contravariantly and covariantly functorial (i.e., there exist natural pull-back and trace morphisms) for arbitrary separated $S$-morphisms of such schemes. Those pull-back and trace morphisms satisfy a projection formula.

The first aim of this paper is to construct the following Chern class map and cycle class map for $m,r \ge 0$:
\[\xymatrix{
K_m(X) \ar[rrd]^-{\hspace{15pt}\sfc_{r,m} \;\;\text{(Chern class)}} \ar@{.>}[d]_{\text{(missing)}} \\
 \CH^r(X,m) \ar[rr]^-{\cl_X^{r,m}}_-{\text{(cycle class)}} && H^{2r-m}_\et(X,\fT_n(r)_X),
}\]
where $K_m(X)$ denotes the algebraic $K$-group \cite{Q} and $\CH^r(X,m)$ denotes the higher Chow group \cite{B2}. The Chern classes with values in higher Chow groups have not been defined in this arithmetic situation for the lack of a product structure on them, which we do not deal with in this paper. Therefore we will construct the above Chern class map and cycle class map independently. By the general framework due to Gillet \cite{Gi}, the existence of $\sfc_{r,m}$ is a rather direct consequence of the product structure on $p$-adic \'etale Tate twists $\{\fT_n(r)_X\}_{r \ge 0}$ and the Dold-Thom isomorphism (cf.\ Theorem \ref{thm3-1}). On the other hand, the construction of $\cl_X^{r,m}$ is more delicate, because the \'etale cohomology with $\fT_n(r)$-coefficients does not satisfy homotopy invariance. To overcome this difficulty, we introduce a version of $p$-adic \'etale Tate twists with log poles along horizontal normal crossing divisors (cf.\ \S\ref{sect2}, see also \cite{JS} p.\ 517), and prove a $p$-adic analogue of the usual homotopy invariance (cf.\ Corollary \ref{cor3-1}) and a semi-purity property (cf.\ Theorem \ref{cor4-1}) for this new coefficient. These results will enable us to pursue an analogy with Bloch's construction of cycle class maps to obtain $\cl_X^{r,m}$ (cf.\ \S\ref{sect5}). This `higher' cycle class map will be a fundamental object to study in `higher' higher classfield theory \cite{Sai}. We will mention a local behavior of $\cl_X^{r,m}$ in Remark \ref{rem5-1} below.
\par
The second aim is of this paper is to relate the $p$-adic \'etale Tate twists with the finite part of Galois cohomology \cite{BK2}, using the Fontaine-Jannsen conjecture proved by Hyodo, Kato and Tsuji (\cite{HK}, \cite{K4}, \cite{T1}, cf.\ \cite{Ni2}). Assume here that $A$ is a $p$-adic integer ring and that $X$ is projective over $A$ with strict semistable reduction. Let $K$ be the fraction field of $A$ and put $\Xgf:=X \otimes_A K =X[p^{-1}]$. We define
\begin{align*}
 H^i(X,\fT_\qp(r)_X) &:=\qp \otimes_{\zp} \varprojlim_{n \ge 1} \ H^i_\et(X,\fT_n(r)_X), \\
 H^i(\Xgf,\qp(r)) &:=\qp \otimes_{\zp} \varprojlim_{n \ge 1} \ H^i_\et(\Xgf,\mu_{p^n}^{\otimes r}),
\end{align*} 
where $\mu_{p^n}$ denotes the \'etale sheaf of $p^n$-th roots of unity on $\Xgf$, i.e, the usual Tate twist on $\Xgf$. We have a natural restriction map $H^i(X,\fT_\qp(r)_X) \to H^i(\Xgf,\qp(r))$ and a canonical descending filtration $F^\bullet$ on $H^i(\Xgf,\qp(r))$ resulting from the Hochschild-Serre spectral sequence for the covering $\Xggf:=\Xgf \otimes_K \ol K \to \Xgf$ (cf.\ \eqref{eq8-hs}). We define a (not necessarily exhaustive) filtration $F^\bullet$ on $H^i(X,\fT_\qp(r)_X)$ as the inverse image of $F^\bullet$ on $H^i(\Xgf,\qp(r))$, which induces obvious inclusions for $m \ge 0$
\[\xymatrix{ \gr_F^m H^i(X,\fT_\qp(r)_X) \; \ar@<-1.5pt>@{^{(}->}[r] & \gr_F^m H^i(\Xgf,\qp(r)) \simeq H^m(K, H^{i-m}(\Xggf,\qp(r))). }\]
Here $H^*(K,-)$ denotes the continuous Galois cohomology of the absolute Galois group $G_K=\Gal(\ol K/K)$ defined by Tate \cite{Ta}. We will prove that $\gr_F^1 H^i(X,\fT_\qp(r)_X)$ agrees with the finite part $H^1_f(K, H^{i-1}(\Xggf,\qp(r)))$ under the assumptions that $p$ and $r$ are sufficiently large and that the monodromy-weight conjecture \cite{M} holds for the log crystalline cohomology of the reduction of $X$ in degree $i-1$ (see Theorem \ref{thm8-2} below, cf.\ \cite{Na}). This result is an extension of the $p$-adic point conjecture (\cite{Sch}, \cite{LS}, \cite{N1}) to the semistable reduction case and gives an `unramified' version of results of Langer \cite{L} and Nekov\'a\v{r} \cite{N2} relating the log syntomic cohomology of $X$ with the geometric part $H^1_g(K, H^{i-1}(\Xggf,\qp(r)))$.
\par
There is an application of the above results as follows. Let $K$ be a number field and let $V$ be a proper smooth geometrically integral variety over $K$. Put $\ol V := V \otimes_K \ol K$. Let $i$ and $r$ be non-negative integers with $2r \ge i+1$, and let $p$ be a prime number. The \'etale Chern characters (cf.\ \cite{So1}) induce $p$-adic regulator maps
\begin{align*}
  \reg_p^{2r-i-1,r} & : K_{2r-i-1}(V)_{\fo} \lra H^1(K,H^i(\ol V, \qp(r))) \qquad\, \hbox{($2r > i+1$)}, \\
  \reg_p^{0,r} & : K_0(V)_\hom \lra H^1(K,H^{2r-1}(\ol V, \qp(r))) \qquad \hbox{($2r = i+1$)}.
\end{align*}
Here $K_0(V)_\hom$ denotes the homologically trivial part of $K_0(V)$, and $K_m(V)_{\fo}$ denotes the integral part of $K_m(V)$ in the sense of Scholl (see \S\ref{sect7} below). Motivated by the study of special values of $L$-functions, Bloch and Kato \cite{BK2} conjecture that the image of $\reg_p^{2r-i-1,r}$ is contained in the finite part $H^1_f(K,H^i(\ol V, \qp(r)))$ and spans it over $\qp$. In the direction of this conjecture, we will prove the following result, which extends a result of Nekov\'a\v{r} \cite{N2} Theorem 3.1 on $\reg_p^{0,r}$ to the case $2r > i+1$ and extends a result of Niziol \cite{Ni} on the potentially good reduction case to the general case:
\begin{thm}[{{\bf \S\ref{sect7}}}]\label{thm1-1}
Assume $i+1<2r \le 2(p-2)$ and the monodromy-weight conjecture for the log crystalline cohomology of degree $i$ of projective strict semistable varieties over $\bF_p$.
Then the  image of $\reg_p^{2r-i-1,r}$ is contained in $H^1_f(K,H^i(\ol V, \qp(r)))$.
\end{thm}
\noindent
Here a projective strict semistable variety over $\bF_p$ means the reduction of a regular scheme which is projective flat over a $p$-adic integer ring with strict semistable reduction. We use the alteration theorem of de Jong \cite{dJ} in the proof this theorem, and the projective strict semistable varieties concerned in the assumption are those obtained from alterations of scalar extensions of $V$ to the completion of $K$ at its places dividing $p$.
\par
\medskip
This paper is organized as follows. In \S\ref{sect1}, we introduce cohomological and homological logarithmic Hodge-Witt sheaves with horizontal log poles on normal crossing varieties over a field of characteristic $p >0$. In \S\ref{sect2}, we define $p$-adic \'etale Tate twists with horizontal log poles, and construct a localization sequence using this object (Theorem \ref{prop2-2}). In \S\ref{sect3}, we prove the Dold-Thom isomorphisms and define the Chern class maps for $p$-adic \'etale Tate twists. The sections \ref{sect4} and \ref{sect5} will be devoted to the construction of cycle class maps for $p$-adic \'etale Tate twists. In \S\ref{sect6} we will introduce Hodge-Witt cohomology and homology of normal crossing varieties and prove that the monodromy-weight conjecture implies a certain invariant cycle theorem. In \S\ref{sect8}, we establish the comparison between $p$-adic \'etale Tate twists and the finite part of Bloch-Kato. We will prove Theorem \ref{thm1-1} in \S\ref{sect7}. In the appendix, we will formulate a continuous crystalline cohomology and a continuous syntomic cohomology and prove several technical compatibility results which will have been used in \S\ref{sect8}.

\par
\medskip
The author expresses his gratitude to Professors Shuji Saito and Takeshi Tsuji, and Masa-nori Asakura for valuable comments and discussions on the subjects of this paper.
\par
\medskip
\section*{\bf Notation}
\par
\medskip
For an abelian group $M$ and a positive integer $n$, ${}_nM$ and $M/n$ denote the kernel and the cokernel of the map $M \os{\times n}{\lra} M$, respectively. For a field $k$, $\ol k$ denotes a fixed separable closure, and $G_k$ denotes the absolute Galois group $\Gal(\ol k/k)$. For a topological $G_k$-module $M$, $H^*(k,M)$ denote the continuous Galois cohomology groups $H^*_\cont(G_k,M)$ in the sense of Tate \cite{Ta}. If $M$ is discrete, then $H^*(k,M)$ agree with the \'etale cohomology groups of $\Spec(k)$ with coefficients in the \'etale sheaf associated with $M$.
\par
Unless indicated otherwise, all cohomology groups of schemes are taken over the \'etale topology. For a scheme $X$, an \'etale sheaf $\cF$ on $X$ (or more generally an object in the derived category of sheaves on $X_{\et}$) and a point $x \in X$, we often write $H^*_x(X,\cF)$ for $H^*_x(\Spec(\cO_{X,x}),\cF)$. 
For a positive integer $m$ which is invertible on $X$, we write $\mu_m$ for the \'etale sheaf of the $m$-th roots of unity on $X$. For a prime number $p$ which is invertible on $X$ and integers $m,r \ge 0$, we define
\begin{align*}
H^m(X,\zp(r)) &:=\varprojlim_{n \ge 1} \ H^i(X,\mu_{p^n}^{\otimes r}), \\
H^m(X,\qp(r)) &:=\qp \otimes_\zp H^m(X,\zp(r)).
\end{align*}
For a pure-dimensional scheme $X$ and a non-negative integer $q$, we write $X^q$ for the set of all points on $X$ of codimension $q$.

\medskip
\newpage
\section{\bf Logarithmic Hodge-Witt sheaves}\label{sect1}
\medskip
We first fix the following terminology.
\begin{defn}\label{def1-1}
\begin{enumerate}
\item[(1)]
A normal crossing varity over a field $k$ is a pure-dimensional scheme which is separated of finite type over $k$ and everywhere \'etale locally isomorphic to \[ \Spec\big(k[t_0,\dotsc,t_\dY]/(t_0\dotsb t_a)\big) \quad \hbox{for some $0 \le a \le \dY := \dim(Y)$}.\]
\item[(2)]
We say that a normal crossing variety $Y$ is {\it simple} if all irreducible components of $Y$ are smooth over $k$.
\item[(3)]
An admissible divisor on a normal crossing varity $Y$ is a reduced effective Cartier divisor $D$ such that the immersion $D \hra Y$ is everywhere \'etale locally isomorphic to \[\xymatrix{ \Spec\big(k[t_0,\dotsc,t_\dY]/(t_0\dotsb t_a,t_{a+1} \dotsb t_{a+b})\big) \;\ar@{^{(}->}[r] & \Spec\big(k[t_0,\dotsc,t_\dY]/(t_0\dotsb t_a)\big) }\] for some $a,b \ge 0$ with $a+b \le \dY=\dim(Y)$.
\end{enumerate}
\end{defn}
Let $k$ be a field of characteristic $p$.
Let $Y$ be a normal crossing variety over $k$, and let $D$ be an admissible divisor on $Y$. Put $V:=Y-D$ and let $f$ and $g$ be as follows: \[\xymatrix{D \;\ar@{^{(}->}[r]^f & Y & \ar@{_{(}->}[l]_{g \qquad\;\;} \; V=Y-D. }\] For $r \ge 0$, we define \'etale sheaves $\nu_{(Y,D),n}^r$ and $\lam_{(Y,D),n}^r$ on $Y$ as follows:
\begin{align*}
\lam_{(Y,D),n}^r & := \Image \left(\dlog: (g_*\cO_V^\times)^{\otimes r} \to \bigoplus{}_{x \in V^0} \ i_{x*}\logwitt x n r \right) \\
\nu_{(Y,D),n}^r  & := \ker \left( \partial: \bigoplus{}_{x \in V^0} \ i_{x*}\logwitt x n r \to \bigoplus{}_{x \in V^1}\ i_{x*}\logwitt x n {r-1} \right),
\end{align*}
where for $x \in V$, $i_x$ denotes the composite map $x \hra V \hra Y$ and the arrow $\partial$ denotes the sum of boundary maps due to Kato \cite{K1}. We define
\[ \lam_{(Y,D),n}^r = \nu_{(Y,D),n}^r =  0 \quad \hbox{ for } \;\; r < 0. \]
When $D=\emptyset$, we put
\[ \nu_{Y,n}^r := \nu_{(Y,\emptyset),n}^r \quad \hbox{ and } \quad \lam_{Y,n}^r:=\lam_{(Y,\emptyset),n}^r,\]
which have been studied in \cite{S1}. Although we assumed the perfectness of $k$ in \cite{S1}, all the local results are extended to the case that $k$ is not necessarily perfect by the Gersten resolution of $\logwitt Y n r$ for smooth $Y$ due to Shiho \cite{Sh} (cf.\ \cite{GS}).
We have $\nu_{(Y,D),n}^r=g_*\nu_{V,n}^r$ by the left exactness of $g_*$, and $\nu_{D,n}^{r-1} \simeq R^1f^!\nu_{Y,n}^r$ by the purity of $\nu_{Y,n}^r$ (\cite{S1} Theorem 2.4.2). Hence there is a short exact sequence
\begin{equation}\label{eq1-1}
0 \lra \nu_{Y,n}^r \lra \nu_{(Y,D),n}^r \lra f_*\nu_{D,n}^{r-1} \lra 0.
\end{equation}
By this fact, we have \[ \lam_{(Y,D),n}^r=\nu_{(Y,D),n}^r=g_*\logwitt V n r\] if $Y$ is smooth
(loc.\ cit.\ (2.4.9)), which we denote by $\logwitt {(Y,D)} n r$. The following proposition is useful later, where $Y$ is not necessarily smooth:
\begin{prop}\label{lem1-2}
Assume that $Y$ is simple, and let $Y_1,Y_2,\dotsb, Y_q$ be the distinct irreducible components of $Y$.
Then there is an exact sequence on $Y_{\et}$ \begin{align*} 0 \lra \lam_{(Y,D),n}^r  \os{\check{r}^0}\lra \bigoplus_{|I|=1} \ \logwitt {(Y_I,D_I)} n r & \os{\check{r}^1}\lra \bigoplus_{|I|=2} \ \logwitt {(Y_I,D_I)} n r \os{\check{r}^2}{\lra} \\ \dotsb & \os{\check{r}^{q-1}}{\lra} \bigoplus_{|I|=q} \ \logwitt {(Y_I,D_I)} n r \lra 0, \end{align*}
where the notation `$|I|=t$' means that $I$ runs through all subsets of $\{1,2,\dotsc,q\}$ consisting of $t$ elements, and for such $I=\{i_1,i_2,\dotsc,i_t \}$\,{\rm(}$i_j$'s are pair-wise distinct{\rm)}, we put
\[ Y_I :=Y_{i_1} \cap Y_{i_2} \cap \dotsb \cap Y_{i_t} \quad \hbox{ and } \quad D_I := D \times_Y Y_I. \]
The arrow $\check{r}^0$ denotes the natural restriction map. For $(I,I')$ with $I=\{i_1,i_2,\dotsc,i_t\}$ {\rm($i_1<i_2<\dotsb<i_t$)} and $|I'|=t+1$, the $(I,I')$-factor of $\check{r}^t$ is defined as
\[ \begin{cases} 0 & \hbox{{\rm(}if $I \not\subset I'${\rm)}} \\
 (-1)^{t-a} \cdot (\beta_{I'\hspace{-1pt}I})^* & \hbox{{\rm(}if $I'=I \cup \{ i_{t+1} \}$ and $i_1<\dotsb<i_a < i_{t+1} < i_{a+1} < \dotsb < i_t${\rm)}}, \end{cases}  \] where $\beta_{I'\hspace{-1pt}I}$ denotes the closed immersion $Y_{I'} \hra Y_I$.
\end{prop}
\noindent
We need the following lemma to prove this proposition:
\begin{lem}\label{lem1-1}
Assume that $(Y,D,V)$ fits into cartesian squares of schemes
\[ \xymatrix{ D \; \ar@{^{(}->}[r]^{f} \ar[d]_{\iota} \ar@{}[rd]|{\square} & Y \ar[d]_i \ar@{}[rd]|{\square} & \; V \ar@{_{(}->}[l]_g \ar[d]_{i'} \\ \cD \; \ar@{^{(}->}[r] & \cY & \ar@{_{(}->}[l]_h \; \cV}  \]
such that $\cY$ is regular, such that the vertical arrows are closed immersions and such that $Y$, $\cD$ and $Y \cup \cD$ are simple normal crossing divisors on $\cY$, where we put $\cV:=\cY-\cD$. Then the pull-back map $i^*h_*\cO_{\cV}^\times \to g_*\cO_V^\times$ on $Y_\et$ is surjective.
\end{lem}
\begin{pf*}{Proof of Lemma \ref{lem1-1}}
We use the same notation as in Proposition \ref{lem1-2}. For a Cartier divisor $E$ on a scheme $Z$, let $c_1^Z(E) \in H^1_{|E|}(Z,\cO_Z^{\times})$ be the localized 1st Chern class of the invertible sheaf $\cO_Z(E)$. Since $\cY$ is regular, $h_*\cO_{\cV}^\times$ is generated by $\cO_{\cY}^\times$ and local uniformizers of the irreducible components $\{\cD_j\}_{j \in J}$ of $\cD$. Put
\[ D_j:=Y \times_{\cY} \cD_j \qquad (j \in J), \]
which is an admissible divisor on $Y$. Since $(g_*\cO_V^\times)/\cO_Y^{\times} \simeq R^1f^!\cO_Y^{\times}$, it is enough to show that the Gysin map \[\varphi : \bigoplus_{j \in J}\ \bZ_{D_j} \lra R^1f^!\cO_Y^{\times}\] sending $1 \in \bZ_{D_j}$ to $c_1^Y(D_j)$ is bijective on $D_{\et}$.
By \cite{S1} Lemma 3.2.2, there is an exact sequence on $Y_{\et}$
\begin{equation}\label{eq7-3} 0 \lra \cO_Y^{\times} \os{\check{r}^0}{\lra} \bigoplus_{|I|=1} \ \cO_{Y_I}^{\times} \os{\check{r}^1}{\lra} \bigoplus_{|I|=2} \ \cO_{Y_I}^{\times} \os{\check{r}^2}{\lra} \dotsb \os{\check{r}^{q-1}}{\lra} \bigoplus_{|I|=q} \ \cO_{Y_I}^{\times} \lra 0,
\end{equation}
where $\check{r}$'s are defined in the same way as $\check{r}$'s in Proposition \ref{lem1-2}, and $q$ denotes the number of the distinct irreducible components of $Y$. Since $f^!\cO_{Y_I}^{\times}=0$ for any non-empty subset $I \subset \{1,2,\dotsc,q \}$, the exactness of \eqref{eq7-3} implies that of the lower row of the following commutative diagram with exact rows:
\[\xymatrix{
\us{\phantom{|}}0 \ar@<5.2pt>[r] & \displaystyle \bigoplus_{\phantom{|}j \in J\phantom{|}}\ \bZ_{D_j} \ar@<5.2pt>[r]^-{\check{r}^0} \ar[d]_{\varphi} & \displaystyle \bigoplus_{j \in J} \ \bigoplus_{|I|=1} \ \bZ_{D_j \cap Y_I} \ar[d]_{\wr\hspace{-1.5pt}} \ar@<5.2pt>[r]^-{\check{r}^1} & \displaystyle \bigoplus_{j \in J} \ \bigoplus_{|I|=2} \ \bZ_{D_j \cap Y_I} \ar[d]_{\wr\hspace{-1.5pt}} \\
\us{\phantom{|}}0 \ar@<5.2pt>[r] & \us{\phantom{|}}{R^1f^!\cO_Y^{\times}} \ar@<5.2pt>[r]^-{\check{r}^0} & \displaystyle \bigoplus_{|I|=1} \ R^1f_1^!\cO_{Y_I}^{\times} \ar@<5.2pt>[r]^{\check{r}^1} & \displaystyle \bigoplus_{|I|=2} \ R^1f_1^!\cO_{Y_I}^{\times}\,,
}\]
where $D_j \cap Y_I$ is regular for each $j \in J$ and $I \subset \{1,2,\dotsc,q \}$ by the assumption that $Y \cup \cD$ has simple normal crossings on $\cX$. The middle and the right vertical arrows are defined in the same way as for $\varphi$, and bijective by the standard purity for $\cO^\times$ (\cite{Gr} III \S6). Hence $\varphi$ is bijective as well.
\end{pf*}
\begin{pf*}{Proof of Proposition \ref{lem1-2}}
Since the problem is \'etale local on $Y$, we may assume that $(Y,D,V)$ fits into a diagram as in Lemma \ref{lem1-1}. Then there is an exact sequence
\[ \logwitt {(\cY,\cD)} n r \lra \bigoplus_{|I|=1} \ \logwitt {(Y_I,D_I)} n r \os{\check{r}^1}{\lra} \dotsb \os{\check{r}^{q-1}}{\lra} \bigoplus_{|I|=q} \ \logwitt {(Y_I,D_I)} n r \lra 0 \] on $\cY_\et$ by Lemma \ref{lem1-1} and an induction argument on the number of components of $Y$ which is similar as for \cite{S1} Lemma 3.2.2. The assertion follows from this exact sequence.
\end{pf*}
\medskip

\section{\bf $\bs{p}$-adic \'etale Tate twists with log poles}\label{sect2}
\medskip
In \S\S\ref{sect2}--\ref{sect5}, we are mainly concerned with the following setting.
\begin{set}\label{set2-1}
Let $A$ be a Dedekind domain whose fraction field has characteristic $0$ and which has a maximal ideal of positive characteristic. Put \[ S:=\Spec(A). \]
Let $p$ be a prime number which is not invertible in $A$.
Let $X$ be a regular scheme which is flat of finite type over $A$ and whose fibers over the closed points of $S$ of characteristic $p$ are reduced normal crossing divisors on $X$. We write $Y \subset X$ for the union of those fibers.
\end{set}
Let $D \subset X$ be a normal crossing divisor such that $D \cup Y$ has normal crossings on $X$
 ($D$ may be empty).
Put $U:=X-(Y \cup D)$ and $V:=Y-(Y \cap D)$, and consider a diagram of immersions
\[ \xymatrix{ 
V \; \ar@<-1.2pt>@{^{(}->}[r] \ar[d]_g & \ar[d]_h X-D & \ar@<1.2pt>@{_{(}->}[l] \ar[ld]^{j} \; U \\
Y \;  \ar@<-1.2pt>@{^{(}->}[r]^i \ar@{}[ru]|{\square} & X  \\ } \]
Let $n$ and $r$ be positive integers. We first state the Bloch-Kato-Hyodo theorem on the structure of the sheaf $M_n^r:=i^*\hspace{-1.2pt}R^rj_*\mu_{p^n}^{\otimes r}$, which will be useful in this paper.
We define the \'etale sheaf $\cK^M_r$ on $Y$ as
\[ \cK^M_r := (\ijcO_U^\times)^{\otimes n}/J, \]
where $J$ denotes the subsheaf of $(\ijcO_U^\times)^{\otimes n}$ generated by local sections of the form $a_1 \otimes a_2 \otimes \dotsb \otimes a_r$ ($a_1,a_2,\dotsc,a_r \in \ijcO_U^\times$) with $a_s + a_t=0$ or $1$ for some $1 \le s < t \le n$. There is a homomorphism of \'etale sheaves (\cite{BK} 1.2)
\begin{equation}\label{symbol}
    \cK ^M_r \lra M_n^r,
\end{equation}
which is a geometric version of Tate's Galois symbol map. For local sections $a_1,a_2,\dotsc,a_r \in \ijcO_U^\times$, 
we denote the class of $a_1 \otimes a_2 \otimes \dotsb \otimes a_r$ in $\cK^M_r$ by $\{a_1,a_2,\dotsc,a_r\}$, and
denote the image of $\{a_1,a_2,\dotsc,a_r\} \in \cK^M_r$ under \eqref{symbol} again by $\{a_1,a_2,\dotsc,a_r\}$. We define filtrations $\cU^{\bullet}$ and $\cV^{\bullet}$ on $M_n^r$ as follows.
\begin{defn}\label{def:vcyc}
Put $\fp := \ker(\cO_X \to i_*\cO_Y)$ and $1+\fp^q:=\ker(\cO_X^\times \to (\cO_X/\fp^q)^\times)$ for $q \ge 1$.
\begin{enumerate}
\item[(1)]
We define $\UK 0$ as the full sheaf $\cK ^M_n$.
For $q \ge 1$, we define $\UK q \subset \cK ^M_r$ as the image of \,$i^*\hspace{-1.2pt}(1+\fp^q) \otimes (\ijcO_U^\times)^{ \otimes r-1}$.
\item[(2)]
For $q \geq 0$, we define $\UM q$ as the image of $\cU^q\cK ^M_r$ under the map \eqref{symbol}.
\item[(3)]
When $A$ is local and its residue field $k$ has characteristic $p$, we fix a prime element $\pi \in A$ and define $\VM q \subset M_n^r$ as the part generated by $\UM {q+1}$ and the image of $\UKo q \otimes \langle \pi \rangle$ under \eqref{symbol}.
\end{enumerate}
\end{defn}
\noindent
Let $L_{X^\circ}$ be the log structure on $X$ associated with the normal crossing divisor $Y \cup D$ (\cite{K3}), and let $L_{Y^\circ}$ be its inverse image log structure onto $Y_\et$ (loc.\ cit.\ (1.4)). The following theorem is a variant of theorems of Bloch-Kato-Hyodo (\cite{BK} Theorem 1.4, \cite{H1} Theorem 1.6), and the case $D=\emptyset$ corresponds to their theorems.
\begin{thm}\label{thm2-1}
\begin{enumerate}
\item[(1)]
The symbol map $\eqref{symbol}$ is surjective, i.e., $\UM 0=M_n^r$.
\item[(2)]
Assume that $A$ is local and that its residue field $k$ has characteristic $p$. Then there are isomorphisms
\begin{align*}
 M_n^r / \VM 0 & \isom \mlogwitt {Y^{\circ}} n r\,, \\
 \VM 0 / \UM 1 & \isom \mlogwitt {Y^{\circ}} n {r-1}\,,
\end{align*}
where $\mlogwitt {Y^{\circ}} n m$ denotes the image of the logarithmic differential map \[ \dlog : (L_{Y^{\circ}}^\gp)^{\otimes m} \lra \bigoplus_{y \in Y^0}\, i_{y*} \logwitt y n m, \] and for a point $y \in Y$, $i_y$ denotes the natural map $y \hra Y$.
\item[(3)]
Under the same assumption as in {\rm(2)}, let $e$ be the absolute ramification index of $A$, and let $L_k$ be the log structure on $\Spec(k)$ associated with the pre-log structure $\bN \to k$ sending $1 \mapsto 0$. Put $e':=pe/(p-1)$. Then for $1 \le q < e'$, there are isomorphisms
\begin{align*}
 \UMo q / \VMo q & \isom
   \begin{cases}
    \omega_{Y^{\circ}}^{r-1}
     /\cB_{Y^{\circ}}^{r-1}
       \quad & \hbox{{\rm(}$p \hspace{-4pt} \not| \, q${\rm)}},\\
    \omega_{Y^{\circ}}^{r-1}
     /\cZ_{Y^{\circ}}^{r-1} \quad &\hbox{{\rm(}$p \, | \, q${\rm)}},
   \end{cases} \\
  \VMo q / \UMo {q+1} & \isom \; \omega_{Y^{\circ}}^{r-2}/\cZ_{Y^{\circ}}^{r-2}\,.
\end{align*}
Here $\omega_{Y^{\circ}}^m$ denotes the differential module of $(Y,L_{Y^{\circ}})$ over $(\Spec(k),L_k)$ {\rm (\cite{K3} (1.7))}, and $\cB_{Y^{\circ}}^m$ {\rm(}resp.\ $\cZ_{Y^{\circ}}^m${\rm)} denotes the image of $d:\omega_{Y^{\circ}}^{m-1} \ra \omega_{Y^{\circ}}^m$ {\rm(}resp.\ the kernel of $d:\omega_{Y^{\circ}}^m \ra \omega_{Y^{\circ}}^{m+1}${\rm)}.
\item[(4)] Under the same assumption and notation as in {\rm(3)}, we have $\UMo q=\VMo q=0$ for $q \geq e'$.
\end{enumerate}
\end{thm}
\begin{pf}
Note that the irreducible components of $D$ are semistable families around the fibers of characteristic $p$ by the assumption that $Y \cup D$ has normal crossings on $X$. The assertions (1) and (2) are reduced to the case that $X$ is smooth over $S$ and that $D=\emptyset$ (i.e., the Bloch-Kato theorem)  by Tsuji's trick in \cite{T2} Proof of Theorem 5.1 and a variant of Hyodo's lemma \cite{H1} Lemma 3.5, whose details will be explained in a forthcoming paper \cite{KSS}. The assertion (4) follows from \cite{BK} Lemma 5.1.
\par
We prove (3). Let $\pi \in A$ be the fixed prime element. Let $\ol \pi \in L_{Y^{\circ}}$ be the image of $\pi$, and let $[\ol \pi] \subset L_{Y^{\circ}}$ be the subsheaf of monoids generated by $\ol \pi$. The quotient $L_{Y^{\circ}}/[\ol \pi]$ is a subsheaf of monoids of $\cO_Y$ $($with respect to the multiplication of functions$)$ generated by $\cO_Y^\times$ and local equations defining $D$ and irreducible components of $Y$. There is a surjective homomorphism
\[ \delta_m : \cO_Y \otimes \{(L_{Y^{\circ}}/[\ol \pi])^{\gp}\}^{\otimes m} \lra \omega_{Y^{\circ}}^m \]
defined by the local assignment
\[ z \otimes y_1 \otimes \dotsb \otimes y_m \longmapsto  z \cdot \dlog(y_1) \wedge \dotsb \wedge \dlog(y_m), \]
with $z \in \cO_Y$ and each $y_i \in (L_{Y^{\circ}}/[\ol \pi])^{\gp}$. The kernel of $\delta_m$ is generated by local sections of the following forms (cf.\ \cite{H1} Lemma 2.2):
\begin{enumerate}
\item[(1)]
$z \otimes y_1 \otimes \dotsb \otimes y_m$ such that $y_s$ belongs to $\cO_{\Spec(k)}^\times |_Y$ for some $1 \le s \le m$.
\item[(2)]
$z \otimes y_1 \otimes \dotsb \otimes y_m$ such that $y_s=y_t$ for some $1 \le s < t \le m$.
\item[(3)]
$\sum{}_{i=1}^{\ell}\, (a_i \otimes a_i \otimes y_1 \otimes \dotsb \otimes y_{m-1}) - \sum{}_{j=1}^{\ell'}\, (b_j \otimes b_j \otimes y_1 \otimes \dotsb \otimes y_{m-1})$
with each $a_i, b_j \in L_{Y^{\circ}}/[\ol \pi]$ such that the sums $\sum{}_{i=1}^{\ell}\, a_i$ and $\sum{}_{j=1}^{\ell'}\, b_j$ taken in $\cO_Y$ belong to $L_{Y^{\circ}}/[\ol \pi]$ and satisfy ${\sum}_{i=1}^{\ell}\, a_i ={\sum}_{j=1}^{\ell'}\, b_j$.
\end{enumerate}
Hence the assertion follows from the arguments in loc.\ cit.\ p.\ 551.
\end{pf}
We define the \'etale subsheaf $FM_n^r \subset M_n^r$ as the part generated by $\UM 1$ and the image of $(\ihcO_{X-D}^\times)^{\otimes r}$, where $h$ denotes the open immersion $X-D \hra X$. By Theorem \ref{thm2-1}\,(2), Proposition \ref{lem1-2} and the same arguments as in \cite{S2} \S3.4, we obtain the following theorem:
\begin{thm}\label{thm2-2}
There are short exact sequences of sheaves on $Y_\et$
\begin{align*}
0 \lra FM_n^r \lra M_n^r & \os{\sigma}\lra \nu_{(Y,D \cap Y),n}^{r-1} \lra 0,\\
0 \lra \UM 1 \lra FM_n^r & \os{\tau}\lra \lam_{(Y,D \cap Y),n}^r \lra 0,
\end{align*}
where $\sigma$ is induced by the boundary map of Galois cohomology groups due to Kato \cite{K1}, and $\tau$ is given by the local assignment
\[\{ a_1, a_2, \dotsc, a_r \} \longmapsto \dlog(\ol{a_1}\otimes \ol{a_2} \otimes \dotsb \otimes \ol{a_r}). \]
Here $a_1, a_2, \dotsc, a_r$ are local sections of $\ihcO_{X-D}^\times$, and for $a \in \ihcO_{X-D}^\times$, $\ol a$ denotes its residue class in $g_*\cO_V^\times$. 
\end{thm}
Now we define the $p$-adic \'etale Tate twists.
\begin{defn}\label{def2-1}
Assume $r \ge 1$, and let $\cI^\bullet$ be the Godement resolution of $\mu_{p^n}^{\otimes r}$ on $U_\et$. We define a cochain complex $C_n(r)_{(X,D)}^\bullet$ of sheaves on $X_\et$ as
\[ j_*\cI^0 \to j_*\cI^1 \to \dotsb \to j_*\cI^{r-1} \to \ker\big(d:j_*\cI^r \to j_*\cI^{r+1}\big) \os{\sigma_n^r}\to i_*\nu_{(Y,D \cap Y),n}^{r-1}, \]
where $j_*\cI^0$ is placed in degree $0$ and $i_*\nu_{(Y,D \cap Y),n}^{r-1}$ is placed in degree $r+1$.
The last arrow $\sigma_n^r$ is defined as the composite map
\[\xymatrix{ \sigma_n^r : \ker\big(d:j_*\cI^r_U \to j_*\cI^{r+1}_U\big) \ar@{->>}[r] & R^rj_*\mu_{p^n}^{\otimes r} \ar[r]^-{i_*\sigma} & i_*\nu_{(Y,D \cap Y),n}^{r-1}. }\]
We write $\fT_n(r)_{(X,D)}$ for $C_n(r)_{(X,D)}^\bullet$ regarded as an object of $D^b(X_\et,\bZ/p^n)$. When $D=\emptyset$, we denote $C_n(r)_{(X,\emptyset)}^\bullet$ and $\fT_n(r)_{(X,\emptyset)}$ by $C_n(r)_X^\bullet$ and $\fT_n(r)_X$, respectively.
For $r=0$, we define $C_n(0)_{(X,D)}^\bullet:=\zpn$, the constant sheaf $\zpn$ placed in degree $0$.
\end{defn}
\begin{prop}\label{prop2-0}
For $r \ge 0$, $\fT_n(r)_{(X,D)}$ is concentrated in $[0,r]$, and there is a distinguished triangle in $D^b(X_\et,\bZ/p^n)$
\[\xymatrix{ i_*\nu_{(Y,D \cap Y),n}^{r-1}[-r-1] \ar[r] & \fT_n(r)_{(X,D)} \ar[r]^-{t} & \tau_{\le r}Rj_*\mu_{p^n}^{\otimes r} \ar[r]^-{\sigma_n^r[-r]} & i_*\nu_{(Y,D \cap Y),n}^{r-1}[-r]. }\]
Here $\nu_{(Y,D \cap Y),n}^{-1}$ means the zero sheaf when $r=0$.
\end{prop}
\begin{pf}
The first assertion follows from the surjectivity of $\sigma$ in Theorem \ref{thm2-2}. The second assertion is straight-forward.
\end{pf}
\noindent
By this proposition $\fT_n(r)_X$ defined here agrees with that in \cite{S2} \S4 by a unique isomorphism compatible with the identity map of $\mu_{p^n}^{\otimes r}$ on $U$ (loc.\ cit.\ Lemma 4.2.2). The following proposition verifies the existence of a functorial flabby resolution of $\fT_n(r)_{(X,D)}$:
\begin{prop}\label{prop2-1}
The complex $C_n(r)_{(X,D)}^\bullet$ is contravariantly functorial in the pair $(X,D)$. Here a morphism of pairs $(X,D) \to (X',D')$ means a morphism of schemes $f: X \to X'$ satisfying $f(X-D) \subset X'-D'$.
\end{prop}
\begin{pf}
The case $r \le 0$ is clear. As for the case $r \ge 1$, it is enough to show that the map
\[i_*\sigma : R^rj_*\mu_{p^n}^{\otimes r} \lra i_*\nu_{(Y,D \cap Y),n}^{r-1} \]
is contravariant in $(X,D)$. Let $f : (X,D) \to (X',D')$ be a morphism of pairs, and consider the following diagram of immersions: \[ \xymatrix{ Y' \;  \ar@<-1.2pt>@{^{(}->}[r]^{i'} & X' & \ar@<1.2pt>@{_{(}->}[l]_-{j'} \; U':=X'-D',} \] where $Y'$ denotes the union of the fibers of $X' \to S$ of characteristic $p$. By the first exact sequence in Theorem \ref{thm2-2}, $\sigma'$ ($:=\sigma$ for $(X',D')$) is surjective and $\ker(i'_*\sigma')$ maps into $\ker(i_*\sigma)$ under the base-change map
\[ f^* : f^*R^rj'_*\mu_{p^n}^{\otimes r} \lra R^rj_*\mu_{p^n}^{\otimes r}. \]
Hence this map induces a pull-back map
\begin{equation}\label{eq2-3}
f^* : f^*i'_*\nu_{(Y',D' \cap Y'),n}^{r-1} \lra i_*\nu_{(Y,D \cap Y),n}^{r-1}.
\end{equation}
These maps are obviously compatible with $\sigma$'s and satisfy transitivity. Thus we obtain the proposition.
\end{pf}

\begin{cor}\label{cor2-1}
\begin{enumerate}
\item[(1)]
Let $\ccS$ be the category whose objects are $S$-schemes satisfying the conditions in Setting {\rm\ref{set2-1}} for $X$ and whose morphisms are $S$-morphisms. Then the complexes $C_n(r)^\bullet_X=C_n(r)^\bullet_{(X,\emptyset)}$ with $X \in \Ob(\ccS)$ form a complex $C_n(r)^\bullet$ of sheaves on the big site $\ccS_\et$.
\item[(2)]
The Godement resolution $G_n(r)^\bullet_{(X,D)}$ on $X_\et$ of $C_n(r)_{(X,D)}^\bullet$ is contravariantly functorial in $(X,D)$.
\end{enumerate}
\end{cor}
\begin{rem}\label{rem2-2}
The object $\fT_n(r)_{(X,D)}$ is also contravariantly functorial in the pair $(X,D)$, that is, for a morphism of pairs $f : (X,D) \to (X',D')$, there is a unique morphism
\[ f^* : f^*\fT_n(r)_{(X',D')} \lra \fT_n(r)_{(X,D)} \quad \hbox{ in } \;\; D^b(X_\et,\bZ/p^n) \]
that extends the pull-back isomorphism of $\mu_{p^n}^{\otimes r}$ for $U \to U'$ {\rm(}cf.\ \cite{S2} Proposition {\rm4.2.8)}.
\end{rem}
\begin{rem}\label{rem2-3}
Let $\ccS$ and $C_n(r)^\bullet$ be as in Corollary {\rm\ref{cor2-1}}, and let $\fT_n(r)$ be the complex $C_n(r)^\bullet$ regarded as an object of the derived category $D^b(\ccS_\et,\zpn)$. The following facts will be useful later in \S{\rm\ref{sect3'}:}
\begin{enumerate}
\item[(1)]
There exists a unique product structure
\[ \fT_n(q) \otimes^\bL \fT_n(r) \lra \fT_n(q+r) \quad \hbox{ in } \;\; D(\ccS_\et,\zpn) \]
that extends the isomorphism $\mu_{p^n}^{\otimes q} \otimes \mu_{p^n}^{\otimes r} \isom \mu_{p^n}^{\otimes q+r}$ on the big \'etale site $S[p^{-1}]_\Et$, which follows from the same arguments as in \cite{S2} Proposition {\rm4.2.6}.
\item[(2)]
There exists a unique isomorphism
\[ \Gm \otimes^\bL \zpn[-1] \lra \fT_n(1) \quad \hbox{ in } \;\; D(\ccS_\et,\zpn) \]
that extends the canonical isomorphism $\Gm \otimes^\bL \zpn[-1] \lra \mu_{p^n}$ on $S[p^{-1}]_\Et$, which follows from the same arguments as in loc.\ cit.\ Proposition {\rm4.5.1}.
\end{enumerate}
\end{rem}

\begin{thm}\label{prop2-2}
When $D$ is regular, there is a canonical morphism
\[ \fd : \fT_n(r)_{(X,D)} \lra \alpha_*\fT_n(r-1)_D[-1] \quad \hbox{ in } \;\; D^b(X_\et,\bZ/p^n) \]
fitting into a distinguished triangle
\begin{equation}\label{eq2-2} \alpha_*\fT_n(r-1)_D[-2] \os{\alpha_*}\lra \fT_n(r)_X \os{\beta^*}\lra \fT_n(r)_{(X,D)} \os{\fd}\lra \alpha_*\fT_n(r-1)_D[-1], \end{equation}
where $\alpha$ denotes the closed immersion $D \hra X$, and $\beta$ denotes the natural morphism of pairs $(X,D) \to (X,\emptyset)$. The arrow $\alpha_*$ denotes the Gysin morphism \cite{S2} Theorem {\rm 6.1.3}.
\end{thm}
\begin{pf}
The case $r < 0$ immediately follows from the absolute purity \cite{FG}. To prove the case $r \ge 0$, we first construct the morphism $\fd$. Consider a diagram of immersions
\[ \xymatrix{ & & \ar[dl]_-j \ar[d]^u \; U \\
 Y \;  \ar@<-1.2pt>@{^{(}->}[r]^i \ar@{}[rd]|{\square} & X \ar@{}[rd]|{\square} & \ar@<1.2pt>@{_{(}->}[l]_-w \; X[p^{-1}] \\
 \ar@<1.2pt>[u] E \; \ar@<-1.2pt>@{^{(}->}[r]^{i'} & D \ar[u]_\alpha  & \; D[p^{-1}] \ar@<1.2pt>@{_{(}->}[l]_-{\varphi} \ar[u]_v. }  \]
There is a distinguished triangle on $X[p^{-1}]_\et$
\[ v_*\mu_{p^n}^{\otimes r-1}[-2] \os{v_*}\lra \mu_{p^n}^{\otimes r} \os{u^*}\lra Ru_*\mu_{p^n}^{\otimes r} \os{\fd_1}\lra v_*\mu_{p^n}^{\otimes r-1}[-1], \]
where $\fd_1$ is defined as the composite
\[ \fd_1 : Ru_*\mu_{p^n}^{\otimes r} \os{-\delta}\lra v_*Rv^!\mu_{p^n}^{\otimes r}[1] \lisom v_*\mu_{p^n}^{\otimes r-1}[-1] \]
and $\delta$ denotes the connecting morphism of a localization sequence. We used the absolute purity \cite{FG} for the last isomorphism. Applying $Rw_*$, we get a distinguished triangle on $X_\et$
\begin{equation}\label{eq2-1} \alpha_*R\varphi_*\mu_{p^n}^{\otimes r-1}[-2] \os{v_*}\lra Rw_*\mu_{p^n}^{\otimes r} \os{u^*}\lra Rj_*\mu_{p^n}^{\otimes r} \os{\fd_2}\lra \alpha_*R\varphi_*\mu_{p^n}^{\otimes r-1}[-1]. \end{equation}
Consider the following diagram with distinguished rows ($\gamma:=\alpha \circ i'$):
\[ \xymatrix{ & \fT_n(r)_{(X,D)} \ar[r]^t \ar@{.>}[d]_\fd & \tau_{\le r} Rj_*\mu_{p^n}^{\otimes r} \ar[d]^{\fd_3:=\tau_{\le r} (\fd_2)} & \\
\gamma_*\nu_{E,n}^{r-2}[-r-1] \ar[r] & \alpha_*\fT_n(r-1)_{D}[-1] \ar[r] & \tau_{\le r}(\alpha_*R\varphi_*\mu_{p^n}^{\otimes r-1}[-1]) \ar[r]^-{b'} & \gamma_*\nu_{E,n}^{r-2}[-r], } \]
where the lower triangle is distinguished by Proposition \ref{prop2-0} for $(D,\emptyset)$.
Since $\fT_n(r)_{(X,D)}$ is concentrated in $[0,r]$, we see that the composite $b' \circ \fd_3 \circ t$ is zero by Theorem \ref{thm2-2} and a simple computation on symbols. On the other hand, we have
\[ \Hom_{D^b(X_\et,\bZ/p^n)}(\fT_n(r)_{(X,D)},\gamma_*\nu_{E,n}^{r-2}[-r-1])=0, \]
again by the fact that $\fT_n(r)_{(X,D)}$ is concentrated in $[0,r]$. Hence there is a unique morphism $\fd$ fitting into the above diagram (cf.\ \cite{S2} Lemma 2.1.2\,(1)), which is the desired morphism. Finally the triangle \eqref{eq2-2} is distinguished by \eqref{eq2-1} and a commutative diagram with exact rows on $X_\et$ ($\gamma=\alpha \circ i'$)
\[\xymatrix{  & R^rw_*\mu_{p^n}^{\otimes r} \ar[r] \ar[d]^-{\sigma_X}  & R^r j_*\mu_{p^n}^{\otimes r} \ar[r] \ar[d]^-{\sigma_{(X,D)}} \ar[r]^-{\fd_3} & \alpha_*R^{r-1}\varphi_*\mu_{p^n}^{\otimes r-1} \ar[r] \ar[d]^-{\sigma_D} & 0 \\
 0 \ar[r] & i_*\nu_{Y,n}^{r-1} \ar[r] & i_*\nu_{(Y,E),n}^{r-1} \ar[r] & \gamma_*\nu_{E,n}^{r-2} \ar[r] & 0, }\]
where the surjectivity of $\fd_3$ in the upper row follows from Theorem \ref{thm2-1}\,(1) for $(D,\emptyset)$, and the exactness of the lower row follows from \eqref{eq1-1}.
\end{pf}
\begin{rem}\label{rem2-1}
Assume that $A$ is local, and let $A'$ be a Dedekind ring which is finite flat over $A$, unramified at the maximal ideals of $A[p^{-1}]$ and tamely ramified at the maximal ideals of $A$ of characteristic $p$. Then all the definitions and results for the pair $(X,D)$ in this section are extended to the scalar extension $(X \otimes_A A',D \otimes_A A')$. In fact, Theorem {\rm\ref{thm2-1}}\,{\rm(1)} and {\rm(2)} will be proved in \cite{KSS} for this generalized situation. One can check Theorem {\rm\ref{thm2-1}}\,{\rm(3)} and {\rm(4)} for $(X \otimes_A A',D \otimes_A A')$ by the same arguments as for $(X,D)$. See \cite{S2} {\rm\S3.5} for an argument to extend Theorem {\rm\ref{thm2-2}}.
\end{rem}
\medskip
\section{\bf Dold-Thom isomorphism}\label{sect3}
\medskip
Let $S,p$ and $X$ be as in Setting \ref{set2-1}. In this section we prove the Dold-Thom isomorphism for $p$-adic \'etale Tate twists. Let $E$ be a vector bundle of rank $a+1$ on $X$, and let $f : \bP:=\bP(E) \to X$ be the associated projective bundle, which is a projective smooth morphism of relative dimension $a$. Let $\cO(1)_E$ be the tautological invertible sheaf on $\bP$, and let $\xi \in H^2(\bP,\fT_n(1)_\bP)$ be the value of the 1st Chern class $c_1(\cO(1)_E) \in H^1(\bP,\cO^\times_\bP)$ under the connecting map associated with the Kummer distinguished triangle
\[  \cO_\bP^\times \lra \cO_\bP^\times \lra \fT_n(1)_\bP[1] \lra \cO_\bP^\times[1] \]
(cf.\ \cite{S2} Proposition 4.5.1). The composite morphisms
\[\xymatrix{ \fT_n(r-q)_X[-2q] \ar[r]^-{f^*} & Rf_*\fT_n(r-q)_\bP[-2i] \ar[r]^-{- \cup \xi^q } & Rf_*\fT_n(r)_\bP
  \qquad \hbox{($0 \le q \le a$)} }\]
induce a canonical morphism
\[  \gamma_E : \bigoplus_{q=0}^a \ \fT_n(r-q)_X[-2q] \lra Rf_*\fT_n(r)_\bP \;\; \hbox{ in } \;\; D^b(X_\et,\bZ/p^n). \]
\begin{thm}[{\bf Dold-Thom isomorphism}]\label{thm3-1}
$\gamma_E$ is an isomorphism for any $r \in \bZ$.
\end{thm}
\begin{pf}
$\gamma_E$ is an isomorphism outside of $Y$ by \cite{Mi} VI Theorem 10.1. The case $r<0$ follows from this fact.
To prove the case $r \ge 0$, we consider a diagram of schemes
\[\xymatrix{  P \, \ar@{^{(}->}[r]^-{\gamma} \ar@<-.35mm>[d]_g \ar@{}[rd]|{\square} & \bP \ar[d]^f \\ Y \, \ar@{^{(}->}[r]^-{i} & X & \ar@{_{(}->}[l]_-{j} \; X[p^{-1}]. }\]
We have to show that $i^*(\gamma_E)$ is an isomorphism:
\begin{equation}\label{eq3-1}
i^*(\gamma_E) : \bigoplus_{q=0}^a \ i^*\fT_n(r-q)_X[-2q] \isom Rg_*\gamma^*\fT_n(r)_\bP,
\end{equation}
where we identified $i^*Rf_*\fT_n(r)_\bP$ with $Rg_*\gamma^*\fT_n(r)_\bP$ by the proper base-change theorem.
By a standard norm argument (cf.\ \cite{S2} \S10.3) using Bockstein triangles (loc.\ cit.\ \S4.3), we are reduced to the case that $n=1$ and that $\vG(X,\cO_X)$ contains a primitive $p$-th root of unity (see also Remark \ref{rem2-1}).
We need the following lemma:
\begin{lem}\label{lem3-1}
Let $\ol{\xi} \in H^1(P,\lam_{P,1}^1)$ be the image of $\xi$ under the pull-back map
\[ \gamma^* : H^2(\bP,\fT_1(1)_\bP) \lra H^1(P,\lam_{P,1}^1) \]
{\rm (}cf.\ {\rm\cite{S2}} Proposition {\rm4.4.10)}. Then we have the following isomorphisms in $D^b(Y_\et,\bZ/p)$:
{\allowdisplaybreaks
\begin{align}
\tag{1}
\displaystyle \bigoplus_{q=0}^a \ \ol\xi{\,}^q \cup - \; &: \; \bigoplus_{q=0}^a \ \lam_{Y,1}^{r-q}[-q] \isom Rg_*\lam_{P,1}^r \\
\tag{2}
\displaystyle \bigoplus_{q=0}^a \ \ol\xi{\,}^q \cup - \; &: \; \bigoplus_{q=0}^a \ \nu_{Y,1}^{r-q}[-q] \isom Rg_*\nu_{P,1}^r \\
\tag{3}
\displaystyle \bigoplus_{q=0}^a \ \ol\xi{\,}^q \cup - \; &: \; \bigoplus_{q=0}^a \ \omega_Y^{r-q}[-q] \isom Rg_*\omega_P^r \\
\tag{4}
\displaystyle \bigoplus_{q=0}^a \ \ol\xi{\,}^q \cup - \; &: \; \bigoplus_{q=0}^a \ \cZ_Y^{r-q}[-q] \isom Rg_*\cZ_P^r.
\end{align}
}Here $\cZ_Y^m$ {\rm(}resp.\ $\cZ_P^m${\rm)} denotes the kernel of $d:\omega_Y^m \to \omega_Y^{m+1}$ {\rm(}resp.\ $d:\omega_P^m \to \omega_P^{m+1}${\rm)}.
\end{lem}
\begin{pf*}{Proof of Lemma \ref{lem3-1}}
Note that $\ol{\xi}$ agrees with the 1st Chern class of the tautological invertible sheaf on $P=\bP(i^*\hspace{-1.8pt}E)$.
Since the problems are \'etale local on $Y$, we may assume that $Y$ is simple.
If $Y$ is smooth, then (1) and (2) are due to Gros \cite{Gr} I Th\'eor\`eme 2.1.11.
The general case is reduced to the smooth case by \cite{S1} Proposition 3.2.1, Corollary 2.2.7.
As for (3), since we have
\[ Rg_*\omega_P^r \simeq \bigoplus_{q=0}^r \ Rg_*(\Omega_{P/Y}^q \otimes_{\cO_Y} \omega_Y^{r-q}) \simeq \bigoplus_{q=0}^r \ (Rg_*\Omega_{P/Y}^q) \otimes_{\cO_Y}^\bL \omega_Y^{r-q} \]
by projection formula, the assertion follows from the isomorphisms
\[ \ol\xi{\,}^q \cup - : \cO_Y[-q] \isom Rg_*\Omega_{P/Y}^q \quad \hbox{($0 \le q \le a$)}. \]
(4) follows from the same arguments as for \cite{Gr} I (2.2.3).
\end{pf*}
We turn to the proof of \eqref{eq3-1} for $r \ge 0$. The case $r=0$ follows from Lemma \ref{lem3-1}\,(1) with $r=0$.
To prove the case $r > 1$, we use the objects $\bK(r-q)_X \in D^b(Y_\et,\bZ/p)$ and $\bK(r)_\bP \in D^b(P_\et,\bZ/p)$ defined in \cite{S2} Lemma 10.4.1, which fit into distinguished triangles
\begin{align*}
 \bK(r-q)_X[-1] \to \mu' \otimes^\bL i^*\fT_1(r-q-1)_X & \to i^*\fT_1(r-q)_X \to \bK(r-q)_X, \\
 \bK(r)_\bP [-1] \lra g^*\mu' \otimes^\bL \gamma^*\fT_1(r-1)_\bP & \lra \gamma^*\fT_1(r)_\bP \lra \bK(r)_\bP.
\end{align*}
Here $\mu'$ denotes the constant sheaf $i^*j_*\mu_p(\simeq \bZ/p)$ on $Y_\et$ and the central arrows are induced by the product structure of Tate twists.
By induction on $r \ge 0$, our task is to show that the morphism
\begin{equation}\label{eq3-2}
\bigoplus_{q=0}^a \  \xi^q \cup - \; : \; \bigoplus_{q=0}^a \ \bK(r-q)_X[-2q] \lra Rg_*\bK(r)_\bP
\end{equation}
is an isomorphism, where we have used the pull-back morphisms
\[ f^* : \bK(r-q)_X \lra \bK(r-q)_\bP \qquad \hbox{($0 \le q \le a$)} \]
induced by the pull-back morphisms for Tate twists (loc.\ cit.\ Proposition 4.2.8, Lemma 2.1.2\,(2)).
By loc.\ cit.\ Lemma 10.4.1\,(2), we have
\[ {\mathscr H}^m(\bK(r-q)_X) \simeq
 \begin{cases} \mu' \otimes \nu_{Y}^{r-q-2} \quad &\hbox{($m = r-q-1$)} \\
 FM_1^{r-q} \quad &\hbox{($m = r-q$)} \\
 0 \quad &\hbox{(otherwise)} \end{cases} \]
and similar facts holds for $\bK(r)_P$ (see \S\ref{sect2} for $FM_n^{q}$). Therefore \eqref{eq3-2} is an isomorphism by Lemma \ref{lem3-1} and Theorems \ref{thm2-1}\,(3),\,(4) and \ref{thm2-2} with $D=\emptyset$ (see also the projection formula in loc.\ cit.\ 4.4.10). This completes the proof of Theorem \ref{thm3-1}.
\end{pf}
The following corollary \ref{cor3-1} follows immediately from Theorems \ref{prop2-2}, \ref{thm3-1} and the projection formulra (\cite{S2} Corollary 7.2.4), which is a $p$-adic version of homotopy invariance and plays an important role in our construction of cycle class maps (see \S\ref{sect5} below).
\begin{cor}\label{cor3-1}
Let the notation be as in Theorem {\rm\ref{thm3-1}}. Let $E' \subset E$ be a subbundle of rank $a$, and let $\bP'$ be the associated projective bundle. Let $\varphi : \bP' \to \bP$ be the natural closed immersion, and assume that the inverse image $\varphi^*(\cO(1)_E)$ is isomorphic to the tautological invertible sheaf of $\bP'$. Then the composite morphism
\[  \fT_n(r)_X \os{f^*}\lra Rf_*\fT_n(r)_{\bP} \lra Rf_*\fT_n(r)_{(\bP,\bP')} \]
is an isomorphism in $D^b(X_\et,\bZ/p^n)$.
\end{cor}

\medskip

\section{\bf Chern class}\label{sect3'}
The main aim of this section is to construct the Chern class map \eqref{eq3'-3} below. Let $S$ and $p$ be as in Setting \ref{set2-1}, and let $\ccS$ be the category whose objects are $S$-schemes satisfying the conditions in Setting \ref{set2-1} for $X$ and whose morphisms are $S$-morphisms. Let $X_\star$ be a simplicial object in $\ccS$, i.e., a contravariant functor
\[ X_\star : \vD^\op \lra \ccS, \]
where $\vD$ denotes the simplex category. For a morphism $\gamma : [a] \to [b]$ in $\vD$, we often write 
\[ \gamma^X : X_b \lra X_a \qquad \hbox{($X_a:=X_\star([a])$)} \]
for $X_\star(\gamma)$, which is a morphism in $\ccS$. For integers $0 \le i \le a$, let $d^i$ be the coface map in $\vD$:
\[ d^i : [a] \lra [a+1],\qquad j \mapsto \begin{cases} \;\;\, j \quad & (0 \le j < i) \\ j+1 \quad & (i \le j \le a). \end{cases} \]
For integers $0 \le i \le a$, we often write
\[ d_i : X_{a+1} \to X_a \]
for $(d^i)^X$. See \cite{Fr} \S1 for the definition of the small \'etale site $(X_\star)_\et$ on $X_\star$.
\begin{defn}\label{def3'-1}
\begin{enumerate}
\item[(1)]
We define a complex $C_n(r)_{X_\star}^\bullet$ of sheaves on $(X_\star)_\et$ by restricting the complex $C_n(r)^\bullet$ on $\ccS_\et$, cf.\ Definition {\rm \ref{def2-1}}, Corollary {\rm \ref{cor2-1}}. We write $\fT_n(r)_{X_\star}$ for the complex $C_n(r)_{X_\star}^\bullet$ regarded as an object of $D^b((X_\star)_\et,\zpn)$.
\item[(2)]
We define a canonical morphism
\[ \varrho : \Gm[-1] \lra \fT_n(1)_{X_\star} \quad \hbox{ in } \;\; D^b((X_\star)_\et) \]
by the composite morphism
\[ \Gm[-1] \lra \Gm \otimes^\bL \zpn [-1] \isom \fT_n(1)_{X_\star}, \]
where the left arrow denotes the canonical morphism induced by $\bZ \to \zpn$ and the right arrow is the restriction of the isomorphism in Remark {\rm\ref{rem2-3}\,(2)}.
\end{enumerate}
\end{defn}
\noindent
We next review the following basic notions:
\begin{defn}
\begin{enumerate}
\item[(1)]
{\it A vector bundle} over $X_\star$ is a morphism $f : E_\star \to X_\star$ of simplicial schemes such that $f_a : E_a \to X_a$ is a vector bundle for any $a \ge 0$ and such that the commutative diagram
\begin{equation}\label{eq3'-a}
\xymatrix{ E_b \ar[r]^{f_b} \, \ar[d]_{\gamma^E} & X_b \ar[d]^{\gamma^X} \\
 E_a \ar[r]^{f_a} \, & X_a}
\end{equation}
induces an isomorphism $E_b \cong \gamma^{X*}\!E_a:=E_a \times_{X_a} X_b$ of vector bundles over $X_b$ for any morphism $\gamma : [a] \to [b]$ in $\vD$\; {\rm(cf.\ \cite{Gi2} {\it Example} 1.1)}.
\item[(2)]
We say that a morphism $f : X_\star \to Y_\star$ of simplicial objects of $\ccS$ is {\it a regular closed immersion} if $f_a : X_a \to Y_a$ is a regular closed immersion for any $a \ge 0$ and if the diagram
\begin{equation}\label{eq3'-b}
\xymatrix{ X_b \ar[r]^{f_b} \, \ar[d]_{\gamma^X} & Y_b \ar[d]^{\gamma^Y} \\
 X_a \ar[r]^{f_a} \, & Y_a}
\end{equation}
is cartesian for any morphism $\gamma : [a] \to [b]$ in $\vD$. {\it An effective Cartier divisor $X_\star$ on $Y_\star$} is a regular closed immersion $X_\star \to Y_\star$ of pure codimension $1$.
\end{enumerate}
\end{defn}
\noindent
We now define the first Chern classes of effective Cartier divisors and line bundles.
\begin{defn}\label{def3'-3}
\begin{enumerate}
\item[{\rm(1)}]
For an effective Cartier divisor $D_\star$ on $X_\star$, we define the first Chern class $\sfc_1(D_\star) \in H^2_{D_\star}(X_\star,\fT_n(1))$ as the value of the first Chern class $\sfc_1(D_\star) \in H^1_{D_\star}((X_\star)_\zar,\cO^\times)$ under the composite map
\[ 
H^1_{D_\star}((X_\star)_\zar,\cO^\times) \os{\ep^*}\to H^1_{D_\star}(X_\star,\Gm) \os{\varrho}\to H^2_{D_\star}(X_\star,\fT_n(1)), \]
where
$\ep: (X_\star)_\et \to (X_\star)_\zar$ denotes the continuous map of small sites. The arrow $\varrho$ denotes that in Definition {\rm\ref{def3'-1}\,(2)}.
\item[{\rm(2)}]
For a line bundle $L_\star$ on $X_\star$, we define the first Chern class
\[ \sfc_1(L_\star) \in H^2(X_\star,\fT_n(1)) \]
as the value of the isomorphism class $[L_\star] \in H^1((X_\star)_\zar,\cO^\times)$ {\rm({\it cf.}\ \cite{Gi2} {\it Example} 1.1)} under the composite map
\begin{equation}\label{eq3'-1}
 H^1((X_\star)_\zar,\cO^\times) \os{\ep^*}\lra H^1(X_\star,\Gm) \os{\varrho}\lra H^2(X_\star,\fT_n(1)).
\end{equation}
\end{enumerate}
\end{defn}
\par\medskip
The following proposition plays a key role in the proof of the Whitney sum formula in Proposition \ref{prop3-2}\,(3) below.
\begin{prop}\label{prop3'-3}
Let $f : X_\star \hra X'_\star$ be a regular closed immersion of simplicial objects in $\ccS$ of pure codimension $c \ge 1$. Assume that the given morphisms $d_0,d_1 : X'_1 \to X'_0$ are smooth. Then there exists a Gysin morphism
\[ \gys_f : \fT_n(r)_{X_\star} \lra Rf^!\fT_n(r+c)_{X'_\star}[2c] \quad \hbox{ in } \;\; D^+((X_\star)_\et,\zpn) \]
satisfying the following three properties{\rm:}
\begin{enumerate}
\item[(a)]
{\rm ({\it Consistency with the first Chern class})}\;
If $r=1$, then the value of $1 \in \zpn = H^0(X_\star,\fT_n(0))$ under the Gysin map
\[ \gys_f : H^0(X_\star,\fT_n(0)) \lra H^2_X(X'_\star,\fT_n(1)) \]
agrees with the first Chern class $\sfc_1(X_\star)$ in Definition {\rm\ref{def3'-3}\,(1)}.
\item[(b)]
{\rm ({\it Transitivity})}\;
For another regular closed immersion $g : X'_\star \hra X''_\star$ of simplicial objects in $\ccS$ of pure codimension $c' \ge 1$ with $d_0,d_1 : X''_1 \to X''_0$ smooth, the composite morphism
\[ \begin{CD} \fT_n(r)_{X_\star} \os{\gys_f}\lra Rf^!\fT_n(r+c)_{X'_\star}[2c] @>{Rf^!(\gys_g)}>> Rf^!Rg^!\fT_n(r+c+c')_{X''_\star} [2c+2c'] \\ @= R(g \circ f)^!\fT_n(r+c+c')_{X''_\star}[2(c+c')]
 \end{CD}\]
agrees with $\gys_{g \circ f}$\,.
\item[(c)]
{\rm ({\it Projection formula})}\;
The following diagram commutes in $D((X_\star)_\et)${\rm:}
\[ \xymatrix{
Rf^!\fT_n(q)_{X'_\star} \otimes^{\bL} \fT_n(r)_{X_\star} \ar[r]^-{\id \otimes \gys_f} \ar[d]_\pi & Rf^!\fT_n(q)_{X'_\star} \otimes^{\bL} Rf^!\fT_n(n+r)_{X'_\star}[2r] \ar[d]^{\text{product}} \\
\fT_n(q+r)_{X_\star} \ar[r]^-{\gys_f} & Rf^!\fT_n(q+r+c)_{X'_\star}[2c]\,,
} \]
where the left vertical arrow $\pi$ is the composite morphism
\[ \begin{CD}
Rf^!\fT_n(q)_{X'_\star} \otimes^{\bL} \fT_n(r)_{X_\star} @>{f^* \otimes \id}>> \fT_n(q)_{X_\star} \otimes^{\bL} \fT_n(r)_{X_\star} @>{\text{product}}>> \fT_n(q+r)_{X_\star},
\end{CD}\]
where the products mean the restriction of the product structure in Remark {\rm\ref{rem2-3}}.
\end{enumerate}
\end{prop}
\begin{pf}
Put $U_\star:=X_\star \otimes \bZ[p^{-1}]$ and $V_\star:=X'_\star \otimes \bZ[p^{-1}]$. Let $\varphi : U_\star \hra V_\star$ be the regular closed immersion induced by $f$. By the absolute purity \cite{FG} and the spectral sequence
\[ E_1^{a,b}=H^b_{U_a}(V_a,\mu_{p^n}^{\otimes c}) \Lra H^{a+b}_{U_\star}(V_\star,\mu_{p^n}^{\otimes c}), \]
we have
\[ H^{2c}_{U_\star}(V_\star,\mu_{p^n}^{\otimes c}) \simeq \ker(d_0^*-d_1^* : H^{2c}_{U_0}(V_0,\mu_{p^n}^{\otimes c}) \lra H^{2c}_{U_1}(V_1,\mu_{p^n}^{\otimes c})). \]
By the smoothness assumption on $d_0,d_1 : V_1 \to V_0$, the cycle class $\cl_{V_0}(U_0) \in H^{2c}_{U_0}(V_0,\mu_{p^n}^{\otimes c})$ lies in the group on the right hand side, loc.\ cit.\ Proposition 1.1.3. We thus define the cycle class
\[ \cl_{V_\star}(U_\star) \in H^{2c}_{U_\star}(V_\star,\mu_{p^n}^{\otimes c}) \]
as the element corresponding to $\cl_{V_0}(U_0)$. Since $\varphi^*\mu_{p^n,V_\star}^{\otimes r} \simeq \mu_{p^n,U_\star}^{\otimes r}$ on $(U_\star)_\et$\,, the cup product with $\cl_{V_\star}(U_\star)$ defines a Gysin morphism
\[\begin{CD} \gys_\varphi :  \mu_{p^n,U_\star}^{\otimes r} \simeq \varphi^*\mu_{p^n,V_\star}^{\otimes r} @>{\cl_{V_\star}(U_\star) \,\cup \, -}>> R\varphi^! \mu_{p^n,V_\star}^{\otimes r+c}[2c]  \quad \hbox{ in } \;\; D^+((U_\star)_\et,\zpn), \end{CD}\]
which satisfies the three properties (a)\,--\,(c) listed above (see loc.\ cit.\ Proposition 1.2.1 for (b)). We show that there exists a unique morphism
\[ \gys_f : \fT_n(r)_{X_\star} \lra Rf^!\fT_n(r+c)_{X'_\star}[2c] \quad \hbox{ in } \;\; D^+((X_\star)_\et,\zpn) \]
that extends $\gys_\varphi$. Put
\[ Y_\star:=X_\star \otimes \bZ/p\bZ, \quad \fL := \fT_n(r)_{X_\star} \;\; \hbox{ and }  \;\; \fM := Rf^!\fT_n(r+c)_{X'_\star}[2c] . \]
and let $\alpha$ and $\beta$ be as follows:
\[\xymatrix{ U_\star \; \ar@<-1pt>@{^{(}->}[r]^\beta & X_\star & \ar@<1pt>@{_{(}->}[l]_\alpha \; Y_\star\,. }\]
Consider the following diagram in $D^+((X_\star)_\et,\zpn)$ whose lower row is distinguished:
\begin{equation}\label{eq3'-4}
\xymatrix{ & \fL \ar[r]^-t  & \tau_{\leq r}R\beta_*\mu_{p^n}^{\otimes r} \ar[d]_{R\beta_*(\gys_\varphi)} \\
\alpha_*R\alpha^! \fM \ar[r]^-{\alpha_*} & \fM \ar[r]^-{\beta^*} & R\beta_*\beta^*\fM \ar[r]^-{-\delta} & \alpha_*R\alpha^! \fM[1]. }
\end{equation}
Here the upper horizontal arrow is the canonical morphism (cf.\ Propositions \ref{prop2-0} and \ref{prop2-1}), and the lower row is the localization distinguished triangle for $\fM$ (cf.\ \cite{S2} (1.9.2)). We have
\begin{equation}\label{eq3'-5}
 \tau_{\le r} \alpha_*R\alpha^! \fM = 0
\end{equation}
by the purity in loc.\ cit.\ Theorem 4.4.7, which implies that
\[ \Hom_{D^+((X_\star)_\et,\zpn)}(\tau_{\leq r}R\beta_*\mu_{p^n}^{\otimes r}, \alpha_*R\alpha^! \fM[1]) = \Hom_{S((X_\star)_\et)}(R^r\beta_*\mu_{p^n}^{\otimes r}, \alpha_*R^{r+1}\alpha^! \fM). \]
By this fact and the compatibility fact in loc.\ cit.\ Theorem 6.1.1, one can easily check that the composite $(-\delta) \circ R\beta_*(\gys_\varphi)$ is zero in $D^+((X_\star)_\et,\zpn)$. Therefore we obtain a unique morphism $\gys_f$ that extends $\gys_\varphi$ again by \eqref{eq3'-5} and by loc.\ cit.\ Lemma 2.1.2\,(1). The property (a) of $\gys_f$ is straight-forward, and the property (b) follows from the uniqueness of $\gys_f$. The property (c) follows from the same argument as for loc.\ cit.\ Corollary 7.2.4.
\end{pf}
\par\medskip
Following the method of Grothendieck \cite{G} and Gillet \cite{Gi}, we define the Chern classes
\[ \sfc(E_\star)=(\sfc_i(E_\star))_{i \ge 0} \in \bigoplus_{i \ge 0} \ H^{2i}(X_\star,\fT(i)) \]
of a vector bundle $E_\star$ over $X_\star$ as follows. Let $E_\star$ be of rank $a$, and let $f$ be the natural projection $\bP(E_\star) \to X_\star$. Let $\xi \in H^2(\bP(E_\star),\fT_n(1))$ be the value of the first Chern class of the tautological line bundle on $\bP(E_\star)$, cf.\ Definition \ref{def3'-3}\,(2). Noting the Dold-Thom isomorphism
\begin{equation}\label{eq3'-2}
 \bigoplus_{i=1}^a \ H^{2i}(X_\star,\fT_n(i)) \simeq H^{2a}(\bP(E_\star),\fT_n(a)), \quad (b_i)_{i=1}^a \mapsto \sum_{i=1}^a \ f^*(b_i) \cup \xi^{a-i}
\end{equation}
obtained from Theorem \ref{thm3-1} and \cite{Gi} Lemma 2.4, we define
\[ \sfc_0(E_\star):=1 \quad \hbox{ and } \quad \sfc_i(E_\star):=0 \;\; \hbox{ for \; $i > a$}, \]
and define $\sfc_i=\sfc_i(E_\star)$ for $i=1,2,\dotsc,a$ by the equation
\[ \xi^a + f^*(\sfc_1) \cup \xi^{a-1} + \dotsb + f^*(\sfc_{a-1})\cup \xi + f^*(\sfc_a) = 0 \]
in $H^{2a}(\bP(E_\star),\fT_n(a))$.
\begin{prop}\label{prop3-2}
Let $X_\star$ be a simplicial object in $\ccS$ with $d_0,d_1 : X_1 \to X_0$ smooth. Then the Chern classes $\sfc(E_\star)$ of vector bundles $E_\star$ over $X_\star$ satisfy the following properties{\rm:}
\begin{enumerate}
\item[(1)] {\rm(}Normalization{\rm)}\; If $E_\star$ is a line bundle, then we have $\sfc_0(E_\star)=1$ and $\sfc_i(E_\star)=0$ for $i > 1$, and $\sfc_1(E_\star)$ is the value of the first Chern class $\sfc_1(E_\star) \in H^1((X_\star)_\zar,\cO^\times)$ under the map \eqref{eq3'-1}.
\item[(2)] {\rm(}Functoriality{\rm)}\; For a morphism $f : X_\star \to X'_\star$ of simplicial objects in $\ccS$ and a vector bundle $E_\star$ on $X'_\star$, we have \[ \sfc(f^*\! E_\star)=f^*\sfc(E_\star), \]
where $f^*$ on the right hand side denotes the pull-back map obtained from the contravariant functoriality of $\fT_n(r)$, cf.\ Remark {\rm\ref{rem2-2}}.
\item[(3)] {\rm(}Whitney sum{\rm)}\; For a short exact sequence $0 \to E'_\star \to E_\star \to E''_\star \to 0$ of vector bundles on $X_\star$, we have \[ \sfc_r(E_\star)=\sum_{s+t=r} \sfc_s(E'_\star) \cup \sfc_t(E''_\star) \in H^{2r}(X_\star,\fT_n(r)) \]
for each $r \ge 0$, where the cup product is taken with respect to the product structure of $\{ \fT_n(r)_{X_\star} \}_{r \ge 0}$ in the derived category of \'etale sheaves on $X_\star$\,, cf.\ Remark {\rm\ref{rem2-3}\,(1)}.
\end{enumerate}
Moreover, the Chern classes $\sfc(E_\star)$ are characterized by these three properties.
\end{prop}
\begin{pf}
The properties (1) and (2) immediately follow from this definition of Chern classes. The last assertion on the uniqueness follows from the splitting principle of vector bundles, whose details are straight-forward and left to the reader.\par
We prove the property (3) using the arguments of Grothendieck in \cite{G} p.\ 144 Theorem 1\,(iii), as follows. Let $\pi': D'_\star \to X_\star$ and $\pi'' : D''_\star \to X_\star$ be the (simplicial) flag schemes of $E'_\star$ and $E''_\star$, respectively, and put
\[ D_\star:=D'_\star \times_{X_\star} D''_\star \,, \]
which is identified with the flag scheme of $\pi''^*E'_\star$ over $D''_\star$. Let $f : D_\star \to X_\star$ be the natural projection. Since the pull-back map
\[ f^* : H^{2i}(X_\star,\fT_n(i)) \lra H^{2i}(D_\star,\fT_n(i)), \]
is injective by \eqref{eq3'-2}, we may replace $(X_\star,E_\star,E'_\star,E''_\star)$ with $(D_\star,f^*\!E_\star,f^*\!E'_\star,f^*\!E''_\star)$ to assume that $E_\star$ has a filtration by subbundles
\[ E_\star=E^0_\star \supset E^1_\star \supset \dotsb \supset E^a_\star = 0 \qquad \hbox{($a:=\rank(E)$)} \]
such that the quotient $E^i_\star/E^{i+1}_\star$ is a line bundle over $X_\star$ for $0 \le i \le a-1$ and such that $E_\star^b= E'_\star$ for $b=\rank(E'')$. Now let
\[ g :  X'_\star := \bP(E_\star) \lra X_\star \]
be the projective bundle associated with $E_\star$.
Let $L_\star$ be the tautological line bundle over $X'=\bP(E_\star)$ and let $s : X'_\star \to g^*\!E_\star \otimes L_\star=:F_\star$ be the section induced by the canonical inclusion $(L_\star)^\vee \hra g^*\!E_\star$:
\[ s : X'_\star \os{1}\lra \bA^1_{X'_\star} \lra g^*\!E_\star \otimes L_\star=F_\star \,. \]
Put $F^i_\star:=g^*\!E^i_\star \otimes L$ and $V^i_\star:=s^{-1}(F^i_\star)$ for $0 \le i \le a$. Then $V^a_\star$ is empty because $s$ does not vanish. Moreover $V^i_\star$ is a simplicial object in $\ccS$ and the section
\[ s_i : V^i_\star \os{s}\lra F^i_\star|_{V^i_\star} \lra (F^i_\star/F^{i+1}_\star)|_{V^i_\star} \]
meets the zero section transversally for $0 \le i \le a-1$, cf.\ \cite{G} p.\ 147. On the other hand, to prove the Whitney sum formula, it is enough to show
\[  \prod_{i=0}^{m+n-1} \ \sfc_1(F^i_\star/F^{i+1}_\star) = 0. \]
By these facts, we are reduced to the following simplicial analogue of loc.\ cit.\ p.141 Lemma 2:
\begin{lem}\label{lem2-1}
Let $X_\star$ be a simplicial object in $\ccS$ with $d_0,d_1:X_1 \to X_0$ smooth. Let $E_\star$ be a vector bundle of rank $a$ over $X_\star$. Let $(E^i_\star)_{0 \le i \le a}$ be a descending filtration on $E_\star$ such that $E^i_\star$ is a subbundle of rank $a-i$ with $E^0_\star=E_\star$ and such that the quotient line bundles $E^i_\star/E^{i+1}_\star$ exist. For $1 \le i \le a$, put
\[ \xi_i:=\sfc_1(E^{i-1}_\star/E^i_\star) \in H^2(X_\star,\fT_n(1)). \]
Let $s : X_\star \to E_\star$ be a section of $E_\star \to X_\star$\,. For $0 \le i \le a$, put
\[ V^i_\star : = s^{-1}(E^i_\star), \qquad L^i_\star:=(E^i_\star/E^{i+1}_\star)|_{V^i_\star} \quad \hbox{{\rm(}restriction of $E^i_\star/E^{i+1}_\star$ onto $V^i_\star${\rm)}}, \]
and let $s_i : V^i_\star \to L^i_\star$ be the section induced by $s$. Assume the following condition{\rm:}
\begin{enumerate}
\item[$\bullet$]
$V^i_\star$ is a simplicial object in $\ccS$ for $0 \le i \le a$, and $s_i$ intersects the zero section transversally for $0 \le i \le a-1$.
\end{enumerate}
Then we have
\[ \cl_{X_\star}(V^a_\star) = \prod_{i=1}^a \ \xi_i \, \quad \hbox{ in } \quad H^{2a}(X_\star,\fT_n(a))\,. \]
Here $\cl_{X_\star}(V^a_\star)$ denotes the value of $1$ under the Gysin map
\[ \zpn=H^0(V^a_\star,\fT_n(0)) \lra H^{2a}(X_\star,\fT_n(a)) \]
and the product on the right hand side means the cup product with respect to the product structure on $\fT_n(*)_{X_\star}$\,, cf.\ Remark {\rm\ref{rem2-3}\,(1)}.
\end{lem}
\noindent
One can easily check this lemma by the properties of the Gysin morphisms in Proposition \ref{prop3'-3} and the arguments in loc.\ cit.\ p.141 Lemma 2. This completes the proof of Proposition \ref{prop3-2}.
\end{pf}

Now let $X$ be a scheme which belongs to $\ccS$. Applying the construction of Chern classes to the case $X_\star=\BGL_r/X$ and $E_\star=$ universal rank $r$ bundle over $\BGL_r/X$, we obtain a Chern class
\[ \sfc_r(E_\star) \in H^{2r}(\BGL_r/X,\fT_n(r)), \]
which is called the {\it universal rank $r$ Chern class}. On the other hand, let $I_n(r)_X^\bullet$ be an injective resolution of the complex $C_n(r)_X^\bullet$ on $X_\et$ defined in Definition \ref{def2-1}, and consider the following complex of abelian sheaves on $X_\zar$:
\[\dotsb \lra \ep_*I_n(r)_X^{2r-2} \os{d^{2r-2}}\lra \ep_*I_n(r)_X^{2r-1} \os{d^{2r-1}}\lra \ker(d^{2r}: \ep_*I_n(r)_X^{2r} \to \ep_*I_n(r)_X^{2r+1}). \]
Here $\ep : X_\et \to X_\zar$ denotes the continuous map of small sites, and we regarded this sequence as a chain complex with the most right term placed in degree $0$. We apply the Dold-Puppe construction \cite{DP} to this complex to obtain a sheaf of simplicial abelian groups, which we denote by $\sfK(\fT_n(r),2r)$. For a closed subset $Z$ of $X$ and a non-negative integer $i \ge 0$, we define the Chern class map
\begin{equation}\label{eq3'-3}
 \sfc_{r,i}^Z : K_i^Z(X) \lra H^{2r-i}_Z(X,\fT_n(r))
\end{equation}
as the following composite map (cf.\ \cite{Gi} Definition 2.22):
{\allowdisplaybreaks
\begin{align*}
K^Z_i\!(X) \
 & \lra \ H_Z^{-i}(X_\zar,\bZ \times \bZ_\infty \BGL(\cO_X)) \\
 & \os{\pr}\lra \ H_Z^{-i}(X_\zar,\bZ_\infty \BGL(\cO_X)) \\
 & \hspace{-9pt} \os{\pi_i(\bZ_\infty\sfc_r)}\lra H_Z^{-i}(X_\zar,\bZ_\infty \sfK(\fT_n(r),2r)) \phantom{\Big|} \cong H^{2r-i}_Z(X_\et,\fT_n(r)).
\end{align*}
}Here $\bZ_\infty$ denotes the Bousfield-Kan completion \cite{BoK}, and we have used the universal rank $r$ Chern class $\sfc_r(E_\star)$ to define the arrow $\pi_i(\bZ_\infty\sfc_r)$. See \cite{Gi} Proposition 2.15 for the first arrow and see loc.\ cit.\ p.\ 226 for the last isomorphism. The map $\sfc_{r,0}^X$ agrees with $\sfc_r$ for $X_\star=X$ (constant simplicial scheme) defined before.
\begin{thm}\label{thm3-2}
\begin{enumerate}
\item[(1)]
$\sfc_{r,i}^Z$ is contravariantly functorial in the pair $(X,Z)$, that is, for a morphism $f :X \to X'$ in $\ccS$ and a closed subset $Z' \subset X'$ with $f^{-1}(Z') \subset Z$, there is a commutative diagram
\[\xymatrix{ K_i^{Z'}\!(X') \ar[r]^-{\sfc_{r,i}^{Z'}} \ar[d]_{f^*} & H^{2r-i}_{Z'}(X',\fT_n(r)) \ar[d]^{f^*} \\ K_i^Z(X) \ar[r]^-{\sfc_{r,i}^Z} & H^{2r-i}_Z(X,\fT_n(r)). }\]
\item[(2)]
$\sfc_{r,i}^Z$ is additive for $i >0$.
\item[(3)]
The induced Chern character
\[ \ch : \bigoplus_{i \ge 0} \ K_i(X) \lra \prod_{i,r \ge 0} \ H^i(X,\fT_{\qp}(r)) \]
with $H^i(X,\fT_{\qp}(r)):=\qp \otimes_{\zp} \varprojlim{}_{n \ge 1} \ H^i(X,\fT_n(r))$ is a ring homomorphism.
\end{enumerate}
\end{thm}
\begin{pf}
(1) follows from the functoriality results in Proposition \ref{prop3-2}\,(2) and Proposition \ref{prop2-1}. The assertion (2) follows from Proposition \ref{prop3-2}\,(3) and \cite{Gi} Lemma 2.26. See loc.\ cit.\ Definition 2.34 and Proposition 2.35 for the assertion (3).
\end{pf}
\medskip

\section{\bf Purity along log poles}\label{sect4}
\medskip
This and the next section are devoted to the construction of cycle class maps from higher Chow groups to $p$-adic \'etale Tate twists.
Let $k$ be a field of characteristic $p>0$, and let $Y$ be a normal crossing variety over $k$. Let $c$ be a positive integer.
Let $D$ be a non-empty admissible divisor on $Y$, and let $Z$ be a reduced closed subscheme of $Y$ which has codimension $\ge c$
  and contained in $D$. Let $i : Z \hra Y$ be the natural closed immersion.
In this section, we first prove the following purity result:
\begin{thm}[{\bf Purity}]\label{thm4-1}
We have $R^q i^!\nu_{(Y,D),n}^r = 0$ for $q \le c$.
\end{thm}
\noindent
See Theorem \ref{cor4-1} below for a consequence of this theorem. Note that the assertion for $q=c$ does not follow directly from Gros' purity \cite{Gr} II Th\'eor\`eme 3.5.8, even when $Y$ is smooth. The case that $Y$ is smooth has been considered and proved independently by Mieda (\cite{Mi1} Proof of Theorem 2.4, \cite{Mi2}). Although the proof of Theorem \ref{thm4-1} given below is essentially a variant of Gros' proof of \cite{Gr} II Th\'eor\`eme 3.5.8, we include detailed computations for the convenience of the reader.
\par\medskip
\begin{pf*}{Proof of Theorem \ref{thm4-1}}
Since the problem is \'etale local, we may assume that $Y$ is simple and that there exists a pair $(\cY,\cD)$ satisfying the following (1) and (2):
\begin{enumerate}
\item[(1)] $\cY$ is smooth over $k$ and contains $Y$ as a normal crossing divisor,
\item[(2)] $\cD$ is a normal crossing divisor on $\cY$ such that $Y \cup \cD$ has normal crossings on $\cY$ and such that $Y \cap \cD=D$.
\end{enumerate}
Then we have a short exact sequence on $\cY_\et$
\[ 0 \lra \logwitt {(\cY,\cD)} n {r+1} \lra \logwitt {(\cY,Y \cup \cD)} n {r+1} \lra \nu_{(Y,D),n}^r \lra 0 \]
(a variant of \eqref{eq1-1}), which reduces the assertion to the case that $Y$ is smooth over $k$ and that $D$ has normal crossings on $Y$. Furthermore, we may assume that $Z$ is {\it regular and of pure codimension} $c$ by a standard devissage argument. In what follows, we prove
\begin{equation}\label{eq4-1}
R^q i^!\logwitt {(Y,D)} n r = 0 \quad \hbox{ for } \quad q \ne c+1,
\end{equation}
assuming that $Y$ is smooth over $k$ (and that $Z$ is regular and contained in $D$).
By the short exact sequences on $Y_\et$
\begin{align*}
 0 \lra \logwitt {(Y,D)} {n-1} r \lra \logwitt {(Y,D)} n r \lra \Omega_{(Y,D),\log}^r \lra 0, \\
 0 \lra \Omega_{(Y,D),\log}^r \lra \cZ_Y^r(\log D) \os{1-C}{-\hspace{-7pt}\lra} \Omega_Y^r(\log D) \lra 0 \quad
\end{align*}
(cf.\ \eqref{eq1-1} and \cite{S1} Corollary 2.2.5\,(2), Lemma 2.4.6), it is enough to show that
\begin{align}
\label{eq4-2} R^q i^!\cZ_Y^r(\log D) = 0 = R^q i^!\Omega_Y^r(\log D) \quad \hbox{ for } \quad q \ne c, \quad \\
\label{eq4-3} 1-C : R^c i^!\cZ_Y^r(\log D) \lra R^c i^!\Omega_Y^r (\log D) \quad \hbox{ is injective}.
\end{align}
The assertion \eqref{eq4-2} follows from the smoothness of $Y$ and the fact that $\Omega_Y^r(\log D)$ and $\cZ_Y^r(\log D)$ are locally free over $\cO_Y$ and $(\cO_Y)^p$, respectively. \par
We prove \eqref{eq4-3} in what follows. Since the problem is \'etale local on $Z$, we may assume the following condition:
\begin{enumerate}
\item[($\star$)] {\it $Y$ is affine, and there exists a regular sequence $t_1,\dotsb,t_c \in \vG(Y,\cO_Y)$ such that the ideal $(t_1,\dotsb,t_c) \subset \cO_Y$ defines $Z$ and such that $t_1,\dotsb,t_a$\,{\rm(}$1 \le a \le c${\rm)} are uniformizers of the irreducible components of $D$. Moreover, $\Omega^1_Y$ is free over $\cO_Y$ and has a basis $\{ dt_\lam \}_{\lam \in \vL}$ which contains $dt_1,\dotsc,dt_c$.}
\end{enumerate}
Let $\tau$ be the natural open immersion
\[\xymatrix{ \tau : Y[t_1^{-1},\dotsc,t_c^{-1}] \; \ar@{^{(}->}[r] &Y.} \]
For $j=1,\dotsc,c$, let $\sigma_j$ be the natural open immersion
\[\xymatrix{ \sigma_j : Y[t_1^{-1},\dotsc,t_{j-1}^{-1},t_{j+1}^{-1},\dotsc,t_c^{-1}] \; \ar@{^{(}->}[r]  &Y. }\]
We recall here the following standard fact (cf.\ \cite{Gr} II (3.3.6)):
\begin{lem}\label{lem4-1}
We have isomorphisms of sheaves on $Z_\et$
\begin{align}
\label{eq-4.2.1}
 R^c i^!\Omega_Y^r (\log D)
& \simeq \tau_*\tau^*\Omega_Y^r(\log D) \Big/\textstyle \sum_{j=1}^c\ \sigma_{j*}\sigma_j^*\Omega_Y^r (\log D) \\
\notag
& = \tau_*\tau^*\Omega_Y^r \Big/\textstyle \Big(\sum_{j=1}^a \ \sigma_{j*}\sigma_j^*\Omega_Y^r (\log D)
  +\sum_{j=a+1}^c \ \sigma_{j*}\sigma_j^*\Omega_Y^r\Big), \\
\label{eq-4.2.2}
R^c i^!\cZ_Y^r (\log D)
& \simeq \tau_*\tau^*\cZ_Y^r(\log D) \Big/\textstyle \sum_{j=1}^c\ \sigma_{j*}\sigma_j^*\cZ_Y^r (\log D) \\
\notag
& = \tau_*\tau^*\cZ_Y^r \Big/\textstyle \Big(\sum_{j=1}^a \ \sigma_{j*}\sigma_j^*\cZ_Y^r (\log D)
  +\sum_{j=a+1}^c \ \sigma_{j*}\sigma_j^*\cZ_Y^r\Big),
\end{align}
where we regarded the sheaves on the right hand side as sheaves on $Z_\et$ naturally.
\end{lem}
We define an ascending filtration $\Fil_m$\,($m\ge 0$) on $\tau_*\tau^*\Omega_Y^r$ as
\[ \Fil_m (\tau_*\tau^*\Omega_Y^r):= \bigg\{ \frac{1}{(t_1t_2\dotsb t_c)^{pm}}\, \omega \in \tau_*\tau^*\Omega_Y^r \ \bigg| \ \omega \in \Omega_Y^r \bigg\}.\]
Let $\Fil_\bullet(R^c i^!\Omega_Y^r (\log D))$ be the induced filtration, and let $\Fil_\bullet(R^c i^!\cZ_Y^r (\log D))$ be its inverse image under the canonical map
\begin{equation}\label{eq4-4}
R^c i^!\cZ_Y^r (\log D) \lra R^c i^!\Omega_Y^r (\log D)) \end{equation}
induced by the natural inclusion $\cZ_Y^r (\log D) \hra \Omega_Y^r (\log D))$. Note that
\[ 
 R^c i^!\cZ_Y^r (\log D) = \bigcup_{m=0}^\infty \  \Fil_m(R^c i^!\cZ_Y^r (\log D)). \]
We prove here the following lemma.
\begin{lem}\label{lem4-2} \begin{enumerate}
\item[(1)] The map \eqref{eq4-4} is injective.
\item[(2)] $\Fil_m(R^c i^!\cZ_Y^r (\log D))$ is generated by elements of the form
\[ \bigg[ \frac{1}{(t_1t_2\dotsb t_c)^{pm}}\, \omega \bigg] \quad \hbox{ with } \;\; \omega \in \cZ^r_Y, \]
where $[-]$ denotes the residue class in $\Fil_m(R^c i^!\cZ_Y^r (\log D))$ via \eqref{eq-4.2.2}.
\item[(3)] The kernel of the projection via \eqref{eq-4.2.1}
\[ \Fil_m(\tau_*\tau^*\Omega_Y^r) \lra \Fil_m(R^c i^!\Omega_Y^r (\log D)) \]
agrees with the subgroup
\[  \Bigg\{ \sum_{j=1}^c \ \frac{1}{(t_1\dotsb t_{j-1}t_{j+1}\dotsb t_c)^{pm}}\, \omega_j \ \Bigg|
 \ \begin{array}{ll} \omega_j \in \Omega_Y^r(\log D_j) & \hbox{{\rm(}$1 \le j \le a${\rm)}} \\
 \omega_j \in \Omega_Y^r \phantom{\Big(} & \hbox{{\rm(}$a < j \le c${\rm)}} \end{array} \Bigg\}, \]
where $D_j \subset Y$ denotes the regular divisor defined by $t_j$ for $1 \le j \le a$.
\end{enumerate} \end{lem}
\begin{pf} We prove (1). By the short exact sequence
\[ 0 \lra \cZ_Y^r (\log D) \lra \Omega_Y^r (\log D) \os{d}\lra  \cB_Y^{r+1}(\log D) \lra 0, \]
we have a long exact sequence on $Z_\et$
\[ \dotsb \lra R^{c-1}i^! \cB_Y^{r+1}(\log D) \lra R^c i^!\cZ_Y^r (\log D) \lra R^c i^!\Omega_Y^r (\log D) \lra \dotsb. \]
The assertion follows from the smoothness of $Y$ and the fact that
 $\cB_Y^{r+1}(\log D)$ is locally free over $(\cO_Y)^p$. The assertion (2) is a consequence of (1), and (3) is straight-forward.
\end{pf}
\noindent
We turn to the proof of \eqref{eq4-3}, and compute the map $1-C$ in \eqref{eq4-3}:
\begin{equation}\label{eq4-5} 1-C : R^c i^!\cZ_Y^r (\log D) \lra R^c i^!\Omega_Y^r (\log D)),  \end{equation}
using the filtration $\Fil_\bullet$. Put
\[ \gr_m(-) := \Fil_m(-)/\Fil_{m-1}(-). \]
Since we have
\[ C \bigg( \frac{1}{(t_1t_2\dotsb t_c)^{pm}}\, \omega \bigg)=\frac{1}{(t_1t_2\dotsb t_c)^m}\, C(\omega)
 \quad \hbox{ for } \quad \omega \in \cZ_Y^r, \]
the Cartier operator preserves $\Fil_\bullet$ (by Lemma \ref{lem4-2}\,(2)) and induces a map
\[ \gr_m(C) : \gr_m(R^c i^!\cZ_Y^r (\log D)) \lra \gr_m(R^c i^!\Omega_Y^r (\log D)), \]
which is the zero map for $m \ge 2$. Hence we have
\[ \hbox{the kernel of \eqref{eq4-5}} \subset \Fil_1(R^c i^!\cZ_Y^r (\log D)). \]
Let $D_j \subset Y$ be as we defined in Lemma \ref{lem4-2}\,(3) for $1 \le j \le a$.
Fix an arbitrary $\omega \in \cZ_Y^r$, and assume that
\[ x:=\bigg[\frac 1{(t_1t_2\dotsb t_c)^p}\, \omega \bigg] \in \Fil_1(R^c i^!\cZ_Y^r (\log D)) \] belongs to the kernel of $\eqref{eq4-5}$.
To show \eqref{eq4-3}, we have to prove that
\begin{equation}\label{eq4-6}  \omega=\sum_{j=1}^c \ t_j^p \, \eta_j \;\;
 \hbox{ for some } \begin{cases}\eta_j \in \Omega_Y^r (\log D_j) & \hbox{($1 \le j \le a$)} \\ \eta_j \in \Omega_Y^r & \hbox{($a < j \le c$)} \end{cases}
\end{equation}
in $\Omega_Y^r (\log D)$, which is equivalent to that $x=0$ in $R^c i^!\cZ_Y^r (\log D)$ by Lemma \ref{lem4-2}\,(1) and \eqref{eq-4.2.1}. Since $x=C(x)$ in $\Fil_1(R^c i^!\Omega_Y^r (\log D))$ by assumption, we have
\[ \frac 1{(t_1t_2\dotsb t_c)^p}\, \omega-\frac 1{t_1t_2\dotsb t_c}\, C(\omega) =
 \sum_{j=1}^c \ \frac {1}{(t_1t_2\dotsb t_{j-1}t_{j+1}\dotsb t_c)^p} \, \alpha_j \quad \hbox{ in } \;\; \tau_*\tau^*\Omega_Y^r \]
for some $\alpha_j \in \Omega_Y^r (\log D_j)$\,($1 \le j \le a$) and $\alpha_j \in \Omega_Y^r$ ($a< j \le c$) by Lemma \ref{lem4-2}\,(3). This implies
\begin{equation}\label{eq4-7}
\omega =  (t_1t_2\dotsb t_c)^{p-1}C(\omega) +\sum_{j=1}^c \ t_j^p \, \alpha_j  \quad \hbox{ in } \;\; \Omega_Y^r,
\end{equation}
and our task is to prove that
\begin{equation}\label{eq4-7'}
(t_1t_2\dotsb t_c)^{p-1}C(\omega) = \sum_{j=1}^c \ t_j^p \, \zeta_j \;\;
 \hbox{ for some } \begin{cases} \zeta_j \in \Omega_Y^r (\log D_j) & \hbox{($1 \le j \le a$)} \\ \zeta_j \in \Omega_Y^r & \hbox{($a < j \le c$)}. \end{cases}
\end{equation}
Take a basis $\{ dt_\lam \}_{\lam \in \vL}$ of $\Omega^1_Y$ over $\cO_Y$ which contain $dt_1,\dotsc,dt_c$
 (see the condition ($\star$)).
Fix an ordering on $\vL$ and let $J$ be the set of all $r$-tuples $\bs{\lam}=(\lam_1,\lam_2,\dotsc,\lam_r)$ of elements of $\vL$ satisfying \[ \lam_1 < \lam_2 < \dotsb < \lam_r. \]
For $\bs{\lam}=(\lam_1,\lam_2,\dotsc,\lam_r) \in J$, put
\[ dt_{\bs{\lam}}:=dt_{\lam_1}\wedge dt_{\lam_2} \wedge \dotsb \wedge dt_{\lam_r} \in \Omega^r_Y. \]
Using the basis $\{dt_{\bs{\lam}}\}_{\bs{\lam} \in J}$ of $\Omega^r_Y$ over $\cO_Y$, we decompose $C(\omega)$ as
\begin{equation}\label{eq4-7''} C(\omega) = \epsilon + \delta \qquad( \epsilon, \delta \in \Omega^r_Y). \end{equation}
Here $\epsilon$ is an $\cO_Y$-linear combination of $dt_{\bs{\lam}}$'s which contain $dt_j$ for some $1 \le j \le a$,
 and $\delta$ is an $\cO_Y$-linear combination of $dt_{\bs{\lam}}$'s which contain none of $dt_1,\dotsc,dt_a$.
We have
\begin{equation}\label{eq4-7'''}
(t_1t_2\dotsb t_c)^{p-1}\epsilon = \sum_{j=1}^a \ t_j^p\, \beta_j \; \hbox{ for some $\beta_j \in \Omega_Y^r (\log D_j)$ \; ($1 \le j \le a$)},
\end{equation}
which implies \eqref{eq4-7'} if $\delta=0$. Otherwise, we proceed as follows. Put \[ y := (t_2t_3\dotsb t_c)^{p-1}\delta \in \Omega_Y^r. \] Since $d\omega=0$, the equalities \eqref{eq4-7} and \eqref{eq4-7'''} imply
\[ t_1^{p-2}\, dt_1 \wedge y = t_1^{p-1} \, dy+\sum_{j=1}^c \ t_j^p \, \beta'_j \;\;
 \hbox{ for some } \begin{cases} \beta'_j \in \Omega_Y^r (\log D_j) & \hbox{($1 \le j \le a$)} \\ \beta'_j \in \Omega_Y^r & \hbox{($a < j \le c$)}. \end{cases}
 \]
By this equality and the assumption on $\delta$, we see that
\begin{equation}\label{eq4-8}
 t_1^{p-2}\, dt_1 \wedge y = \sum_{j=1}^c \ t_j^p \, \beta''_j \;\;
 \hbox{ for some } \begin{cases} \beta''_1 \in \Omega_Y^r (\log D_1) \; & \\ \beta''_j \in \Omega_Y^r & \hbox{($2 \le j \le c$)}. \end{cases} \end{equation}
To proceed with the proof of \eqref{eq4-7'}, we need the following lemma:
\begin{lem}\label{lem4-3}
Let $f$ and $g$ be the following $\cO_Y$-linear maps, respectively{\rm:}
{\allowdisplaybreaks
\begin{align*}
f : \Omega_Y^r & \lra \Omega_Y^{r+1}, \quad z \mapsto dt_1 \wedge z, \\
g : \Omega_Y^r(\log D_1) & \lra \Omega_Y^{r+1}(\log D_1), \quad z \mapsto \frac{dt_1}{t_1} \wedge z.
\end{align*}
}Then the following holds{\rm :}
\begin{enumerate}
\item[(1)]
For $m \ge 1$, we have
\[ \Image(f) \cap \bigg(t_1^m \,\Omega^{r+1}_Y(\log D_1) + \sum_{j=2}^c \, t_j^m \,\Omega^{r+1}_Y \bigg) = \; t_1^m \,\Image(g) + \sum_{j=2}^c \, t_j^m \,\Image(f). \]
\item[(2)]
For $m,n \ge 0$, we have
\[ t_1^n \,\Omega^{r+1}_Y \cap \bigg(\sum_{j=2}^c \, t_j^m \,\Image(f) \bigg) = \sum_{j=2}^c \, t_1^n\,t_j^m \,\Image(f). \]
\item[(3)]
For $m \ge 0$, the sequence
\[\xymatrix{t_1^m \, \Omega_Y^{r-1} \ar[r]^-{f}& t_1^m\, \Omega_Y^{r} \ar[r]^-{f} & t_1^m\, \Omega_Y^{r+1} }\]
is exact.
\end{enumerate}
\end{lem}
\begin{pf}
These assertions follow from linear algebra over $\cO_Y$. The details are straight-forward and left to the reader.
\end{pf}
\noindent
By \eqref{eq4-8} and Lemma \ref{lem4-3}\,(1) for $m=p$, we have
\begin{equation}\label{eq4-9-0}
 t_1^{p-2}\, dt_1 \wedge y =  t_1^{p-1} \, dt_1 \wedge \gamma_1 + \sum_{j=2}^c \ t_j^p \, dt_1 \wedge \gamma_j
\end{equation}
for some $\gamma_1 \in \Omega_Y^r(\log D_1)$ and some $\gamma_j \in \Omega_Y^r$ ($2 \le j \le c$). If $p=2$, this equality implies
\begin{equation}\label{eq4-9}
 t_1^{p-2}\, dt_1 \wedge y =  t_1^{p-2} \, dt_1 \wedge \bigg(t_1 \gamma_1 + \sum_{j=2}^c \ t_j^p \,\gamma_j\bigg).
\end{equation}
If $p \ge 3$, we obtain \eqref{eq4-9} from \eqref{eq4-9-0} by replacing $\gamma_2,\dotsc,\gamma_c$ suitably in $\Omega_Y^r$ ($2 \le j \le c$), where we used Lemma \ref{lem4-3}\,(2) for $(m,n)=(p,p-2)$. Noting that $t_1 \gamma_1$ and $t_j^p \,\gamma_j$\,($2 \le j \le c$) belong to $\Omega^r_Y$, we have
\[ t_1^{p-2}\,y = t_1^{p-2}dt_1 \wedge \vt + t_1^{p-1} \gamma_1 + \sum_{j=2}^c \ t_j^p\,(t_1^{p-2}\gamma_j) \; \hbox{ for some $\vt \in \Omega^{r-1}_Y$}  \]
by \eqref{eq4-9} and Lemma \ref{lem4-3}\,(3) for $m=p-2$. Thus we have
\[  (t_1t_2\dotsb t_c)^{p-1}\delta = t_1^{p-1}\,y = t_1^p\bigg(\frac{dt_1}{t_1} \wedge \vt + \gamma_1\bigg) + \sum_{j=2}^c \ t_j^p\,(t_1^{p-1}\gamma_j), \]
which implies \eqref{eq4-7'} (see also \eqref{eq4-7''}, \eqref{eq4-7'''}). This completes the proof of \eqref{eq4-6}, \eqref{eq4-3} and Theorem \ref{thm4-1}.
\end{pf*}
\begin{thm}\label{cor4-1}
Let $S,p$ and $X$ be as in Setting {\rm\ref{set2-1}}, and let $c$ be a positive integer. Let $D$ be a non-empty normal crossing divisor on $X$, and let $Z$ be a closed subset of $X$ which has codimension $\ge c$. Put $U:=X-D$. Then we have \begin{equation}\label{eq4-10} H^q_Z(X,\fT_n(r)_{(X,D)}) \simeq  \begin{cases} 0 & \hbox{{\rm(}$q < r+c${\rm)}} \\ H^{r+c}_{Z \cap U}(U,\fT_n(r)_U) \quad & \hbox{{\rm(}$q = r+c${\rm)}}. \end{cases}\end{equation} In particular, when $Z$ has pure codimension $r$ on $X$, we have \begin{equation}\label{eq4-11} H^q_Z(X,\fT_n(r)_{(X,D)}) \simeq \begin{cases} 0 & \hbox{{\rm(}$q < 2r${\rm)}} \\ \bZ/p^n[Z^0 \cap U] \quad & \hbox{{\rm(}$q = 2r${\rm)}}, \end{cases}\end{equation} where $\bZ/p^n[Z^0 \cap U]$ means the free $\bZ/p^n$-module generated by the set $Z^0 \cap U$.
\end{thm}
\begin{pf}
Admitting \eqref{eq4-10}, we obtain \eqref{eq4-11} from the purity of $\fT_n(r)_U$ (\cite{S2} Theorem 4.4.7, Corollary 4.4.9) and the absolute purity of $\mu_{p^n}^{\otimes r}$ on $U[p^{-1}]$ (\cite{FG}). To show \eqref{eq4-10}, we divide the problem into the following 4 cases:
\begin{enumerate}
\item[(1)] $Z \subset D \cap Y$ \qquad \qquad \qquad\qquad\;\; (2) $Z \subset D$ and $Z \not \subset Y$
\item[(3)] $Z$ is arbitrary and $q<r+c$ \qquad \, (4) $Z$ is arbitrary and $q=r+c$.
\end{enumerate}
The case (1) follows from Theorems \ref{thm2-1}, \ref{thm4-1} (see also Remark \ref{rem2-1}) and the same arguments as in \cite{S2} Theorem 4.4.7. The case (2) follows from the case (1) and the same arguments as in loc.\ cit.\ Corollary 4.4.9. The case (3) also follows from similar arguments as for the previous cases. To prove the case (4), write $Z$ as the union of closed subsets
\[ Z = Z_1 \cup Z_2, \]
where $Z_1$ has pure codimension $c$ on $X$, $Z_2$ has codimension $\ge c+1$ on $X$ and we suppose that no irreducible components of $Z_2$ are contained in $Z_1$. Since $Z_1 \cap Z_2$ has codimension $\ge c+2$ on $X$, the assertion in the cases (1)--(3) and a Mayer-Vietoris long exact sequence
\begin{align*} \dotsb & \lra H^{r+c}_{Z_1 \cap Z_2}(X,\fT_n(r)_{(X,D)}) \lra \bigoplus_{j=1,2}\, H^{r+c}_{Z_j}(X,\fT_n(r)_{(X,D)}) \lra H^{r+c}_Z(X,\fT_n(r)_{(X,D)}) \\ & \lra H^{r+c+1}_{Z_1 \cap Z_2}(X,\fT_n(r)_{(X,D)}) \lra \dotsb \end{align*}
imply isomorphisms
\[  H^{r+c}_Z(X,\fT_n(r)_{(X,D)}) \simeq \bigoplus_{j=1,2}\, H^{r+c}_{Z_j}(X,\fT_n(r)_{(X,D)}) \simeq H^{r+c}_{Z_1}(X,\fT_n(r)_{(X,D)}). \]
Similarly, we have $H^{r+c}_{Z\cap U }(U,\fT_n(r)_U) \simeq H^{r+c}_{Z_1 \cap U}(U,\fT_n(r)_U)$ by the purity of $\fT_n(r)_U$ (\cite{S2} Corollary 4.4.9). Hence we may assume that $Z$ has pure codimension $c$ on $X$. Moreover we may assume that no irreducible components of $Z$ are contained in $D$ by the cases (1)+(2) and a similar devissage argument. Then noting that $Z \cap D$ has codimension $\ge c+1$ on $X$, we obtain the case (4) from the cases (1)+(2) and a long exact sequence
\begin{align*} \dotsb & \lra H^{r+c}_{Z \cap D}(X,\fT_n(r)_{(X,D)}) \lra H^{r+c}_Z(X,\fT_n(r)_{(X,D)}) \lra H^{r+c}_{Z \cap U}(U,\fT_n(r)_U) \\ & \lra H^{r+c+1}_{Z \cap D}(X,\fT_n(r)_{(X,D)}) \lra \dotsb. \end{align*}
This completes the proof of Theorem \ref{cor4-1}.
\end{pf}
\begin{defn}\label{def4-1}
Let $X$ and $D$ be as in Theorem {\rm\ref{cor4-1}}, and let $C$ be a cycle on $U:=X-D$ of codimension $r$.
Let $W \subset U$ be the support of $C$ and let $\ol W$ be its closure in $X$. Then we define the cycle class 
\[ \cl_X(C) \in H^{2r}_{\ol W}(X,\fT_n(r)_{(X,D)}) \]
as the inverse image of $\cl_U(C)$ {\rm(}\cite{S2} {\rm 5.1.2)} under the isomorphism in Theorem {\rm\ref{cor4-1}}
\[ H^{2r}_{\ol W}(X,\fT_n(r)_{(X,D)}) \isom H^{2r}_W(U,\fT_n(r)_U). \]
\end{defn}
\medskip

\section{\bf Cycle class map}\label{sect5}
\medskip
Let $S,p$ and $X$ be as in Setting \ref{set2-1}. In this section, we construct a cycle class morphism
\begin{equation}\label{eq5-1} \cl_X : \bZ(r)_X^\et \otimes^{\bL}\bZ/p^n \lra \fT_n(r)_X  \quad \hbox{ in } \;\; D^-(X_\et,\bZ/p^n), \end{equation} following the method of Bloch \cite{B3} \S4 (cf.\ \cite{GL} \S3, see also  Remark \ref{rem5-2} below). Here $\bZ(r)_X^\et$ denotes the \'etale sheafification on $X$ of the presheaf of cochain complexes \[  U \longmapsto z^r(U,\star)[-2r], \] and $z^r(U,\star)$ denotes Bloch's cycle complex \[\begin{matrix} \dotsb \lra & z^r(U,q) & \os{d_q}{\lra} & z^r(U,q-1) & \os{d_{q-1}}{\lra} \dotsb \os{d_1}{\lra} & z^r(U,0). \\ & \hbox{{\scriptsize(degree $-q$)}} & & \hbox{{\scriptsize(degree $-q+1$)}} & & \hbox{{\scriptsize(degree $0$)}}  \end{matrix}\] We review the definition of this complex briefly, which will be useful later. Let $\vD^\star$ be the standard cosimplicial scheme over $\Spec(\bZ)$: \[ \nvD q := \Spec(\bZ[t_0,t_1,\dotsc, t_q]/(t_0+t_1+\dotsb+ t_q=1)). \] A face of $\nvD q$ (of codimension $a \ge 1$) is a closed subscheme defined by the equation \[ t_{i_1}=t_{i_2}=\dotsb=t_{i_a}=0 \;\;\; \hbox{ for some \; $0 \le i_1 <i_2 < \dotsb< i_a \le q$}. \]
Now $z^r(U,q)$ is defined as the free abelian group generated by the set of all integral closed subschemes on $U \times \nvD q$ of codimension $r$ which meet all faces of $U \times \nvD q$ properly. Here a face of $U \times \nvD q$ means the product of $U$ and a face of $\nvD q$. Noting that the faces of codimension $1$ are effective Cartier divisors, we define the differential map $d_q$ as the alternating sum of pull-back maps along the faces of codimension $1$, which defines the above complex $z^r(U,\star)$. 
\par
We fix a projective completion $\bvD q$ of $\nvD q$ as follows:
\[ \bvD q := \Proj(\bZ[T_0,T_1,\dotsc, T_q,T_\infty]/(T_0+T_1+\dotsb+T_q=T_\infty)). \]
Let $\vH q \subset \bvD q$ be the hyperplane at infinity:
\[ \vH q : T_\infty=0. \]
The following proposition plays a key role in our construction of the morphism \eqref{eq5-1}.
\begin{prop}\label{lem5-1}
Let $q$ and $r$ be integers with $q,r \ge 0$, and let $U$ be \'etale of finite type over $X$.
Let $\vS^{r,q}$ be the set of all closed subsets on  $U \times \nvD q$ of pure codimension $r$ which meet the faces of $U \times \nvD q$ properly.
For $W \in \vS^{r,q}$, let $\ol W$ be the closure of $W$ in $U \times \bvD q$. 
Then{\rm:}
\begin{enumerate}
\item[(1)]
There is an isomorphism
\[ \cl^{r,q} : z^r(U,q) \otimes \bZ/p^n \isom \varinjlim_{W \in \vS^{r,q}} \ H^{2r}_{\ol W}\left(U \times \bvD q,\fT_n(r)_{(U \times \bvD q,U \times \vH q)}\right) \]
sending a cycle $C \in z^r(U,q)$ to the cycle class $\cl_{U \times \bvD q}(C)$
 {\rm(}see Definition {\rm\ref{def4-1}}{\rm)}.
\item[(2)]
For $W \in \vS^{r,q}$, the natural morphism
\[ \tau_{\le 2r}\,R\vG_{\ol W}\left(U \times \bvD q,\fT_n(r)_{(U \times \bvD q,U \times \vH q)}\right)
 \lra H^{2r}_{\ol W}\left(U \times \bvD q,\fT_n(r)_{(U \times \bvD q,U \times \vH q)}\right)[-2r] \]
 is an isomorphism in the derived category of $\bZ/p^n$-modules.
\item[(3)]
Let $C$ be a cycle which belongs to $z^r(U,q)$, and let $W$ be the support of $C$ {\rm(}note that $W$ belongs to $\vS^{r,q}${\rm)}.
Let $\ol i : U \times \bvD {q-1} \hra U \times \bvD q$ be the closure of a face map $i : U \times \nvD {q-1} \hra U \times \nvD q$. Then the pull-back map
\[ \ol i{}^* : H^{2r}_{\ol W}\left(U \times \bvD q,\fT_n(r)_{(U \times \bvD q,U \times \vH q)}\right)
  \to H^{2r}_{\ol i{}^{-1}(\ol W)}\left(U \times \bvD {q-1},\fT_n(r)_{(U \times \bvD {q-1}, U \times \vH {q-1})}\right) \]
sends the cycle class $\cl_{U \times \bvD q}(C)$ to $\cl_{U \times \bvD {q-1}}(i^*C)$, where $i^*C$ denotes the pull-back of the cycle $C$ along $i$.
\end{enumerate}
\end{prop}
\begin{pf}
(1) and (2) follow from Theorem \ref{cor4-1} and the definition of $z^r(U,q)$.
We show (3). By Theorem \ref{cor4-1}, it is enough to show that the pull-back map
\[ i^* : H^{2r}_W\left(U \times \nvD q,\fT_n(r)_{U \times \nvD q}\right)
  \lra H^{2r}_{i^{-1}(W)}\left(U \times \nvD {q-1},\fT_n(r)_{U \times \nvD {q-1}}\right) \]
sends the cycle class $\cl_{U \times \nvD q}(C)$ to $\cl_{U \times \nvD {q-1}}(i^*C)$. To show this, put
\[ \cU:=U \times \nvD q, \quad \cD:=i(U \times \nvD {q-1}) \subset \cU
 \quad \hbox{ and } \quad \cV:=\cU - \cD \] and let $t \in \vG(\cU,\cO)$ be a defining equation of $\cD$. Noting that $i^{-1}(W)=W \cap \cD$, consider the following diagram:
\[\xymatrix{ H^{2r}_W(\cU,\fT_n(r)_\cU) \ar[d]_{\alpha} \ar[r]^-{i^*} &
 H^{2r}_{W \cap \cD}(\cD,\fT_n(r)_\cD) \ar[d]^{i_*}_{\wr\hspace{-2pt}} \\
 H^{2r+1}_{W \cap \cV}(\cV,\fT_n(r+1)_\cV) \ar[r]^-{-\delta} \ar[d] &
 H^{2r+2}_{W \cap \cD}(\cU,\fT_n(r+1)_\cU) \ar[d]_-{\wr\hspace{-2pt}} \\
 \us{y \in W^0}\bigoplus \, H^{2r+1}_y(\cV,\fT_n(r+1)_\cV) \ar@<5pt>[r]^-{-\delta}  &
 \us{x \in (W \cap \cD)^0}\bigoplus \, H^{2r+2}_x(\cU,\fT_n(r+1)_\cU) \\
 \ar[u]^-g_-{\hspace{-2pt}\wr} \us{y \in W^0}\bigoplus \, \kappa(y)^\times/p^n \ar@<5pt>[r]^-{\div} &
 \ar[u]_-{g'}^-{\wr\hspace{-2pt}} \us{x \in (W \cap \cD)^0}\bigoplus \, \bZ/p^n.  }\]
Here we defined $\alpha$ by sending $\omega \in H^{2r}_W(\cU,\fT_n(r)_\cU)$ to $\{ t|_{\cV} \} \cup \omega|_{\cV}$, where $\{ t|_{\cV} \}$ denotes the symbol in $H^1(\cV,\fT_n(1)_\cV)$ arising from $t|_{\cV} \in \vG(\cV,\cO^\times)$. The arrows $i_*$, $g$ and $g'$ are Gysin isomorphisms (cf.\ \cite{S2} Theorem 4.4.7). The arrows $\delta$ are boundary maps of localization long exact sequences, and the central square commutes obviously. The arrow $\div$ denotes the divisor map, and the bottom square commutes by the compatibility in loc.\ cit.\ Theorem 6.1.1 (for $x$ with $\ch(x)=p$) and \cite{JSS} Theorem 1.1.1 (for $x$ with $\ch(x)\ne p$).
The top square commutes by the following equalities ($\omega \in H^{2r}_W(\cU,\fT_n(r)_\cU)$):
\[ i_* \circ i^*(\omega)=\cl_{\cU}(\cD) \cup \omega=-\delta(\{ t|_{\cV} \}) \cup \omega=-\delta(\{ t|_{\cV} \} \cup \omega|_\cV)=-\delta\circ \alpha(\omega), \]
where the first (resp.\ second) equality follows from the projection formula in \cite{S2} Corollary 7.2.4 (resp.\ the same compatibility as mentioned before).
Hence the assertion follows from the transitivity of Gysin maps for $W \cap \cD \hra \cD \hra \cU$ (loc.\ cit.\ Corollary 6.3.3) and the fact that the divisor map $\div$ sends $t|_y \in \kappa(y)^\times/p^n$ to the cycle $i^*[y]$, where $[y]$ means the cycle on $\cU$ given by the closure of $y$. This completes the proof.
\end{pf}
We construct the morphism \eqref{eq5-1}.
For $U$ as in the proposition, there is a diagram of cochain complexes concerning $\bullet$ (see Corollary \ref{cor2-1}\,(2) for $G_n(r)^\bullet$):
\[\xymatrix{ \us{\phantom{a}}z^r(U,q) \otimes \bZ/p^n[-2r] \ar@<2.4pt>[r]^-{\cl^{r,q}}_-\sim & \displaystyle  \varinjlim_{W \in \vS^{r,q}} \ H^{2r}_{\ol W}\left(U \times \bvD q,\fT_n(r)_{(U\times \bvD q,U \times \vH q)}\right)[-2r] \\
\us{\phantom{(}}\vG \left(U \times \bvD q,G_n(r)^\bullet_{(U \times \bvD q,U \times \vH q)}\right) & \ar@<-2.4pt>[l]_-{\beta^{r,q}} \displaystyle \varinjlim_{W \in \vS^{r,q}} \ \tau_{\le 2r}\, \vG_{\ol W}\left(U \times \bvD q,G_n(r)^\bullet_{(U\times \bvD q,U \times \vH q)}\right)  \ar[u]_{\alpha^{r,q}}. }\]
Here $\alpha^{r,q}$ and $\beta^{r,q}$ are natural maps of complexes, which are obviously contravariant for the face maps $U \times \nvD {q-1} \hra U \times \nvD q$. The top arrow is also contravariant for these face maps by Proposition \ref{lem5-1}\,(3). Hence we get a diagram of double complexes concerning $(\star,\bullet)$
\[\xymatrix{ \us{\phantom{a}}z^r(U,\star) \otimes \bZ/p^n[-2r] & & \ar@<-2.4pt>[ll]_-{(\cl^{r,\star})^{-1} \circ \alpha^{r,\star}} \displaystyle \varinjlim_{W \in \vS^{r,\star}} \ \tau_{\le 2r}\, \vG_{\ol W}\left(U \times \bvD \star,G_n(r)^\bullet_{(U \times \bvD \star,U \times \vH \star)}\right) \ar[d]^{\beta^{r,\star}} \\
 \us{\phantom{a}}\vG(U,G_n(r)^\bullet_U) \ar@<2.4pt>[rr]^-{\text{canonical}} &  & \us{\phantom{(}}\vG \left(U \times \bvD \star,G_n(r)^\bullet_{(U \times \bvD \star,U \times \vH \star)}\right) , }\]
where the differential maps in the direction of $\star$ are the alternating sums of pull-back maps along the faces of codimension $1$, and the bottom map is the inclusion to the factor at $\star=0$. The top and the bottom arrows are quasi-isomorphisms on the associated total complexes by Proposition \ref{lem5-1}\,(1), (2) and Corollary \ref{cor3-1}, respectively.
We thus obtain the desired cycle class morphism \eqref{eq5-1} in $D^-(X_\et,\bZ/p^n)$ by sheafifying the diagram of total complexes.
\begin{rem}\label{rem5-1}
Assume that $S$ is local. Then admitting the Bloch-Kato conjecture for norm residue symbols \cite{BK} {\rm\S3} {\rm(\cite{Vo1}, \cite{Vo2}, \cite{We})}, one can easily see that the morphism \eqref{eq5-1} induces isomorphisms {\rm(}cf.\ \cite{S2} Conjecture {\rm1.4.1)}
\begin{align*}
 \tau_{\le r}\,(\bZ(r)_X^\et \otimes^{\bL}\bZ/p^n) & \isom \fT_n(r)_X  \qquad\qquad \hbox{ in } \;\; D^b(X_\et,\bZ/p^n), \\
 \tau_{\le r}\,(\bZ(r)_X^\zar \otimes^{\bL}\bZ/p^n) & \isom \tau_{\le r} R\ve_* \fT_n(r)_X  \quad\hspace{0.6pt} \hbox{ in } \;\; D^b(X_\zar,\bZ/p^n),
\end{align*}
by the Suslin-Voevodsky theorem \cite{SV} {\rm(}cf.\ \cite{GL}{\rm)} and similar arguments as in \cite{SS} {\rm\S A.2}. Here $\bZ(r)_X^\zar$ denotes the complex of Zariski sheaves on $X$ defined as $U \mapsto z^r(U,\star)[-2r]$, and $\ve$ denotes the continuous map of sites $X_\et \to X_\zar$.
\end{rem}
\begin{rem}\label{rem5-2}
If $X$ is smooth over $S$, then $\bZ(r)_X^\zar$ is concentrated in $\le r$ by a result of Geisser \cite{Ge} Corollary {\rm4.4}. Hence his arguments in loc.\ cit.\ {\rm\S6}, Proof of Theorem {\rm 1.3} gives an alternative construction of \eqref{eq5-1} in this case.
\end{rem}
\medskip
\section{\bf de Rham-Witt cohomology and homology}\label{sect6}
\medskip
In this section we define two kinds of de Rham-Witt complexes for normal crossing varieties playing the roles of cohomology and homology, and relate them with the modified de Rham-Witt complex of Hyodo \cite{H2}. The main results of this section are Theorems \ref{thm6-1} and \ref{thm6-2} below.
\par
Let $k$ be a field of characteristic $p>0$, and let $Y$ be a normal crossing variety over $k$.
Let $U \subset Y$ be a dense open subset which is smooth over $k$, and let $\sigma : U \hra Y$ be the natural open immersion.
\begin{defn}
We define the complex $(\cohwitt Y n {\bullet},d)$ of \'etale sheaves on $Y$ as the differential graded subalgebra of $\sigma_*\witt U n {\bullet}$ generated by
\[ \Wn\cO_Y \quad \hbox{ and } \quad \lam_{Y,n}^1 (\subset \sigma_*\witt U n 1), \]
which we call {\it the cohomological de Rham-Witt complex of $Y$}. It is easy to see that $\cohwitt Y n {\bullet}$ does not depend on the choice of $U$. For $n=1$, $\cohwitt Y 1 {\bullet}$ is the same as the complex $\vL_Y^\bullet$ defined in \cite{S1} {\rm\S3.3}.
\end{defn}
\noindent
By the relations
\begin{equation}\label{eq6-1}
\begin{cases}
 V(aR(x))=V(a)x \in \witt U {n+1} r  \quad& \hbox{($a \in \Wn\cO_U$, $x \in \logwitt U {n+1} r$)} \\
 F(ax)=F(a)R(x) \in \witt U n r & \hbox{($a \in W_{n+1}\cO_U$, $x \in \logwitt U {n+1} r$)} \\
 R(ax)=R(a)R(x) \in \witt U n r & \hbox{($a \in W_{n+1}\cO_U$, $x \in \logwitt U {n+1} r$)} \\
\end{cases} \end{equation}
 (cf.\ \cite{I} I Th\'eor\`eme 2.17, Proposition 2.18), the operators $V,F,R$ on $\witt U n r$ induce operators
\begin{align*} V & : \cohwitt Y n r \lra \cohwitt Y {n+1} r, \\
F & : \cohwitt Y {n+1} r \lra \cohwitt Y n r, \\
R & : \cohwitt Y {n+1} r \lra \cohwitt Y n r, \end{align*}
which satisfy relations
\begin{align}
\label{eq6-2}
 &FV=VF=p, \quad FdV=d, \quad dF=pFd, \quad  Vd=pdV, \\
\label{eq6-2'}
 &RV=VR, \qquad RF=FR, \qquad Rd=dR.
\end{align}
The local structure of $\cohwitt Y n r$ can be written in terms of the usual Hodge-Witt sheaves of the strata of $Y$.
\begin{prop}\label{prop6-1}
\begin{enumerate}
\item[(1)]
Assume that $Y$ is simple, and let $Y_1,Y_2,\dotsb, Y_q$ be the irreducible components of $Y$.
Then there is an exact sequence on $Y_\et$
\[ 0 \lra \cohwitt Y n r \os{\check{r}^0}\lra \bigoplus_{|I|=1} \ \witt {Y_I} n r \os{\check{r}^1}\lra \bigoplus_{|I|=2} \ \witt {Y_I} n r \os{\check{r}^2}\lra \dotsb \os{\check{r}^{q-1}}\lra \bigoplus_{|I|=q} \ \witt {Y_I} n r \lra 0, \]
where the notation is the same as in Proposition {\rm\ref{lem1-2}}.
\item[(2)]
Assume that $Y$ is embedded into a smooth $k$-variety $\cY$ as a normal crossing divisor.
Let $i : Y \hra \cY$ be the closed immersion. Then there is a short exact sequence on $\cY_\et$
\[ 0 \lra \witt \cY n r(-\log Y) \lra \witt \cY n r \os{i^*}\lra i_* \cohwitt Y n r \lra 0. \]
\end{enumerate}
\end{prop}
\begin{pf}
(1) follows from the contravariant functoriality of Hodge-Witt sheaves of smooth varieties and the same arguments as for \cite{S1} Proposition 3.2.1. The assertion (2) follows from (1) and \cite{M} Lemma 3.15.1.
\end{pf}
We next introduce homological Hodge-Witt sheaves.
\begin{defn}
For $r \ge 0$, we define the \'etale sheaf $\homwitt Y n r$ on $Y$ as the $\Wn\cO_Y$-submodule of $\sigma_*\witt U n r$ generated by $\nu_{Y,n}^r (\subset \sigma_*\witt U n r)$.
Similarly as for $\cohwitt Y n r$, the sheaf $\homwitt Y n r$ does not depend on the choice of $U$.
For $n=1$, $\homwitt Y 1 r$ agrees with the sheaf $\varXi_Y^r$ defined in \cite{S1} {\rm\S2.5}.
\end{defn}
\noindent
Since $\lam_{Y,n}^r \subset \nu_{Y,n}^r$ by definition, we have $\cohwitt Y n r \subset \homwitt Y n r$,
which agree with $\witt Y n r$ when $Y$ is smooth (cf.\ \cite{S1} Remark 3.1.4).
The following local description will be useful later.
\begin{prop}\label{lem6-1}
\begin{enumerate}
\item[(1)]
Under the same setting and notation as in Proposition {\rm\ref{prop6-1}\,(1)}, there is a canonical ascending filtration $\F_a$\,{\rm ($a \ge 0$)} on $\homwitt Y n r$ satisfying \[ \F_0(\homwitt Y n r)=0 \quad \hbox{ and } \quad \gr^\F_a\,\homwitt Y n r \simeq \bigoplus_{|I|=a} \ \witt {Y_I} n {r-a+1}, \]
where the isomorphism for $a \ge 2$ depends on the fixed ordering of the irreducible components of $Y$.
\item[(2)]
Under the same setting and notation as in Proposition {\rm\ref{prop6-1}\,(2)}, there is a short exact sequence on $\cY_\et$
\[ 0 \lra \witt \cY n {r+1} \lra \witt \cY n {r+1} (\log Y) \os{\varrho}\lra i_* \homwitt Y n r \lra 0, \]
where $\varrho$ is induced by the Poincar\'e residue mapping
\[ \varrho_0 : \witt \cY n {r+1} (\log Y) \lra (i\sigma)_*\witt U n r . \]
\end{enumerate}
\end{prop}
\begin{pf}
We first recall the following fact due to Mokrane \cite{M} Proposition 1.4.5. Under the setting of (2), let $\F_a\witt \cY n m(\log Y)$\,($a \ge 0$) be the image of the product
\[ \witt \cY n a(\log Y) \times \witt \cY n {m-a} \lra \witt \cY n m(\log Y). \]
Put $\F_{-1}\witt \cY n m(\log Y):=0$. Then there are isomorphisms
\begin{equation}\label{eq6-10} \gr^\F_a \, \witt \cY n m (\log Y) \simeq  \begin{cases} \witt \cY n m & \hbox{($a=0$)} \\ \bigoplus_{|I|=a} \ \witt {Y_I} n {m-a} & \hbox{($a \ge 1$)}, \end{cases} \end{equation}
which are induced by Poincar\'e residue mappings for $a \ge 1$ and depend on the fixed ordering of the irreducible components of $Y$ for $a \ge 2$.
We prove (2). Because the complement $\vS:=Y-U$ has codimension $\ge 2$ on $\cY$ and the sheaves $\witt \cY n {r+1}$ and $\witt \cY n {r+1}(\log Y)$ are finitely successive extensions of locally free $\cO_{\cY}$-modules of finite rank (\cite{I} I Corollaire 3.9), the sequence \[ 0 \lra \witt \cY n {r+1} \lra \witt \cY n {r+1} (\log Y) \os{\varrho_0}\lra(i\sigma)_*\witt U n r \]
is exact by \eqref{eq6-10} for $\cY-\vS$. Hence it is enough to show that $\Image(\varrho_0)=i_* \homwitt Y n r$.
Let $\F_a(\nu_{Y,n}^r)$\,($a \ge 0$) be the filtration on $\nu_{Y,n}^r$ in \cite{S1} Proposition 2.2.1, which is defined under the assumption that $Y$ is simple and satisfies \[ \F_0(\nu_{Y,n}^r)=0 \quad \hbox{ and } \quad \gr^\F_a\,\nu_{Y,n}^r \simeq \bigoplus_{|I|=a} \ \logwitt {Y_I} n {r-a+1} \quad \hbox{(non-canonically)}. \]
Hence comparing this local description of $\nu_{Y,n}^r$ with \eqref{eq6-10} for $m=r+1$, we see that $\Image(\varrho_0)=i_* \homwitt Y n r$. As for (1), we define the desired filtration $\F_a(\homwitt Y n r)$ as the $\Wn\cO_Y$-submodule generated by $\F_a(\nu_{Y,n}^r)$. This filtration satisfies the required properties by (2) and again by \eqref{eq6-10}.
\end{pf}
Since the residue mapping $\varrho_0$ in Proposition \ref{lem6-1}\,(2) commutes with the operators $d,V,F,R$ (cf.\ \eqref{eq6-1}), these operators on $\witt U n {\bullet}$ induce operators
\begin{align*}
d & : \homwitt Y n r \lra \homwitt Y n {r+1}, \\
V & : \homwitt Y n r \lra \homwitt Y {n+1} r, \\
F & : \homwitt Y {n+1} r \lra \homwitt Y n r, \\
R & : \homwitt Y {n+1} r \lra \homwitt Y n r \end{align*}
by Proposition \ref{lem6-1}\,(2) (without the assumptions in Proposition \ref{lem6-1}\,(1) or (2)), which satisfy the relations listed in \eqref{eq6-2} and \eqref{eq6-2'}.
We call the resulting complex $(\homwitt Y n \bullet,d)$ {\it the homological de Rham-Witt complex of $Y$}.
\begin{rem}\label{rem6-1}
There are natural injective homomorphisms of complexes
\[ \bigoplus_{r=0}^{\dim (Y)} \lam_{Y,n}^r[-r] \hra \cohwitt Y n \bullet, \qquad \bigoplus_{r=0}^{\dim (Y)} \nu_{Y,n}^r[-r] \hra \homwitt Y n \bullet, \]
where the differentials on the complexes on the left hand side are defined as zero. Indeed the differentials on $\witt U n \bullet$ are zero on $\logwitt U n \bullet$.
\end{rem}
\noindent
The following Proposition \ref{prop6-2} and Theorem \ref{thm6-1} explain a fundamental relationship between $\cohwitt Y n \bullet$ and $\homwitt Y n \bullet$.
\begin{prop}\label{prop6-2}
The complex $\homwitt Y n \bullet$ is a $\cohwitt Y n {\bullet}$-module, that is, there is a natural $\Wn\cO_Y$-bilinear pairing
\begin{align}
\label{eq6-3}
\qquad \cohwitt Y n r \times \homwitt Y n s \lra \homwitt Y n {r+s} & \qquad \hbox{{\rm(}$r,s \ge 0${\rm)}} \\
\intertext{satisfying the Leibniz rule}
\notag
d(x \cdot y) = (dx) \cdot y + (-1)^r x \cdot dy & \qquad \hbox{{\rm(}$x \in \cohwitt Y n r$, $y \in \homwitt Y n s${\rm)}}. \\
\intertext{Moreover this pairing is compatible with $R$ and satisfies relations}
\notag
\begin{array}{r} \phantom{\big(} F(x \cdot y)=F(x) \cdot F(y)  \\ \phantom{\Big(} x \cdot V(y) = V(F(x) \cdot y) \\
\phantom{\big(} V(x) \cdot y = V(x \cdot F(y))  \end{array}
 & \qquad \hbox{{\rm(}$x \in \cohwitt Y n r$, $y \in \homwitt Y n s${\rm)}}. \end{align}
\end{prop}
\begin{pf}
The pairing \eqref{eq6-3} is induced by the product of $\witt U n \bullet$ and the biadditive pairing
\[ \lam_{Y,n}^r \times \nu_{Y,n}^s \lra \nu_{Y,n}^{r+s} \]
defined in \cite{S1} Definition 3.1.1. The properties of the pairing \eqref{eq6-3} follow from
 the corresponding properties of the product of $\witt U n \bullet$ (\cite{I} I Th\'eor\`eme 2.17, Proposition 2.18).
\end{pf}
\par
\begin{rem}\label{rem6-2}
When $k$ is perfect and $Y$ is the special fiber of a regular semistable family $X$ over a discrete valuation ring with residue field $k$, we have 
\[ \cohwitt Y n \bullet \subset \mwitt Y n \bullet \subset \homwitt Y n \bullet \]
by \cite{S1} Proposition {\rm 4.2.1}. Here $\mwitt Y n \bullet$ denotes the modified de Rham-Witt complex associated with $X$ \cite{H2}. The complex $\homwitt Y n \bullet$ does not in general have a product structure unless $Y$ is smooth.
\end{rem}
\begin{thm}\label{thm6-1}
Assume that $k$ is perfect and that $Y$ is proper over $k$. Put $b:=\dim(Y)$ and $\Wn:=\Wn(k)$. Then{\rm:}
\begin{enumerate}
\item[(1)]
There is a canonical trace map $\tr_n : \bH^{2b}(Y,\homwitt Y n \bullet) \to \Wn$, which is bijective if $Y$ is geometrically connected over $k$.
\item[(2)]
The pairing \[\xymatrix{ \bH^i(Y,\cohwitt Y n \bullet) \times \bH^{2b-i}(Y,\homwitt Y n \bullet) \ar[r]^-{\eqref{eq6-3}} &\bH^{2b}(Y,\homwitt Y n \bullet) \ar[r]^-{\tr_n} & \Wn }\] is a non-degenerate pairing of finitely generated $\Wn$-modules for any $i \ge 0$.
\end{enumerate}
\end{thm}
\begin{pf}
(1) Let $f_n$ be the canonical morphism $\Wn Y \to \Wn$, where $\Wn Y$ means the scheme consisting of the topological space $Y$ and the structure sheaf $\Wn\cO_Y$.
Since $\bH^{2b}(Y,\homwitt Y n \bullet) \simeq H^b(Y,\homwitt Y n b)$,  our task is to show that
\begin{equation}\label{eq6-4} \homwitt Y n b \isom f_n^!\Wn \quad \hbox{ in } \;\; D^+_{\qc}(\Wn Y,\Wn\cO_Y), \end{equation}
where the subscript $\qc$ means the triangulated subcategory consisting of quasi-coherent cohomology sheaves
 and $f_n^!$ means the twisted inverse image functor of Hartshorne \cite{Ha}.
If $Y$ is smooth, \eqref{eq6-4} is a theorem of Ekedahl \cite{E}.
In the general case, \eqref{eq6-4} is reduced to the smooth case by
Proposition \ref{lem6-1}\,(2) and the same arguments as for \cite{S1} Proposition 2.5.9.

(2) The case $n=1$ follows from \cite{S1} Corollary 2.5.11, Proposition 3.3.5. The case $n \ge 2$ follows from a standard induction argument on $n$ and the following lemma:
\end{pf}
\begin{lem}\label{lem6-2}
\begin{enumerate}
\item[(1)]
There are injective homomorphisms
\[ \xymatrix{ \ul{p} : \cohwitt Y {n-1} r \hra \cohwitt Y n r
 \quad \hbox{ and } \quad
 \ul{p} : \homwitt Y {n-1} r \hra \homwitt Y n r }\]
induced by the multiplication by $p$ on $\cohwitt Y n r$ and $\homwitt Y n r$, respectively.
\item[(2)] The following natural projections of complexes are quasi-isomorphisms{\rm:}
\[ \cohwitt Y n {\bullet}/\ul{p}(\cohwitt Y {n-1} {\bullet}) \lra \vL_Y^\bullet
 \quad \hbox{ and } \quad
\homwitt Y n {\bullet}/\ul{p}(\homwitt Y {n-1} {\bullet}) \lra \varXi_Y^\bullet. \]
\end{enumerate}
\end{lem}
\begin{pf}
When $Y$ is smooth, the assertions are due to Illusie \cite{I} I Proposition 3.4 and Corollaire 3.15.
In the general case, (1) follows from that for the dense open subset $U \subset Y$ we fixed before.
To prove (2), we may assume that $Y$ is simple. Then the assertions follow from those for the strata of $Y$ and Propositions \ref{prop6-1}\,(1) and \ref{lem6-1}\,(1).
\end{pf}
In the rest of this section, we work under the following setting.
Let $A$ be a discrete valuation ring with perfect residue field $k$, and let $X$ be a regular semistable family over $A$. Put $Y:=X \otimes_A k$. Recall that
\[ \cohwitt Y n \bullet \subset \mwitt Y n \bullet \subset \homwitt Y n \bullet \]
by Remark \ref{rem6-2}.
We define a Frobenius endomorphism $\varphi$ on these complexes by $p^mF$ on degree $m \ge 0$.
There is a short exact sequence of complexes with Frobenius action
\begin{equation}\label{eq6-5}
 0 \lra \mwitt Y n {\bullet}(-1)[-1] \os{\frac{dt}{t}\wedge}{\lra} \mtwitt Y n \bullet \lra \mwitt Y n \bullet \lra 0
\end{equation}
by \cite{H2} (1.4.3), where the complex $\mwitt Y n {\bullet}(-1)$ means the complex $\mwitt Y n {\bullet}$ with Frobenius endomorphism $p\cdot\varphi$. The monodromy operator
\[ N : \mwitt Y n \bullet \lra \mwitt Y n {\bullet}(-1) \]
is defined as the connecting morphism associated with this sequence.
\begin{prop}\label{prop6-3}
There is a short exact sequence of complexes with Frobenius action
\begin{equation}\label{eq6-6}
 0 \lra \cohwitt Y n \bullet \lra \mtwitt Y n \bullet \lra \homwitt Y n {\bullet}(-1)[-1] \lra 0
\end{equation}
fitting into a commutative diagram of complexes
\[\xymatrix{ & \mwitt Y n {\bullet}(-1)[-1] \ar[d]_{\frac{dt}{t}\wedge} \ar@{_{(}->}[rd] & \\
 \cohwitt Y n \bullet \ar[r] \ar@{_{(}->}[rd]  & \mtwitt Y n \bullet \ar[r] \ar[d] & \homwitt Y n {\bullet}(-1)[-1] \\
 & \mwitt Y n \bullet. & }\]
Consequently, we obtain a complex of \,$\Wn(k)$-modules with Frobenius action for $i \ge 0$
\begin{equation}\label{eq6-7}
 \bH^i(Y,\cohwitt Y n \bullet) \to \bH^i(Y,\mwitt Y n \bullet) \os{N}\lra \bH^i(Y,\mwitt Y n \bullet)(-1)
 \to \bH^i(Y,\homwitt Y n \bullet)(-1). \end{equation}
\end{prop}
\begin{pf}
Let $\sigma : U \hra Y$ be as before and define a homomorphism
 $\partial : \mtwitt Y n r \to \sigma_*\witt U n {r-1}$ as the composite
\[ \partial : \mtwitt Y n r \lra \sigma_*\sigma^*\mtwitt Y n r \isom
 \sigma_*(\witt U n r \oplus \witt U n {r-1}) \os{\pr_2}\lra \sigma_*\witt U n {r-1}, \]
where the second arrow is obtained from the fact that the sequence \eqref{eq6-5} splits on $U$.
It is easy to see that $\partial$ satisfies
\[ d \partial = \partial d \qquad \hbox{and} \qquad  \partial F = p F \partial, \]
which induces a map of complexes $\partial^\bullet : \mtwitt Y n \bullet \lra \sigma_* \witt U n {\bullet}(-1)[-1]$.
To show that $\partial^\bullet$ induces the exact sequence \eqref{eq6-6}, we may assume the assumption in Proposition \ref{prop6-1}\,(2). Then we have an isomorphism of complexes
\[ \mtwitt Y n \bullet \simeq \witt \cY n \bullet (\log Y)/\witt \cY n \bullet (-\log Y) \qquad\hbox{(\cite{H2} p.\ 247 Lemma)} \]
and the assertion follows from Propositions \ref{prop6-1}\,(2) and \ref{lem6-1}\,(2).
The commutativity of the diagram in the proposition is straight-forward and left to the reader.
\end{pf}

Now we prove that the monodromy-weight conjecture for log crystalline cohomology implies an invariant cycle `theorem' (cf.\ \cite{I2} 2.4.5), which will be useful in \S\S\ref{sect8}--\ref{sect7} below. Let $A$ be a discrete valuation ring with finite residue field $k$, and $X$ be a regular scheme which is projective flat over $A$ with strict semistable reduction. Put $K_0:=\Frac(W(k))$ and $Y:=X \otimes_A k$. For integer $i \ge 0$, we define
\begin{align}
\label{eq8-C} C^i:=\qp \otimes_{\zp} \varprojlim_{n \ge 1} \ \bH^i(Y,\cohwitt Y n \bullet), \\
\label{eq8-D} D^i:=\qp \otimes_{\zp} \varprojlim_{n \ge 1} \ \bH^i(Y,\mwitt Y n \bullet), \\
\label{eq8-E} E^i:=\qp \otimes_{\zp} \varprojlim_{n \ge 1} \ \bH^i(Y,\homwitt Y n \bullet),
\end{align}
which are finite-dimensional vector spaces over $K_0$. The content of the monodromy-weight conjecture due to Mokrane is the following:
\begin{conj}[{\bf\cite{M} Conjecture 3.27}]\label{conj6-1}
The filtration on $D^i$ induced by the weight spectral sequence {\rm(}loc.\ cit.\ {\rm3.25}, cf.\ \cite{Na} {\rm(2.0.9;$p$))}
agrees with the monodromy filtration induced by the monodromy operator $N:D^i \to D^i(-1)$ {\rm(}\cite{M} {\rm3.26)}.
\end{conj}
\begin{thm}\label{thm6-2}
Let $i \ge 0$ be an integer, and assume Conjecture {\rm\ref{conj6-1}} for $D^i$. Then the following sequence of $F$-isocrystals induced by \eqref{eq6-7} is exact{\rm:}
\[\xymatrix{ C^i \ar[r] & D^i \ar[r]^-{N} & D^i(-1) \ar[r] & E^i(-1). }\]
\end{thm}
\begin{pf}
Suppose that $\#(k) = p^a$. Then $\varphi^a$ is $K_0$-linear and we use the notion of weights concerning $\varphi^a$.
Because the groups
\[ \bH^i(Y,\cohwitt Y n \bullet), \quad \bH^i(Y,\mwitt Y n \bullet),\quad \bH^i(Y,\homwitt Y n \bullet),
 \quad \bH^i(Y,\mtwitt Y n \bullet) \]
 are finite for any $n \ge 1$, there is a diagram of $F$-isocrystals with exact rows by \eqref{eq6-5} and Proposition \ref{prop6-3}
\[ \xymatrix{ & \tD{}^i \ar@{=}[d] \ar[r] & D^i \ar[r]^-{N} & D^i(-1) \ar[r]\ar@{..>}[rrd] & \tD{}^{i+1} \ar@{=}[d] & \\
  C^i \ar[r]\ar@{..>}[rru] & \tD{}^i \ar[r]^-{\partial} & E^{i-1}(-1) \ar[r] & C^{i+1} \ar[r] & \tD{}^{i+1} \ar[r]^-{\partial} & E^i(-1),} \]
where we put
\[ \tD^i:=\qp \otimes_{\zp} \varprojlim_{n \ge 1} \ \bH^i(Y,\mtwitt Y n \bullet).  \]
By the assumption on $D^i$, $\ker(N)$ has weights $\le i$ and $\Coker(N)$ has weights $\ge i+2$. Hence for the reason of weights, it is enough to show the following lemma, where we do not need Conjecture \ref{conj6-1}:
\begin{lem}\label{lem6-4}
For integers $i \ge 0$, $C^i$ has weights $\le i$ and $E^i$ has weights $\ge i$.
\end{lem}
\noindent
{\it Proof of Lemma \ref{lem6-4}.}
When $Y$ is smooth, the lemma is a consequence of the Katz-Messing theorem \cite{KM}.
In the general case, $C^i$ has weights $\le i$ by the spectral sequence of $F$-isocrystals obtained from Proposition \ref{prop6-1}\,(1)
\[ E_1^{s,t}=\bigoplus_{|I|=s+1} \ D^t_I \Lra C^{s+t}, \]
where we put
\[ D^t_I :=\qp \otimes_{\zp} \varprojlim_{n \ge 1} \ \bH^t(Y_I,\witt {Y_I} n \bullet).  \]
Moreover $E^i$ has weights $\ge i$ by the duality result in Theorem \ref{thm6-1}.
Thus we obtain Lemma \ref{lem6-4} and Theorem \ref{thm6-2}.
\end{pf}
\begin{rem}
One can show Lemma {\rm\ref{lem6-4}} for $E^i$ alternatively by the Katz-Messing theorem and the spectral sequence obtained from Proposition {\rm\ref{lem6-1}\,(1)}
\[ E_1^{s,t}=\bigoplus_{|I|=1-s} \ D^{2s+t}_I(s) \Lra E^{s+t}. \]
\end{rem}
\medskip
\section{\bf Comparison with the finite part of Bloch-Kato}\label{sect8}
\medskip
In this section we relate the cohomology of $p$-adic \'etale Tate twists with the finite part of Bloch-Kato \cite{BK2}.
Let $p$ be a prime number, and let $K$ be a $p$-adic local field, i.e., a finite extension of $\qp$. Let $O_K$ be the integer ring of $K$, and let $k$ be the residue field of $O_K$. In this section we consider Tate twists only in the sense of $G_K$-$\zp$-modules (and not in the sense of $F$-crystals). For a topological $G_K$-$\zp$-module $M$ and an integer $r \in \bZ$, we define
\[ M(r) := \begin{cases} M \otimes_\zp \zp(1)^{\otimes r} \quad &\hbox{(if $r \ge 0$),} \\
\Hom_\zp(\zp(1)^{\otimes (-r)},M) \quad &\hbox{(if $r < 0$),} \end{cases}  \]
where $\zp(1)$ denotes the topological $G_K$-module $\varprojlim{}_{n\ge 1} \ \mu_{p^n}(\ol K)$.
\par
Let $X$ be a regular scheme which is projective flat over $O_K$ with strict semistable reduction, and put $\Xgf :=X \otimes_{O_K} K$ and $Y:=\cX \otimes_{O_K} k$.  Let $j : \Xgf \hra X$ be the natural open immersion, and let $i$ and $r$ be non-negative integers. Put $\cV^i:=H^i(\Xggf,\qp)$ (see also Notation) and
\[ H^{i+1}(X,\fT_{\qp}(r)) := \qp \otimes_{\zp} \varprojlim_{n \ge 1} \ H^{i+1}(X,\fT_n(r)_X), \]
which are finite-dimensional over $\qp$. Let $F^\bullet$ be the filtration on $H^{i+1}(\Xgf,\qp(r))$ resulting from the Hochschild-Serre spectral sequence
\begin{equation}\label{eq8-hs}
 E_2^{a,b}=H^a(K,\cV^b(r)) \Lra H^{a+b}(\Xgf,\qp(r)).
\end{equation}
Let $F^q H^{i+1}(X,\fT_{\qp}(r))$ be the inverse image of $F^q H^{i+1}(\Xgf,\qp(r))$ under the canonical map $j^* : H^{i+1}(X,\fT_{\qp}(r)) \to H^{i+1}(\Xgf, \qp(r))$. We have
\[ \varPhi^{i,r} := \gr_F^1 H^{i+1}(X,\fT_{\qp}(r)) \hra \gr_F^1 H^{i+1}(\Xgf, \qp(r)) = H^1(K,\cV^i(r)) \]
by definition. We relate $\varPhi^{i,r}$ with the finite part of Bloch-Kato \cite{BK2}:
\[  H^1_f(K,\cV^i(r)) := \ker(\iota : H^1(K,\cV^i(r)) \to H^1(K,\cV^i\otimes_{\qp} B_{\crys})). \]
where $B_\crys$ denotes the period ring of crystalline representations defined by Fontaine \cite{Fo1}, and $\iota$ is induced by the natural inclusion $\qp(r) \hra B_\crys$. The main result of this section is as follows
\begin{thm}\label{thm8-2}
Assume $0 \le r \le p-2$ and Conjecture {\rm\ref{conj6-1}} for $D^i:=D^i(Y)$. Then we have \[ \varPhi^{i,r} \subset H^1_f(K,\cV^i(r)). \] Further if $r \ge \dim(\Xgf)-p+3$, then we have $\varPhi^{i,r} = H^1_f(K,\cV^i(r))$.
\end{thm}
\begin{rem}
When $X$ is smooth over $O_K$, Theorem {\rm\ref{thm8-2}} is a part of the $p$-adic point conjecture raised by Schneider \cite{Sch} and proved by Langer-Saito \cite{LS} and Nekov\'a\v{r} \cite{N1}.
\end{rem}
We first prepare a key commutative diagram \eqref{eq8-01} below to prove Theorem \ref{thm8-2}. Put
\begin{align*}
H^{i+1}(X,\tau_{\le r}Rj_*\qp(r)) &:= \qp \otimes_{\zp} \varprojlim_{n \ge 1} \ H^{i+1}(X,\tau_{\le r}Rj_*\mu_{p^n}^{\otimes r})
\end{align*}
and let $F^q H^{i+1}(X,\tau_{\le r}Rj_*\mu_{p^n}^{\otimes r})$ be the inverse image of $F^q H^{i+1}(\Xgf,\qp(r))$ under the map
\[ j^* : H^{i+1}(X,\tau_{\le r}Rj_*\mu_{p^n}^{\otimes r}) \to H^{i+1}(\Xgf, \qp(r)). \]
To prove Theorem \ref{thm8-2}, we construct the following commutative diagram for $0 \le r \le p-2$:
\begin{equation}\label{eq8-01}
\xymatrix{  H^1(K,\cV^i(r))\ar[d]_\iota & F^1H^{i+1}(X,\tau_{\le r}Rj_*\qp(r)) \ar[rd]^-{\sigma'} \ar[l]  \ar[d]_\alpha \ar@{}[ld]|-{\text{(A)}} & \\
H^1(K,\cV^i\otimes_{\qp} B_\crys) & \ar[l]^-\delta D^i/ND^i \ar[r]_-\ep \ar@{}@<25pt>[r]_{\text{(B)}\qquad\quad\;\;}  &  E^i, }
\end{equation}
where $E^i:=E^i(Y)$ is as in \eqref{eq8-E}, and $D^i/ND^i$ denotes the cokernel of the monodromy operator $N:D^i \to D^i$. The top arrow is induced by $j^*$ and an edge homomorphism of the spectral sequence \eqref{eq8-hs}. The arrow $\ep$ is induced by the complex \eqref{eq6-7}, and the arrow $\sigma'$ is induced by the composite map
\[ H^{i+1}(X,\tau_{\le r}Rj_*\mu_{p^n}^{\otimes r}) \os{\sigma}\lra H^{i+1-r}(Y,\nu_{Y,n}^{r-1}) \lra \bH^i(Y,\homwitt Y n \bullet) \quad\hbox{($n \ge 1$).} \]
See Theorem \ref{thm2-2} for $\sigma$ and see Remark \ref{rem6-1} for the last map. The constructions of the arrows $\alpha$ and $\delta$ and the commutativity of the square (A) deeply rely on the $p$-adic Hodge theory. The commutativity of the triangle (B) is rather straight-forward once we define $\alpha$ (see Proposition \ref{prop8-2} below). The proof of Theorem \ref{thm8-2} will be given after Proposition \ref{prop8-2}.
\begin{rem}\label{rem8-1}
In his paper \cite{L}, Langer considered the diagram {\rm(A)} assuming that $O_K$ is absolutely unramified {\rm(}loc.\ cit.\ Proposition {\rm2.9)}. We remove this assumption and the assumption $(*)$ stated in loc.\ cit.\ p.\ {\rm191}, using a continuous version of crystalline cohomology {\rm(}see Appendix {\rm A} below{\rm)}. For this purpose, we include a detailed construction of $\alpha$ and a proof of the commutativity of {\rm(A)}.
\end{rem}
\def\logA{(X,M)}
\def\logAn{(X_n,M_n)}
\def\logB{(\ol X,\ol M)}
\def\logBn{(\ol X_n,\ol M_n)}
\def\logC{(\cE,M_\cE)}
\def\logCn{(\cE_n,M_{\cE_n})}
\def\logD{(W,M_W)}
\def\logDn{(\Wn,M_{\Wn})}
\def\logE{(S,M_S)}
\def\logEn{(S_n,M_{S_n})}
\def\logF{(V,M_V)}
\def\logFn{(V_n,M_{V_n})}
\def\logY{(Y,M_Y)}
\def\logZ{(Z,M_Z)}
\def\logZn{(Z_n,M_{Z_n})}
\def\cSb{\cS_\bullet}
\def\cSbt{\cS_\bullet\,\wt{}\,}
\def\cSbp{\cS_\bullet'}
\def\cSn{\cS_n}
\def\cSnt{\cS_n\hspace{1pt}\wt{}\,}
\def\cSnp{\cS_n'}
\def\cSqpt{\cS_\qp\hspace{-5pt}\wt{}\;\,}
\par\medskip
\noindent
{\bf Construction of $\bs{\delta}$.}
We first construct the map $\delta$. There is an exact sequence of topological $G_K$-$\qp$-modules
\begin{equation}\label{eq8-1} 0 \lra \cV^i \otimes_\qp B_\crys \lra \cV^i \otimes_\qp B_\st \os{1\otimes N}\lra \cV^i \otimes_\qp B_\st \lra 0, \end{equation}
where $N$ denotes the monodromy operator on $B_\st$. Taking continuous Galois cohomology groups, we obtain a long exact sequence
\begin{align*}
0 & \lra (\cV^i \otimes_\qp B_\crys)^{G_K} \lra (\cV^i \otimes_\qp B_\st)^{G_K} \os{1\otimes N}\lra (\cV^i \otimes_\qp B_\st)^{G_K} \\
& \lra H^1(G_K,\cV^i \otimes_\qp B_\crys) \lra \dotsb. \end{align*}
By the $C_\st$-conjecture proved by Hyodo, Kato and Tsuji (\cite{HK}, \cite{K4}, \cite{T1}, cf.\ \cite{Ni2}), the last arrow of the top row is identified with the map $N:D^i \to D^i$ (see \eqref{eq8-03} below for the construction of $C_\st$-isomorphisms). We thus define the desired map $\delta$ as the connecting map of this sequence.
\par\medskip\noindent
{\bf $\bs{C_\st}$-isomorphism.}
Before constructing $\alpha$, we give a brief review of the $C_\st$-isomorphisms, which will be useful later (\cite{K4}, \cite{T1} \S4.10). See \cite{K3} for the general framework of log structures and log schemes. Put
\[ W:=W(k) \quad \hbox{ and } \quad K_0:=\Frac(W)=W[p^{-1}]. \]
For a log scheme $\logZ$ and an integer $n \ge 1$, we define
\[ \logZn := \logZ \times_{\Spec(\bZ)} \Spec(\bZ/p^n), \]
where we regarded $\Spec(\bZ)$ and $\Spec(\bZ/p^n)$ as log schemes by endowing them with the trivial log structures and the fiber product is taken in the category of log schemes. Let $M$ be the log structure on $X$ associated with the normal crossing divisor $Y$. Let $\ol {O_K}$ be the integral closure of $O_K$ in $\ol K$, and let $\ol M$ be the log structure on $\ol X:=X \otimes_{O_K} \ol {O_K}$ define by base-change (see \eqref{eqA-01} below). Put $\ol Y := Y \otimes_k \ol k$. For $n \ge 1$ and $s \ge 0$, let $\cSnt(s)_{\logB} \in D^+(\ol Y_\et,\bZ/p^n)$ be Tsuji's version of syntomic complex of $\logB$ (\S\ref{sectA-5pre}).
Let $\ol j: \Xggf \hra \ol X$ be the natural open immersion. For $m,s \ge 0$, we define the syntomic cohomology of $\logB$ as
\[ H^m_\syn(\logB,\cSqpt(s)):=\qp \otimes_{\zp} \varprojlim_{n \ge 1} \ H^m\big(\ol Y,\cSnt(s)_{\logB}\big). \]
The following isomorphism due to Tsuji (\cite{T1} Theorem 3.3.4\,(2)) plays a crucial role:
\begin{equation}\label{eq8-02} H^m_\syn(\logB,\cSqpt(s)) \simeq H^m(\ol X,\tau_{\le s}R\ol j_*\qp(s)). \end{equation}
The right hand side is isomorphic to $\cV^m(s)=H^m(\Xggf,\qp(s))$ if $s \ge \dim(\Xgf)$ or $s \ge m$. We now introduce the crystalline cohomology of $\logB$ over $W$:
\[ H^m_\crys(\logB/W)_\qp:=\qp \otimes_{\zp} \varprojlim_{n \ge 1} \ H^m_\crys(\logBn/\Wn), \]
where $H^m_\crys(\logBn/\Wn)$ denotes the crystal cohomology of $\logBn$ over $\Wn$ (see \S\ref{sectA-1} below).
We define a $G_K$-equivariant homomorphism
\[ g^m : \cV^m  \lra D^m \otimes_{K_0} B_\st \]
as the composite map ($d:=\dim(\Xgf)$)
{\allowdisplaybreaks
\begin{align}
\label{eq8-03} \cV^m &\hspace{-2.8pt}\os{\text{\;(1)}}\isom H^m_\syn(\logB,\cSqpt(d))\otimes_\qp \qp(-d) \\
\notag &\os{\text{(2)}}\lra H^m_\crys(\logB/W)_\qp \otimes_\qp \qp(-d) \\
\notag &\hspace{-2.8pt}\os{\text{\;(3)}}\isom (D^m \otimes_{K_0} B_\st^+)^{N=0} \otimes_\qp \qp(-d) \\
\notag &\os{\text{(4)}}\lra D^m \otimes_{K_0} B_\st\,,
\end{align}
where (1) is obtained from the isomorphism \eqref{eq8-02} with $s=d$, the arrow (2) is the $p^d$-times of the canonical map induced by \eqref{eqA-cant} below, and (3) is the crystalline interpretation of $(D^m \otimes_{K_0} B_\st^+)^{N=0}$ due to Kato (\cite{K4} \S4, cf.\ Proposition \ref{propA-1}\,(3) below). The map (4) is induced by the natural inclusion $\qp(-d) \hra B_\st$. The $B_\st$-linear extension of $g^m$ is an isomorphism by \cite{T2} Theorem 4.10.2:
\begin{equation}\label{eq8-04} \cV^m \otimes_\qp B_\st \isom D^m \otimes_{K_0} B_\st\,, \end{equation}
which we call the $C_\st$-isomorphism. Now we have a commutative diagram of $G_K$-modules
\begin{equation}\label{eq8-cst}
\xymatrix{ H^m_\syn(\logB,\cSqpt(s)) \ar[r]_-{{\rm(1)}} \ar[d]_{\eqref{eq8-02}} & H^m_\crys(\logB/W)_\qp \ar[r]_-{{\rm(2)}} \ar@{}@<-2pt>[r]^-\sim & (D^m \otimes_{K_0} B_\st^+)^{N=0} \ar@<-15pt>@{^{(}->}[d] \\
\cV^m(s) \;\ar@{^{(}->}[r]^-{\qp(s) \hra B_\st} & \cV^m \otimes_{\qp} B_\st \ar[r]_-{\eqref{eq8-04}} \ar@{}@<-2pt>[r]^-\sim & D^m \otimes_{K_0} B_\st \;\;\quad}
\end{equation}
for $m,s \ge 0$, where the arrow (1) is the $p^r$-times of the canonical map \eqref{eqA-cant}, and (2) is the crystalline interpretation of $(D^m \otimes_{K_0} B_\st^+)^{N=0}$ mentioned in \eqref{eq8-03}. This diagram commutes by the construction of \eqref{eq8-04} and the compatibility facts in \cite{T1} Corollaries 4.8.8 and 4.9.2.
\begin{rem}
We modified the $C_\st$-isomorphisms in \cite{T1} by constant multiplications to make them consistent with the canonical map from syntomic cohomology of Kato to {\rm(}continuous{\rm)} crystalline cohomology {\rm(}see the map {\rm(1)} of \eqref{eq8-2} and the equality \eqref{eqA-trans} below{\rm)}. Under this modification, the element $t_\syn$ in \cite{T1} p.\ {\rm378} is replaced with $p^{-1}t_\syn$, and Corollary {\rm4.8.8} of loc.\ cit.\ remains true.
\end{rem}
\par\medskip
\noindent
{\bf Construction of $\bs{\alpha}$.}
We start the construction of $\alpha$ assuming $0 \le r \le p-2$, which will be complete in Corollary \ref{cor8-1} below. For $n \ge 1$, let $\cS_n(r)_{\logA} \in D^+(Y_\et,\bZ/p^n)$ be the syntomic complex of Kato (see \S\ref{sectA-5pre} below). Since $0 \le r \le p-2$, there is an isomorphism
\begin{equation}\label{eq8-0}  (\tau_{\le r}Rj_*\mu_{p^n})|_Y \simeq \cS_n(r)_{\logA} \quad \hbox{ in } \;D^+(Y_\et,\bZ/p^n) \end{equation}
for each $n \ge 1$ (\cite{T2} Theorem 5.1), and $H^{i+1}(X,\tau_{\le r}Rj_*\qp(r))$ is isomorphic to the syntomic cohomology
\[ H^{i+1}_\syn(\logA,\cS_\qp(r)):=\qp \otimes_{\zp} \varprojlim_{n \ge 1} \ H^{i+1}(Y,\cS_n(r)_{\logA}) \]
by the proper base-change theorem. Therefore we consider $H^{i+1}_\syn(\logA,\cS_\qp(r))$ instead of $H^{i+1}(X,\tau_{\le r}Rj_*\qp(r))$, which is the major reason we assume $r \le p-2$. We first prove:
\begin{lem}\label{lem8-1}
The kernel of the composite map
\[ \xi : H^m_\syn(\logA,\cS_\qp(r)) \lra H^m_\crys(\logA/W)_\qp \lra H^m_\crys(\logB/W)_\qp \]
agrees with that of the composite map
\[ \eta : H^m_\syn(\logA,\cS_\qp(r)) \os{\eqref{eq8-0}}\simeq H^m(X,\tau_{\le r}Rj_*\qp(r)) \lra \cV^m(r)=H^m(\Xggf,\qp(r)). \]
\end{lem}
\begin{pf}
The map $\eta$ factor through \eqref{eq8-02} with $s=r$, and the map $\xi$ factors through the map (1) in \eqref{eq8-cst} with $s=r$ by \eqref{eqA-trans} below. Hence the diagram
\[\xymatrix{ H^m_\syn(\logA,\cS_\qp(r)) \ar[r]^{\xi} \ar[d]_\eta & H^m_\crys(\logB/W)_\qp \ar@<-5pt>@{^{(}->}[d] \\
\cV^m(r) \;\ar@{^{(}->}[r] & D^m \otimes_{\qp} B_\st }\]
commutes by \eqref{eq8-cst}, which implies the assertion.
\end{pf}
\noindent
Let $H^m_\cocrys(\logA/W)$ be the continuous crystalline cohomology of $\logA$ over $W$ (see \S\ref{sectA-2} below). Put
\begin{align*}
H^m_\syn(\logA,\cS_\qp(r))^0 &:= \ker(\eta : H^m_\syn(\logA,\cS_\qp(r)) \lra \cV^m(r)), \\
H^m_\cocrys(\logA/W)_\qp & :=\qp \otimes_{\zp} H^m_\cocrys(\logA/W), \\
H^m_\cocrys(\logA/W)_\qp^0 &:= \ker(H^m_\cocrys(\logA/W)_\qp \to H^m_\crys(\logB/W)_\qp).
\end{align*}
There are canonical maps
{\allowdisplaybreaks
\begin{align}
\label{eq8-2} H^{i+1}_\syn(\logA,\cS_\qp(r))^0 &\os{\text{(1)}}\lra H^{i+1}_\cocrys(\logA/W)_\qp^0 \\
\notag &\os{\text{(2)}}\lra H^1(K,H^i_\crys(\logB/W)_\qp) \\
\notag &\hspace{-2.8pt}\os{\text{\;(3)}}\isom H^1(K,(D^i \otimes_{K_0} B_\st^+)^{N=0}) \\
\notag &\os{\text{(4)}}\lra H^1(K,\cV^i \otimes_\qp B_\crys),
\end{align}
}where (1) is the canonical map obtained from Lemma \ref{lem8-1} and Proposition \ref{propA-5}\,(2) below, the arrow (3) is the isomorphism (2) in \eqref{eq8-cst} for $m=i$, and (4) is obtained from the $C_\st$-isomorphism \eqref{eq8-04} and \eqref{eq8-1}. The most crucial map (2) is obtained from an edge homomorphism of the spectral sequence in Theorem \ref{thmA-1} (see also Corollary \ref{corA-2} below). To proceed with the construction of $\alpha$, we need to introduce some notation. 
\begin{defn}\label{def8-1}
Put $S:=\Spec(O_K)$, and fix a prime $\pi$ of $O_K$. Let $M_S$ be the log structure on $S$ associated with the closed point $\{\pi=0\}$. Put $V:=\Spec(W[t])$, and let $M_V$ be the log structure on $V$ associated with the divisor $\{t=0\}$.
Let \[ \iota : \logE \lra \logF \] be the exact closed immersion over $W$ sending $t \mapsto \pi$. For each $n \ge 1$, let $\logCn$ be the PD-envelope of the exact closed immersion $\iota_n : \logFn \hra \logEn$ compatible with the canonical PD-structure on $p\Wn$. We define
\[ H^m_\crys(\logA/\logC)_\qp:=\qp \otimes_{\zp} \varprojlim_{n \ge 1} \ H^m_\crys(\logAn/\logCn), \]
which is an $R_\cE:=\varprojlim{}_{n \ge 1} \ \vG(\cE_n,\cO_{\cE_n})$-module endowed with a $\varphi_\cE$-semilinear endomorphism $\varphi$ and a $K_0$-linear endomorphism $N$ satisfying $N\varphi=p\varphi N$. Here $\varphi_\cE$ is a certain canonical Frobenius operator on $R_\cE$. See \cite{T1} p.\ {\rm253} for $\varphi_\cE$, loc.\ cit.\ \S{\rm4.3} for $N$ and $\varphi$, and see loc.\ cit.\ Remark {\rm4.3.9} for $N$.
\end{defn}
\noindent
The following lemma plays a key role in our construction of $\alpha$, where we do not need to assume $0 \le r \le p-2$ (compare with the assumption $(*)$ in \cite{L} p.\ {\rm191}).
\begin{lem}\label{lem8-2}
\begin{enumerate}
\item[(1)]
There exists a long exact sequence
{\begin{align*}
 \dotsb & \to H^i_\cocrys(\logA/W)_\qp \to H^i_\crys(\logA/\logC)_\qp \os{N}\to H^i_\crys(\logA/\logC)_\qp \\
 & \os{\partial}\to H^{i+1}_\cocrys(\logA/W)_\qp \to \dotsb 
\end{align*}} and the image of $\partial$ agrees with $H^{i+1}_\cocrys(\logA/W)_\qp^0$.
\item[(2)]
There is a commutative diagram with exact lower row
\[\xymatrix{ & H^i_\crys(\logA/\logC)_\qp \ar@{->>}[r]^-\partial \ar[d] & H^{i+1}_\cocrys(\logA/W)_\qp^0 \ar[d] \\
D^i \ar[r]^-N & D^i \ar[r]^-{(-1)^{i+1}\delta} & H^1(K,\cV^i \otimes_\qp B_\crys), }\]
where the right vertical arrow is the composite of {\rm(2)\,--\,(4)} in \eqref{eq8-2}, and the left vertical arrow is the specialization maps with $t \mapsto 0$ {\rm(}cf.\ \cite{T1} {\rm(4.4.4))}.
\end{enumerate}
\end{lem}
\begin{pf}
(1) follows from Proposition \ref{propA-20} and Corollaries \ref{corA-1}\,(1) and \ref{propA-3} below.  As for (2), the exactness of the lower row has been explained in the construction of $\delta$. See Corollary \ref{corA-2} below for the commutativity of the square.
\end{pf}
\begin{cor}\label{cor8-1}
The composite map \eqref{eq8-2} factors as follows{\rm:}
\[ H^{i+1}_\syn(\logA,\cS_\qp(r))^0 \lra D^i/ND^i \os{\delta}\hra H^1(K,\cV^i \otimes_\qp B_\crys). \]
We define $\alpha$ as the first arrow of this decomposition.
\end{cor}
\par\medskip
\noindent
{\bf Commutativity of the diagram (\ref{eq8-01}).}
We first prove that the square (A) in \eqref{eq8-01} is commutative. By the construction of $\alpha$ and $\delta$, it is enough to show
\begin{thm}\label{prop8-1}
The following square is commutative{\rm:}
\begin{equation}\label{eq8-3}
\xymatrix{H^{i+1}_\syn(\logA,\cS_\qp(r))^0 \ar[d]_{\eqref{eq8-2}} \ar[r]\ar@{}@<-2pt>[r]^-\sim & F^1H^{i+1}(X,\tau_{\le r}Rj_*\qp(r)) \ar[d] \\
 H^1(K,\cV^i \otimes_\qp B_\crys) & \ar[l]_-{B_\crys \hookleftarrow \qp(r)} H^1(K,\cV^i(r)),}
\end{equation}
where the right vertical arrow is the top arrow in \eqref{eq8-01}.
\end{thm}
\begin{pf}
Put $d:=\dim(\Xgf)$. We define a $G_K$-homomorphism
\[ \beta^{i,r} : \cV^i(r) \lra \begin{cases} H^i_\crys(\logB/W)_\qp (r-d) \quad & \hbox{(if $r < d$)} \\
 H^i_\crys(\logB/W)_\qp  \quad & \hbox{(if $r \ge d$)} \end{cases} \]
as the following composite map if $r \ge d$:
\[\xymatrix{ \beta^{i,r} :  \cV^i(r) & \ar[l]_-{\eqref{eq8-02}}^-\sim H^i_\syn(\logB,\cSqpt(r)) \ar[rr]^-{\text{\eqref{eq8-cst}\,(1)}} && H^i_\crys(\logB/W)_\qp, }\]
and as the following map if $r < d$:
\[\xymatrix{ \beta^{i,r} :  \cV^i(r) \simeq \cV^i(d) \otimes_\qp \qp(r-d) \ar[rr]^-{\beta^{i,d} \otimes \id} && H^i_\crys(\logB/W)_\qp\otimes_\qp \qp(r-d). }\]
We next define a homomorphism
\[ \gamma^{i,r} : H^{i+1}_\syn(\logA,\cS_\qp(r))^0 \lra \begin{cases} H^1(K,H^i_\crys(\logB/W)_\qp (r-d)) \;\; & \hbox{(if $r < d$)} \\
 H^1(K,H^i_\crys(\logB/W)_\qp)  \quad & \hbox{(if $r \ge d$)} \end{cases}   \]
as the following composite map if $r < d$:
\begin{align*}
H^{i+1}_\syn(\logA,\cS_\qp(r))^0 & \os{\text{(i)}}\lra H^1(K,H^i_\crys(\logB/W)_\qp) \\
& \,\;\simeq\;\hspace{1.1pt} H^1(K,H^i_\crys(\logB/W)_\qp\otimes_\qp\qp(d-r)\otimes\qp(r-d)) \\
& \os{\text{(ii)}}\lra H^1(K,H^i_\crys(\logB/W)_\qp(r-d)),
\end{align*}
where (i) is the composite of (1) and (2) in \eqref{eq8-2}, and (ii) is induced by the inclusion $\qp(d-r) \hra H^0_\crys(\logB/W)_\qp$ (cf.\ \eqref{eqA-incl} below) and the cup product of crystalline cohomology. If $r \ge d$, we define $\gamma^{i,r}$ as the composite of (1) and (2) in \eqref{eq8-2}.

We first prove the commutativity of \eqref{eq8-3} assuming $r <d$. In this case, the following triangle commutes by Theorem \ref{propA-4} below:
\begin{equation}\label{eq8-4}
\xymatrix{H^{i+1}_\syn(\logA,\cS_\qp(r))^0 \ar[d] \ar[rd]^{\gamma^{i,r}} \\
 H^1(K,\cV^i(r)) \ar[r]^-{\beta^{i,r}} & H^1(K,H^i_\crys(\logB/W)_\qp (r-d)),}
\end{equation}
where the left vertical arrow is the composite of the top horizontal and the right vertical arrows of \eqref{eq8-3}. Moreover the following square commutes by the construction of the $C_\st$-isomorphism \eqref{eq8-04} (compare with \eqref{eq8-cst}):
\[ \xymatrix{\cV^i(r) \ar[d]_{\qp(r) \hra B_\crys} \ar[rr]^-{\beta^{i,r}} && H^i_\crys(\logB/W)_\qp (r-d) \ar[d] \\
 \cV^i \otimes_\qp B_\crys \ar[rr]^-{\eqref{eq8-1},\;\eqref{eq8-04}}_-\sim && (D^i \otimes_\qp B_\st)^{N=0},}\]
where the right vertical arrow is given by the isomorphism (2) of \eqref{eq8-cst} and the inclusion $\qp(r-d) \hra B_\crys$.  The commutativity of \eqref{eq8-3} follows from these commutative diagrams and the fact that the composite map \eqref{eq8-2} agrees with the composite of $\gamma^{i,r}$ and the map
\[ H^1(K,H^i_\crys(\logB/W)_\qp (r-d)) \lra H^1(K,\cV^i \otimes_\qp B_\crys) \]
induced by the right vertical arrow and the bottom isomorphism in the previous diagram.

When $r \ge d$, the following triangle commutes by Corollary \ref{thmA-2} below (see also Proposition \ref{propA-5}\,(2)):
\[ \xymatrix{H^{i+1}_\syn(\logA,\cS_\qp(r))^0 \ar[d] \ar[rd]^{\gamma^{i,r}} \\
 H^1(K,\cV^i(r)) \ar[r]^-{\beta^{i,r}} & H^1(K,H^i_\crys(\logB/W)_\qp),} \]
where the left vertical arrow is defined in the same way as for that of \eqref{eq8-4}. The commutativity of \eqref{eq8-3} follows from this diagram and the diagram \eqref{eq8-cst}. This completes the proof of Theorem \ref{prop8-1}.
\end{pf}
\begin{prop}\label{prop8-2}
The triangle {\rm(B)} in \eqref{eq8-01} is commutative.
\end{prop}
\begin{pf}
There is a commutative diagram by the definition of $\partial$ in Lemma \ref{lem8-2}\,(1) (cf.\ Proposition \ref{propA-20})
\[\xymatrix{ H^i_\crys(\logA/\logC)_\qp \ar[r]^-{t \mapsto 0} \ar[d]_{\partial} & D^i\phantom{|} \ar@<-2pt>[d]^{\frac{dt}{t} \wedge } \\
H^{i+1}_\cocrys(\logA/W)_\qp \ar[r]^-{t \mapsto 0} & \wt{D}^{i+1},\phantom{|}}\]
where $\wt{D}^{i+1}=\wt{D}^{i+1}(Y)$ is as in the proof of Theorem \ref{thm6-2}, the left (resp.\ right) vertical arrow is as in Lemma \ref{lem8-2}\,(1) (resp.\ given by the product with $\dlog(t)$, cf.\ \eqref{eq6-5}), and the horizontal arrows are the specialization maps by $t \mapsto 0$. Therefore in view of the constructions of $\alpha$ and $\sigma'$, it is enough to check that the following diagram of canonical morphisms commutes in $D^+(Y_\et,\bZ/p^n)$ for each $n \ge 0$:
\[\xymatrix{ \cS_n(r)_{\logA} \ar[r] \ar[d]_\sigma & \bJ^{[r]}_{\logAn/\Wn} \ar[r]^-s
 & \mtwitt Y n {\bullet \ge r} \ar[r] \ar[d]_{\eqref{eq6-6}} & \mtwitt Y n {\bullet} \ar[d]^{\eqref{eq6-6}} \\
\nu_{Y,n}^{r-1}[-r] \ar[rr]^{\text{Remark\;\ref{rem6-1}}} && \homwitt Y n {\bullet \ge r-1}[-1] \ar[r] & \homwitt Y n {\bullet}[-1], }\]
where $\bJ^{[r]}_{\logAn/\Wn}$ is as in \S\ref{sectA-5pre} below, and $s$ denotes the specialization map by $t \mapsto 0$. The right square commutes obviously in the category of complexes of sheaves. As for the left square, we have only to show the commutativity on the diagram of the $r$-th cohomology sheaves because $\homwitt Y n {\bullet \ge r-1}[-1]$ is concentrated in degrees $\ge r$ and $\cS_n(r)_{\logA}$ is concentrated in degrees $\le r$ by \eqref{eq8-0}. Since we have
\[ \cH^r(\cS_n(r)_{\logA}) \os{\eqref{eq8-0}}\simeq (R^rj_*\mu_{p^n}^{\otimes r})|_Y \]
and the sheaf on the right hand side is generated by symbols (cf.\ the case $D=\emptyset$ of Theorem \ref{thm2-1}\,(1)), one can check this commutativity of the diagram of sheaves by a straight-forward computation on symbols using the compatibility in \cite{T1} Proposition 3.2.4\,(2).
\end{pf}
\par\medskip
\noindent
{\bf Proof of Theorem \ref{thm8-2}.} We have the following commutative diagram by \eqref{eq8-01}:
\[\xymatrix{ F^1H^{i+1}(X,\fT_{\qp}(r)) \ar[r] \ar[rd] & F^1H^{i+1}(X,\tau_{\le r}Rj_*\qp(r)) \ar[r]^-{\sigma'} \ar[d] \ar[rd]_\alpha \ar@{}@<-10pt>[r]_-{\text{(B)}} \ar@{}[rdd]|-{\text{(A)}} & E^i \\
 & H^1(K,\cV^i(r))\ar[rd]_\iota & \ar[d]^\delta \Coker(D^i \os{N}\to D^i) \phantom{\big|} \ar@{^{(}->}[u]_\ep \\
 & & \hspace{-20pt} H^1(K,\cV^i \otimes_{\qp} B_\crys), }\]
where the left triangle commutes obviously. The kernel of $\iota$ is $H^1_f(K,\cV^i(r))$ by definition. The top row is a complex by the definitions of $\sigma'$ and $\fT_n(r)_X$. The map $\ep$ is injective by Theorem \ref{thm6-2} and Conjecture \ref{conj6-1} for $D^i$. We obtain the first assertion
\begin{equation}\label{eq8-5} \varPhi^{i,r} \subset H^1_f(K,\cV^i(r)) \end{equation}
by a simple diagram chase on this diagram. To obtain the second assertion, we repeat the same arguments for $\varPhi^{2d-i,d-r+1}$, where we put $d:=\dim(\Xgf)$. The map
\[ \ep : \Coker(N : D^{2d-i} \to D^{2d-i}) \lra E^{2d-i} \]
is injective by Conjecture \ref{conj6-1} for $D^i$, Theorems \ref{thm6-2} and \ref{thm6-1} and the Poincar\'e duality between $D^i$ and $D^{2d-i}$ \cite{H2}. Hence we have
\begin{equation}\label{eq8-6} \varPhi^{2d-i,d-r+1} \subset H^1_f(K,\cV^{2d-i}(d-r+1)). \end{equation}
The results \eqref{eq8-5} and \eqref{eq8-6} imply $\varPhi^{i,r} = H^1_f(K,\cV^i(r))$, because $\varPhi^{i,r}$ and $\varPhi^{2d-i,d-r+1}$ (resp.\ $H^1_f(K,\cV^i(r))$ and $H^1_f(K,\cV^{2d-i}(d-r+1))$) are the exact annihilators of each other under the non-degenerate pairing of the Tate duality
\[ H^1(K,\cV^i(r)) \times H^1(K,\cV^{2d-i}(d-r+1)) \lra \qp \]
by \cite{S1} Corollary 10.6.1 (resp.\ \cite{BK2} Proposition 3.8).
\hfill \qed
\bigskip

\section{\bf Image of $\bs{p}$-adic regulators}\label{sect7}
\medskip
In this section, we prove Theorem \ref{thm1-1}. We first review the definitions of the finite part of Galois cohomology \cite{BK2} and the integral part of algebraic $K$-groups \cite{Scho} over number fields.
\par
Let $K$ be a number field and let $\cV$ be a topological $G_K$-$\qp$-module which is finite-dimensional over $\qp$. Let $\fo$ be the integer ring of $K$. For a place $v$ of $K$, let $K_v$ (resp.\ $\fo_v$) be the completion of $K$ (resp.\ $\fo$) at $v$. The finite part $H^1_f(K,\cV) \subset H^1(K,\cV)$ is defined as
\[  H^1_f(K,\cV):=\ker\bigg(H^1(K,\cV) \lra \prod_{v\,:\,\text{finite}} \ \frac{H^1(K_v,\cV)}{H^1_f(K_v,\cV)}\bigg), \]
where $v$ runs through all finite places of $K$, and $H^1_f(K_v,\cV)$ is defined as follows:
\[  H^1_f(K_v,\cV) := \begin{cases} \ker(H^1(K_v,\cV) \to H^1(I_v,\cV)) \quad & \hbox{(if $v \hspace{-6pt}\not| p$),} \\
\ker(H^1(K_v,\cV) \to H^1(K_v,\cV\otimes_{\qp} B_{\crys})) \quad & \hbox{(if $v | p$),} \end{cases} \]
where $I_v$ denotes the inertia subgroup of $G_{K_v}$. We next review the definition of the integral part of algebraic $K$-groups. Let $V$ be a proper smooth variety over the number field $K$. First fix a finite place $v$ of $K$. By de Jong's alteration theorem \cite{dJ}, there is a proper generically finite morphism
\[ \pi : V' \lra V_v:=V\otimes_K K_v \]
such that $V'$ has a projective regular model $X'$ with strict semistable reduction over the integer ring of some finite extension of $K_v$. The integral part $K_m(V_v)_{\fo_v} \subset K_m(V_v) \otimes \bQ$ is the kernel of the composite map
\[ K_m(V_v) \otimes \bQ \os{\pi^*}\lra K_m(V')\otimes \bQ \lra \frac{K_m(V')\otimes \bQ}{\hbox{Image of $K_m(X')\otimes \bQ$}}, \]
which is in fact independent of $X'$ (\cite{Scho} \S1). The integral part $K_m(V)_\fo$ is defined as
\[  K_m(V)_\fo :=\ker\bigg(K_m(V) \otimes \bQ \lra \prod_{v\,:\,\text{finite}} \ \frac{K_m(V_v) \otimes \bQ}{K_m(V_v)_{\fo_v}}\bigg), \]
where $v$ runs through all finite places of $K$. If $V$ admits a regular model which is proper flat over the integer ring of $K$, then $K_m(V)_\fo$ agrees with the image of the $K$-group of the model, i.e., the integral part considered by Beilinson.
By these definitions of $H^1_f(K,\cV)$ and $K_{2r-i-1}(V)_\fo$, Theorem \ref{thm1-1} is immediately reduced to Theorem \ref{thm9-1} below, which is an analogue of Theorem \ref{thm1-1} over local fields.

Let $\ell$ and $p$ be prime numbers. We change the setting slightly, and let $K$ be an $\ell$-adic local field, i.e., a finite extension of $\bQ_\ell$. Let $O_K$ be the integer ring of $K$, and let $V$ be a proper smooth variety over $K$. Let $i$ and $r$ be integers with $2r \ge i+1 \ge 1$. There is a $p$-adic regulator map obtained from \'etale Chern character
\[ \reg_p : K_{2r-i-1}(V)_{O_K,\,\hom} \lra H^1(K,H^i(\ol V, \qp(r))), \]
where $\ol V$ denotes $V \otimes_K \ol K$, and $K_{2r-i-1}(V)_{O_K,\,\hom}$ denotes the subspace of $K_{2r-i-1}(V)_{O_K}$ consisting of all elements which vanish in $H^{i+1}(\ol V, \qp(r))$ under the Chern character.
\begin{thm}\label{thm9-1}
Assume $2r>i+1$. If $\ell = p$, assume further that $r \le p-2$ and that Conjecture {\rm\ref{conj6-1}} holds for projective strict semistable varieties over $\bF_p$ in degree $i$. Then $\Image(\reg_p)$ is contained in $H^1_f(K,H^i(\ol V, \qp(r)))$.
\end{thm}
\begin{rem}
When $\ell=p$, we need Conjecture {\rm\ref{conj6-1}} for the reduction of an alteration of $V$. If $\ell \ne p$, we do not need the monodromy-weight conjecture, but use Deligne's proof of the Weil conjecture \cite{D} to show that $\reg_p$ is zero.
\end{rem}
\begin{pf*}{Proof of Theorem \ref{thm9-1}}
We first reduce the problem to the case that $V$ has a regular model which is projective flat over $O_K$ with strict semistable reduction. By de Jong's alteration theorem \cite{dJ}, there exists a proper generically finite morphism $\pi : V' \to V$ such that $V'$ has a projective regular model with strict semistable reduction over the integer ring $O_L$ of $L:=\vG(V',\cO_{V'})$. Then there is a commutative diagram
\[\xymatrix{
 K_{2r-i-1}(V)_{O_K,\,\hom} \ar[r]^-{\reg_p} \ar[d]_{\pi^*} & H^1(K,H^i(\ol V, \qp(r))) \ar[r] \ar[d]_{\pi^*} & \dfrac{H^1(K,H^i(\ol V, \qp(r)))}{H^1_f(K,H^i(\ol V, \qp(r)))} \ar@{^{(}->}[d]_{\pi^*} \\
 K_{2r-i-1}(V')_{O_L,\,\hom} \ar[r]^-{\reg_p} & H^1(L,H^i(\ol {V'}, \qp(r))) \ar[r] & \dfrac{H^1(L,H^i(\ol {V'}, \qp(r)))}{H^1_f(L,H^i(\ol {V'}, \qp(r)))}
}\] where $\ol {V'}$ denotes $V' \otimes_L \ol K$, and the right (and the middle) vertical arrows are split injective by a standard argument using a corestriction map of Galois cohomology and a trace map of \'etale cohomology. By this diagram and the definition of $K_*(V)_{O_K}$, Theorem \ref{thm9-1} for $V$ is reduced to that for $V'$.  Thus we may assume that $V$ has a projective regular model $X$ with strict semistable reduction over $O_K$. Then the case $\ell \ne p$ follows from \cite{N1} II Theorem 2.2 (cf.\ \cite{D}). We prove the case $\ell = p$. Let $Y$ be the reduction of $X$. Put
\[ H^{i+1}(X,\fT_{\qp}(r)) := \qp \otimes_{\zp} \varprojlim_{n \ge 1} \ H^{i+1}(X,\fT_n(r)), \]
and let $K_{2r-i-1}(X)_{\hom}$ (resp.\ $H^{i+1}(X,\fT_{\qp}(r))^0$) be the kernel of the composite map
\begin{align*}
& K_{2r-i-1}(X) \lra K_{2r-i-1}(V) \os{\ch}{\lra} H^{i+1}(V, \qp(r)) \to H^{i+1}(\ol V, \qp(r))\\
&  \hbox{(resp.\ } \;\; H^{i+1}(X,\fT_{\qp}(r)) \lra H^{i+1}(V, \qp(r)) \to H^{i+1}(\ol V, \qp(r)) \;\hbox{).}
\end{align*}
There is a commutative diagram by Corollary \ref{thm3-2}
\[\xymatrix{ K_{2r-i-1}(X)_{\hom} \ar@{->>}[r] \ar[d]_{\ch} & K_{2r-i-1}(V)_{O_K,\,\hom} \ar[d]^{\reg_p} \\
 H^{i+1}(X,\fT_{\qp}(r))^0 \ar[r] & H^1(K,H^i(\ol V, \qp(r))).}\]
The image of the bottom arrow is contained in $H^1_f(K,H^i(\ol V, \qp(r)))$ by Theorem \ref{thm8-2} and Conjecture \ref{conj6-1} for $D^i=D^i(Y)$, which implies Theorem \ref{thm9-1}. This completes the proof of Theorems \ref{thm9-1} and \ref{thm1-1}.
\end{pf*}
\medskip
\appendix
\section{\bf Continuous crystalline cohomology}\label{sectA}
\numberwithin{equation}{subsection}
\def\logAb{(X_\bullet,M_\bullet)}
\def\logAs{(X^\star,M^\star)}
\def\logAbs{(X_\bullet^\star,M_\bullet^\star)}
\def\logAns{(X_n^\star,M_n^\star)}
\def\logAL{(X_{O_L},M_{O_L})}
\def\logALs{(X_{O_L}^\star,M_{O_L}^\star)}
\def\logALbs{(X_{O_L,\bullet}^\star,M_{O_L,\bullet}^\star)}
\def\logALn{(X_{O_L,n},M_{O_L,n})}
\def\logBb{(\ol X_\bullet, \ol M_\bullet)}
\def\logBs{(\ol X{}^\star, \ol M{}^\star)}
\def\logBbs{(\ol X{}_\bullet^\star, \ol M{}_\bullet^\star)}
\def\logCb{(\cE_\bullet,M_{\cE_\bullet})}
\def\logEL{(S_L,M_{S_L})}
\def\logEb{(\ol S,\ol {M_S})}
\def\logP{(\cD,M_\cD)}
\def\logPn{(\cD_n,M_{\cD_n})}
\def\logPns{(\cD_n^\star,M_{\cD_n^\star})}
\def\logZno{(Z_{n+1},M_{Z_{n+1}})}
\def\logZi{(Z^i,M_{Z^i})}
\def\logZs{(Z^\star,M_{Z^\star})}
\def\logZLs{(Z_L^\star,M_{Z_L^\star})}
\def\logZLps{(Z_{L'}^\star,M_{Z_{L'}^\star})}
\def\logZLi{(Z_L^i,M_{Z_L^i})}
\def\logZns{(Z_n^\star,M_{Z_n^\star})}
\def\logZnrs{(Z_{n+r}^\star,M_{Z_{n+r}^\star})}
\def\logT{(T,M_T)}
\def\logTo{(T_1,M_{T_1})}
\def\Wb{W_\bullet}
\def\bZpd{\bZ/p^\bullet}
\def\gkbZpd{G_K\text{-}\bZ/p^\bullet}
\def\Mod{\text{\it Mod}}
\def\qt{\bQ\otimes}
\def\Shv{\text{\it Shv}}
\def\cC{\mathscr C}
\def\cJ{\mathscr J}
\medskip

In this appendix, we formulate continuous versions of crystalline and syntomic cohomology of log schemes, combining the methods of Jannsen, Kato and Tsuji (\cite{J}, \cite{KaV}, \cite{K3}, \cite{T1}, \cite{T2}). The results of this appendix have been used in \S\ref{sect8} of this paper. See \cite{K3} for the general framework of log structures and log schemes.
\par
Let $p$ be a prime number, and let $K$ be a complete discrete valuation field of characteristic $0$ whose residue field $k$ is a perfect field of characteristic $p$. Let $O_K$ be the integer ring of $K$. Put $W:=W(k)$, $\Wn:=\Wn(k)$ ($n \ge 1$) and $K_0:=\Frac(W)$.
Let $X$ be a regular scheme which is projective flat over $O_K$ with semistable reduction, and put
\begin{align*}
    &\Xgf:=X \otimes_{O_K} K, \qquad Y:=X \otimes_{O_K} k \\
\Xggf := & \Xgf \otimes_K \ol K, \quad \ol Y := Y \otimes_k \ol k \quad \hbox{ and } \quad \ol X:=X \otimes_{O_K} \ol {O_K},
\end{align*}
where $\ol {O_K}$ denotes the integral closure of $O_K$ in $\ol K$.  Let $M$ be the log structure on $X$ associated with the normal crossing divisor $Y$. We define a log structure $\ol M$ on $\ol X$ as follows. For a finite field extension $L/K$, put $S_L:=\Spec(O_L)$, and let $M_{S_L}$ be the log structure on $S_L$ associated with its closed point. We denote $(S_K,M_{S_K})$ simply by $\logE$, and define $\logAL$ by base-change in the category of log schemes
\begin{equation}\label{eqA-02}
 \logAL:=\logA \times_{\logE} \logEL \,.
\end{equation}
We then define a log structure $\ol M$ on $\ol Y$ as that associated with the pre-log structure 
\begin{equation}\label{eqA-01}
 \varinjlim_{K \subset L \subset \ol K} \ M_{O_L}|_{\ol X},
\end{equation}
where $L$ runs through all finite field extensions of $K$ contained in $\ol K$, $M_{O_L}|_{\ol X}$ denotes the topological inverse image of $M_{O_L}$ onto $\ol X$, and the inductive limit is taken in the category of \'etale sheaves of monoids on $\ol X$.
\par
For a log scheme $\logZ$ and an integer $n \ge 1$, we define
\[ \logZn := \logZ \times_{\Spec(\bZ)} \Spec(\bZ/p^n) (=: \logZ \otimes \bZ/p^n) \]
where we regarded $\Spec(\bZ)$ and $\Spec(\bZ/p^n)$ as log schemes by endowing them with the trivial log structures and the fiber product is taken in the category of log schemes.
\par
We denote by $\Mod(\bZpd)$ the abelian category of the projective systems $\{F_n \}_{n \ge 1}$ of abelian groups such that $F_n$ is a $\bZ/p^n$-module for each $n \ge 1$. For a profinite group $G$, $\Mod(G\text{-}\bZpd)$ denotes the abelian category of the projective systems $\{F_n \}_{n \ge 1}$ of discrete $G$-modules such that $F_n$ is a $\bZ/p^n$-module.
\par
For a scheme $T$, $\Shv(T_\et,\bZpd)$ denotes the abelian category of the projective systems $\{ \cF_n \}_{n \ge 1}$ of \'etale sheaves on $T$ such that $\cF_n$ is a $\bZ/p^n$-sheaf for each $n \ge 1$. For a profinite group $G$ acting on $T$, $\Shv(T_\et,G\text{-}\bZpd)$ denotes the abelian category of the projective systems $\{ \cF_n \}_{n \ge 1}$ of \'etale $G$-sheaves on $T$ such that $\cF_n$ is a $\bZ/p^n$-sheaf for each $n \ge 1$. We write $D(T_\et,\bZpd)$ (resp.\ $D(T_\et,G\text{-}\bZpd)$) for the derived category of $\Shv(T_\et,\bZpd)$ (resp.\ $\Shv(T_\et,G\text{-}\bZpd)$).
\par
For an additive category $\cC$, we write $\bQ \otimes \cC$ for the $\bQ$-tensor category of $\cC$, i.e., the category whose objects are the same as $\cC$ and such that for objects $A,B \in \cC$, the group of morphisms $\Hom_{\bQ \otimes \cC}(A,B)$ is given by $\bQ \otimes \Hom_{\cC}(A,B)$. We often write $\bQ \otimes A$ for the class of $A \in \cC$ in $\bQ \otimes \cC$ to avoid confusions.
\subsection{Crystalline complexes}\label{sectA-1}
We construct $\bE_{\logAb/\Wb}, \bE_{\logAb/\logCb} \in D^+(Y_\et,\bZpd)$ and $\bE_{\logBb/\Wb} \in D^+(\ol Y_\et,\gkbZpd)$,
where $\Wn$ means $\Spec(\Wn)$ endowed with the trivial log structure, for each $n \ge 1$. See Definition \ref{def8-1} for $\logCn$.
To construct $\bE_{\logAb/\Wb}$, we fix an \'etale hypercovering $\logAs \to \logA$ and a closed immersion $\logAs \hra \logZs$ of simplicial fine log schemes over $W$ such that $\logZi$ is smooth over $W$ for each $i \in \bN$ (cf.\ \cite{HK} (2.18)). Put 
\[ Y^\star:=X^\star \otimes_{O_K} k, \] which is an \'etale hypercovering of $Y$. For $n\ge 1$, let $\logPns$ be the PD-envelope of $\logAns \hra \logZns$ with respect to the canonical PD-structure on $(p) \subset \Wn$ (\cite{K3} Definition 5.4), and we define a complex $\bE_{\logAns/\Wn}$ of sheaves on $Y^\star_\et$ as
\[ \cO_{\cD_n^\star} \os{d}\to \cO_{\cD_n^\star} \otimes_{\cO_{Z_n}^\star} \omega^1_{\logZns/\Wn} \os{d}\to \dotsb \os{d}\to \cO_{\cD_n^\star} \otimes_{\cO_{Z_n^\star}} \omega^q_{\logZns/\Wn} \os{d}\to \dotsb, \] where the first term is placed in degree $0$ and $\omega^q_{\logZns/\Wn}$ denotes the $q$-th differential module of $\logZns$ over $\Wn$ cf.\ \cite{K3} (1.7). See loc.\ cit.\ Theorem 6.2 for $d$. Regarding this complex as a complex of projective systems (with respect to $n \ge 1$) of sheaves on $Y^\star_\et$, we obtain a complex $\bE_{\logAbs/\Wb}$ of objects of $\Shv(Y^\star_\et,\bZpd)$. We then define
\[ \bE_{\logAb/\Wb} := R\theta_*(\bE_{\logAbs/\Wb}) \in D^+(Y_\et,\bZpd), \]
where $\theta : \Shv(Y^\star_\et,\bZpd) \to \Shv(Y_\et,\bZpd)$ denotes the natural morphism of topoi. The resulting object $\bE_{\logAb/\Wb}$ is independent of the choice of the pair $(\logAs,\logZs)$ by a standard argument (cf.\ \cite{KaV} p.\ 212).
\par
To construct $\bE_{\logAb/\logCb}$, we put $V:=\Spec(W[t])$, and define $M_V$ as the log structure on $V$ associated with the divisor $\{t=0\}$. We regard $\logA$ as a log scheme over $\logF$ by the composite map \[ \logA \to \logE \to \logF, \] where the last map is given by $T \mapsto \pi$, the prime element of $O_K$ we fixed in Definition \ref{def8-1} to define $\logCn$. We then apply the same construction as for $\bE_{\logAb/\Wb}$ to the morphism $\logA \to \logF$, i.e, fix an \'etale hypercovering $\logAs \to \logA$ and a closed immersion $i^\star : \logAs \hra \logZs$ of simplicial log schemes over $\logF$ such that $\logZi$ is smooth over $W$ for each $i \in \bN$. We obtain $\bE_{\logAb/\logCb}$ by replacing $\omega^q_{\logZns/\Wn}$ with $\omega^q_{\logZns/\logFn}$.
\par
We construct $\bE_{\logBb/\Wb}$ as follows. Fix an \'etale hypercovering $\logAs \to \logA$ and for each finite extension $L/K$ contained in $\ol K$, fix a closed immersion $\logALs \hra \logZLs$ ($\logALs:=\logAs \times_{\logE} \logEL$) of simplicial fine log schemes over $W$ such that $\logZLi$ is smooth over $W$ for each $i \in \bN$ and $L/K$, and such that for finite extensions $L'/L$ there are morphisms $\tau_{L'/L} : \logZLps \to \logZLs$ satisfying transitivity. For a finite extension $L/K$, let $k_L$ be the residue field of $L$, and put
\[ Y_L^\star:=X^\star \otimes_{O_K} k_L \quad \hbox{ and } \quad  \ol Y{}^\star:=X^\star \otimes_{O_K} \ol k, \] which are \'etale hypercoverings of $Y_L:=Y \otimes_k k_L$ and $\ol Y$, respectively. We define a complex $\bE_{\logALbs/\Wb}$ on $Y_L^\star$ applying the same construction as for $\bE_{\logAbs/\Wb}$ to the embedding $\logALs \hra \logZLs$, whose inverse image onto $\ol Y{}^\star$ yields an inductive system of complexes of objects of $\Shv(\ol Y{}^\star_\et,\gkbZpd)$ with respect to finite extensions $L/K$. We then define
\begin{align*}
 \bE_{\logBbs/\Wb} &:= \varinjlim_{K \subset L \subset \ol K} \ \bE_{\logALbs/\Wb}\big|{}_{\ol Y{}^\star}, \\
 \bE_{\logBb/\Wb} &:= R\ol \theta_*(\bE_{\logBbs/\Wb}) \in D^+(\ol Y_\et,\gkbZpd),
\end{align*}
where $\ol \theta : \Shv(\ol Y{}^\star_\et,\gkbZpd) \to \Shv(\ol Y_\et,\gkbZpd)$ denotes the natural morphism of topoi.

\par
We define the following objects of $\Mod(\bZpd)$:
\begin{align}
\label{eqA-001} H^i_\crys(\logAb/\Wb) &:= R^i\vG(\bE_{\logAb/\Wb}),\\
\label{eqA-002} H^i_\crys(\logAb/\logCb) &:= R^i\vG(\bE_{\logAb/\logCb}),\\
\label{eqA-003} H^i_\crys(\logBb/\Wb) &:= R^i\ol\vG\big(\bE_{\logBb/\Wb}\big),
\end{align}
where $\vG$ and $\ol \vG$ denote the following left exact functors, respectively:
\begin{align}
\label{eqA-021} \vG & := \vG(Y,-) : \Shv(Y_\et,\bZpd) \to \Mod(\bZpd), \\
\label{eqA-022} \ol \vG & := \vG(\ol Y,-) : \Shv(\ol Y_\et,\gkbZpd) \to \Mod(\gkbZpd).
\end{align}
We define `naive' crystalline cohomology groups \[ H^i_\crys(\logA/W), \;\; H^i_\crys(\logA/\logC), \;\;  H^i_\crys(\logB/W)\]  as the projective limit of  \eqref{eqA-001}--\eqref{eqA-003}, respectively.
\addtocounter{athm}{5}
\begin{arem}
Let $\pi_n : \Shv(Y_\et,\bZpd) \to \Shv(Y_\et,\bZ/p^n)$ be the natural functor sending $\{\cF_m\}_{m \ge 1} \mapsto \cF_n$. Since $\pi_n$ is exact, it extends to a triangulated functor \[ \pi_n : D^+(Y_\et,\bZpd) \lra D^+(Y_\et,\bZ/p^n), \] which sends $\bE_{\logAb/\Wb} \mapsto \bE_{\logAn/\Wn}$,
 an object computing the crystalline cohomology of $\logAn/\Wn$, because $\pi_n$ is compatible with the gluing functor $R\theta_*$. Moreover by \cite{J} Proposition {\rm1.1\,(b)}, we have
\[ H^i_\crys(\logAb/\Wb) \simeq \{ H^i_\crys(\logAn/\Wn) \}_{n \ge 1}. \]
We have similar facts for \eqref{eqA-002} and \eqref{eqA-003} as well.
\end{arem}

\subsection{Continuous crystalline cohomology}\label{sectA-2}
We define continuous crystalline cohomology groups as follows:
\begin{align}
\label{eqA-011} H^i_\cocrys(\logA/W) &:=R^i\big(\varprojlim \circ \vG\big) \big(\bE_{\logAb/\Wb}\big),\\
\label{eqA-012} H^i_\cocrys(\logA/\logC) &:=R^i\big(\varprojlim \circ \vG\big) \big(\bE_{\logAb/\logCb}\big),
\end{align}
where $\vG$ is as in \eqref{eqA-021}. Because $\vG$ has an exact left adjoint, it preserves injectives and there exists a spectral sequence
\[ E_2^{a,b}=\Ri a \ H^b_\crys(\logAb/\Wb) \Lra H^{a+b}_\cocrys(\logA/W). \]
Because $\Ri a = 0$ for $a \ge 2$, this spectral sequence breaks up into short exact sequences
\begin{align}
\label{eqA-0} 0 \lra \Ri 1 \ H^{i-1}_\crys(\logAb/\Wb) & \lra H^i_\cocrys(\logA/W) \\
\notag & \lra H^i_\crys(\logA/W) \lra 0.
\end{align}
We have similar exact sequences for \eqref{eqA-012} by the same arguments.
\addtocounter{athm}{3}
\begin{aprop}\label{propA-20}
There exists a distinguished triangle
\[ \bE_{\logAb/\Wb} \os{\can}\lra \bE_{\logAb/\logCb} \os{\nu}\lra \bE_{\logAb/\logCb} \os{\frac{dt}{t} \wedge }\lra \bE_{\logAb/\Wn} [1] \]
in $D^+(Y_\et,\bZpd)$. Consequently there is a long exact sequence
{\begin{align*} \dotsb & \to H^i_\cocrys(\logA/W) \to H^i_\cocrys(\logA/\logC) \os{\nu}\to H^i_\cocrys(\logA/\logC) \\
& \to H^{i+1}_\cocrys(\logA/W) \to \dotsb. \end{align*}}Moreover the map $\nu$ fits into a commutative diagram {\rm(}see Definition {\rm\ref{def8-1}} for $N${\rm)}
\[\xymatrix{
 H^i_\cocrys(\logA/\logC) \ar[r]^-{\nu} \ar[d] & H^i_\cocrys(\logA/\logC) \ar[d] \\
 H^i_\crys(\logA/\logC) \ar[r]^-{N} & H^i_\crys(\logA/\logC)\,.}\]
{\rm(}We will see that the vertical arrows are bijective later in Corollary {\rm\ref{corA-1}\,(1)} below.{\rm)}
\end{aprop}
\begin{pf}
The assertion follows from the same arguments as in the proof of \cite{K4} Lemma 4.2 with ${\mathcal F}=\cO^{\crys}$.
\end{pf}
\subsection{Comparison of projective systems of crystalline cohomology}\label{sectA-3}
For a finite extension $L/K$, let $\logAL$ be as in the beginning of this appendix, and put
\[ H^i_\crys(\logBn/\logCn):=\varinjlim_{K \subset L \subset \ol K} \ H^i_\crys(\logALn/\logCn), \]
where $L$ runs through all finite field extensions of $K$ contained in $\ol K$. We recall here the following comparison facts on projective systems of crystalline cohomology groups, which will be useful later. Note that for $A_\bullet \in \Mod(\bZpd)$, both $\varprojlim \ A_\bullet$ and $R^1\varprojlim \ A_\bullet$ have finite exponents if $\bQ \otimes (A_\bullet) \simeq 0$ in $\qt \Mod(\bZpd)$.
\begin{aprop}\label{propA-1}
Let $i$ be a non-negative integer, and put $D^i_n:=\bH^i(Y,\mwitt Y n \bullet)$. Let $P_n$ be the ring defined in \cite{T1} {\rm\S1.6}, and put $R_{\cE_n}:=\vG(\cE_n,\cO_{\cE_n})$.
\begin{enumerate}
\item[(1)] There is an isomorphism in $\bQ \otimes \Mod(\bZpd)$
\[ \bQ\otimes \big(H^i_\crys(\logAb/\logCb)\big) \simeq \bQ\otimes \{ R_{\cE_n} \otimes_{\Wn} D^i_n \}_{n \ge 1} \]
and $\Ri 1 \ H^i_\crys(\logAb/\logCb)$ has a finite exponent.
\item[(2)] There is an isomorphism in $\bQ \otimes \Mod(\gkbZpd)$
\[ \bQ\otimes \{ H^i_\crys(\logBn/\logCn)\}_{n \ge 1} \simeq \bQ\otimes \{ P_n \otimes_{\Wn} D^i_n \}_{n \ge 1} \]
and $\Ri 1 \ \{ H^i_\crys(\logBn/\logCn)\}_{n \ge 1}$ has a finite exponent.
\item[(3)] There is an isomorphism in $\bQ \otimes \Mod(\gkbZpd)$
\[ \bQ\otimes \big((H^i_\crys(\logBb/\Wb)\big) \simeq \bQ\otimes \big\{(P_n \otimes_{\Wn} D^i_n)^{N=0} \big\}{}_{n \ge 1} \]
and $\Ri 1 \ H^i_\crys(\logBb/\Wb)$ has a finite exponent. Here $N$ acts on $P_n \otimes_{\Wn} D^i_n$ by $N_{P_n} \otimes 1 + 1 \otimes N$ and $N_{P_n}$ denotes the monodromy operator of $P_n$ \cite{T1} p.\ {\rm253}.
\end{enumerate}
\end{aprop}
\begin{pf}
(1) The isomorphism in the assertion follows from \cite{HK} Lemma 5.2. As for the second assertion, since $D^i_n$ is finitely generated over $\Wn$, it is enough to check that the projection $R_{\cE_{n+1}} \to R_{\cE_n}$ is surjective, which implies that the projective system $\{ R_{\cE_n} \otimes_{\Wn} D^i_n \}_{n \ge 1}$ satisfies the Mittag-Leffler condition. This surjectivity follows from the presentation
\[ R_{\cE_n} = W[t,\,t^{e\nu}/\nu! \;(\nu \ge 1)]\otimes_W \Wn \quad (e:=[K : K_0]) \]
obtained from the definition of $\logCn$ (cf.\ the proof of \cite{T1} Proposition 4.4.6).
\par
(2) The isomorphism in the assertion follows from \cite{T1} Proposition 4.5.4. The second assertion follows from the facts that the natural projection $P_{n+1} \to P_n$ is surjective (loc.\ cit.\ Lemma 1.6.7) and that $D^i_n$ is finitely generated over $\Wn$.
\par
(3) See \cite{K4} (4.5) for the isomorphism in the assertion. We show the second assertion. Note that $P_n$ is flat over $\Wn$, because $R_{\cE_n}$is flat over $\Wn$ and $P_n$ is flat over $R_{\cE_n}$ by the above presentation of $R_{\cE_n}$ and \cite{T1} Proposition 4.1.5. Let $\{ D_n\}_{n \ge 1}$ be a projective system of $W$-modules such that $D_n$ is a finitely generated $\Wn$-module for each $n$, and let $N : \{ D_n\}_{n \ge 1} \to \{ D_n\}_{n \ge 1}$ be a nilpotent $W$-endomorphism. Our task is to show that \[ \Ri 1 \ \{ \big(P_n \otimes_{\Wn} D_n \big){}^{N=0} \}_{n \ge 1} = 0. \] Consider a short exact sequence of projective systems
\[ 0 \lra \{(D_n)^{N=0}\}_{n\ge 1} \lra \{D_n\}_{n\ge 1} \lra \{N(D_n)\}_{n\ge 1} \lra 0. \]
Note that $(D_n)^{N=0}$ and $N(D_n)$ are finitely generated over $\Wn$ for each $n$. Let $b$ be the minimal integer for which $N^b=0$ on $\{D_n\}_{n\ge 1}$. By \cite{K4} Lemma 4.3 and the flatness of $P_n$ over $\Wn$, we have a short exact sequence for each $n \ge 1$
\[ 0 \to (P_n \otimes_{\Wn} (D_n)^{N=0})^{N=0} \to (P_n \otimes_{\Wn} D_n)^{N=0} \to (P_n \otimes_{\Wn} N(D_n))^{N=0} \to 0, \]
which yields a short exact sequence of projective systems with respect to $n \ge 1$. Since $N^{b-1}=0$ on $\{N(D_n)\}_{n\ge 1}$, we may assume $b=1$ by induction on $b \ge 1$. Now let $B_n$ be as in \cite{T1} \S1.6 and let $A_\crys$ be as in loc.\ cit.\ \S1.1. Then we have isomorphisms
\[ (P_n \otimes_{\Wn} D_n)^{N=0} \os{(1)}\simeq (P_n)^{N=0} \otimes_{\Wn} D_n \os{(2)}\simeq B_n \otimes_{\Wn} D_n \os{(3)}= (A_\crys/p^n) \otimes_{\Wn} D_n, \]
where (1) follows from the assumption $b=1$ and the flatness of $P_n$ over $\Wn$. The isomorphism (2) (resp.\ (3)) follows from \cite{T1} Corollary 1.6.6 (resp.\ the definition of $B_n$ in loc.\ cit.\ \S1.6). Thus $\Ri 1 \ \{ (P_n \otimes_{\Wn} D_n)^{N=0} \}_{n \ge 1}$ is zero, and we obtain the assertion.
\end{pf}
\begin{acor}\label{corA-1}
\begin{enumerate} \item[(1)] For $i \ge 0$, we have
\[ H^i_\crys(\logA/\logC)_\qp \simeq H^i_\cocrys(\logA/\logC)_\qp\,. \] 
\item[(2)] The torsion subgroups of
\[ H^i_\crys(\logB/\logC) \quad \hbox{ and } \quad H^i_\crys(\logB/W) \]
have finite exponents for any $i \ge 0$.
\end{enumerate}
\end{acor}
\begin{pf}
(1) follows from Proposition \ref{propA-1}\,(1) and the remark after \eqref{eqA-0}. Since $\varprojlim{}_{n \ge 1} \ P_n$ is $p$-torsion free, the assertion (2) follows from the isomorphisms in Proposition \ref{propA-1}\,(2) and (3). 
\end{pf}
\begin{arem}\label{remA-1}
If $K$ is absolutely unramified, then we have
\[ H^i_\crys(\logAn/\Wn) \simeq \bH^i(Y,\mtwitt Y n \bullet), \]
which is finitely generated over $\Wn$ by Theorem {\rm\ref{thm6-1}\,(2)} and \eqref{eq6-6} {\rm(}see also \cite{H2} {\rm(1.4.3), (2.4.2))}. Consequently, the projective system $H^i_\crys(\logAb/\Wb)$ satisfies the Mittag-Leffler condition and we have a long exact sequence
\begin{align*} \dotsb & \to H^i_\crys(\logA/W)_\qp \to H^i_\crys(\logA/\logC)_\qp \os{N}\to H^i_\crys(\logA/\logC)_\qp \\
& \to H^{i+1}_\crys(\logA/W)_\qp \to \dotsb, \end{align*}
which removes the assumption $(*)$ in \cite{L} p.\ {\rm191}. On the other hand, the author does not know if $H^i_\crys(\logAb/\Wb)$ satisfies the Mittag-Leffler condition even up to torsion, when $K$ is not absolutely unramified. 
\end{arem}
The following corollary has been used in the proof of Lemma \ref{lem8-2}\,(1):
\begin{acor}\label{propA-3}
In the following commutative diagram of canonical maps, the arrows {\rm(3)} and {\rm(4)} are injective{\rm:}
\[\xymatrix{
 H^i_\cocrys(\logA/W)_\qp \ar[r]^-{{\sf(1)}} \ar[d]_{{\sf(2)}} & H^i_\crys(\logA/\logC)_\qp \; \ar[d]^{{\sf(3)}} \\
 H^i_\crys(\logB/W)_\qp \ar[r]^-{{\sf(4)}} & H^i_\crys(\logB/\logC)_\qp\,.
}\]
In particular, the kernel of {\rm(1)} agrees with that of {\rm(2)}.
\end{acor}
\begin{pf}
The injectivity of (3) follows from Proposition \ref{propA-1}\,(1) and (2) and the injectivity of the natural maps $R_{\cE_n} \to P_n$ for $n \ge 1$ (\cite{T1} Proposition 4.1.5). The injectivity of (4) follows from Proposition \ref{propA-1}\,(2) and (3). 
\end{pf}

\subsection{continuous-Galois crystalline cohomology}\label{sectA-4}
We define the continuous-Galois crystalline cohomology as follows:
\[ H^i_\dccrys(\logA/W):=R^i\big(\varprojlim \vG_\Gal \ol \vG\big)\big(\bE_{\logBb/\Wb}\big), \]
where $\ol \vG$ is as in \eqref{eqA-022}, and $\vG_\Gal$ denotes the functor taking $G_K$-invariant subgroups:
\[ \vG_\Gal:=\vG(G_K,-) : \Mod(\gkbZpd) \lra \Mod(\bZpd). \]
There is a natural map
\begin{equation}\label{eqA-1} H^i_\cocrys(\logA/W) \lra H^i_\dccrys(\logA/W) \end{equation}
by definition.
\stepcounter{athm}
\begin{athm}\label{thmA-1}
There exists a Hochschild-Serre spectral sequence
\[ E_2^{a,b}=H^a(K,H^b_\crys(\logB/W)_\qp) \Lra H^{a+b}_\dccrys(\logA/W)_\qp. \]
\end{athm}
\begin{pf}
Because $\ol \vG$ has an exact left adjoint functor, it preserves injectives and there exists a spectral sequence
\[ E_2^{a,b}=H^a(K,H^b_\crys(\logBb/\Wb)) \Lra H^{a+b}_\dccrys(\logA/W), \]
where $H^*(K,\{F_n\}_{n \ge 1})$ for a projective system $\{F_n\}_{n \ge 1}$ of discrete $G_K$-modules denotes the continuous Galois cohomology of $G_K$ in the sense of Jannsen \cite{J} \S2. By Corollary \ref{corA-1}\,(2) and loc.\ cit.\ Theorem 5.15\,(c), we have
\[ H^a(K,H^b_\crys(\logB/W))_\qp \simeq H^a(K,H^b_\crys(\logB/W)_\qp) \]
for $a,b \ge 0$. It remains to show that the canonical map
\[ H^a(K,H^b_\crys(\logB/W)_\qp) \lra H^a(K,H^b_\crys(\logBb/\Wb))_\qp \]
(loc.\ cit.\ Proof of  Theorem 2.2) is bijective for $a,b \ge 0$. Put
\begin{align*}
 L_\bullet^b &:= \{P_n \otimes_{\Wn} D_n^b\}_{n\ge 1}\,, \qquad L^b:=\qp \otimes_\zp \varprojlim \ L_\bullet^b = B_\st^+ \otimes_{K_0} D^b\,, \\
 T_\bullet^b &:= H^b_\crys(\logBb/\Wb) \quad \hbox{ and } \quad T^b:=\qp \otimes_\zp H^b_\crys(\logB/W)\,.
\end{align*}
Note that we have $H^a(K,L^b) \simeq H^a(K,\varprojlim \ L_\bullet^b)_\qp$ by Corollary \ref{corA-1}\,(2) and \cite{J} Theorem 5.15\,(c).
We have a short exact sequence in $\bQ \otimes \Mod(\gkbZpd)$
\[ 0 \lra \bQ \otimes (T_\bullet^b) \lra \bQ \otimes (L_\bullet^b) \os{N}\lra \bQ \otimes (L_\bullet^b) \lra 0 \]
by \cite{K4} Lemma 4.3, and a short exact sequence of topological $G_K$-modules
\stepcounter{equation}
\begin{equation}\label{eqA-2} 0 \lra T^b \lra L^b \os{N}\lra L^b \lra 0 \end{equation}
by Proposition \ref{propA-1}\,(3),
which yield a commutative diagram with exact rows
\[ \begin{CD}
\dotsb \lra \, @. H^{a-1}(K,L^b)  @. \, \lra \, @. H^a(K,T^b)  @. \, \lra \, @. H^a(K,L^b)  @. \, \os{N}\lra \, @. H^a(K,L^b) @. \,\lra \dotsb\phantom{,}  \\
@. @V{f_{a-1}}VV @. @VVV @. @V{f_a}VV @. @V{f_a}VV \\
\dotsb \lra \, @. H^{a-1}(K,L_\bullet^b)_\qp @. \, \lra \, @. H^a(K,T_\bullet^b)_\qp @. \, \lra \, @. H^a(K,L_\bullet^b)_\qp @. \, \os{N}\lra \, @. H^a(K,L_\bullet^b)_\qp @. \,\lra \dotsb,
\end{CD}\]
where $f_a$ are bijective for all $a \ge 0$ by \cite{J} Theorem 2.2. Hence the assertion follows from the five lemma.
\end{pf}
\stepcounter{athm}
\begin{acor}\label{corA-2}
Put
\[H^{i+1}_\cocrys(\logA/W)_\qp^0:=\ker(H^{i+1}_\cocrys(\logA/W)_\qp \to H^{i+1}_\crys(\logB/W)_\qp). \]
Then the canonical homomorphism
\[ H^{i+1}_\cocrys(\logA/W)_\qp^0 \lra H^1(K,H^i_\crys(\logB/W)_\qp), \]
induced by \eqref{eqA-1}, fits into a commutative diagram
\[\xymatrix{ H^i_\crys(\logA/\logC)_\qp \ar@{->>}[r]^-\partial \ar[d] & H^{i+1}_\cocrys(\logA/W)_\qp^0 \ar[d] \\
(B_\st^+ \otimes_{K_0} D^i)^{G_K} \ar[r]^-{(-1)^{i+1}\delta} & H^1(K,H^i_\crys(\logB/W)_\qp), }\]
where $\partial$ is the connecting map in Proposition {\rm\ref{propA-20}} {\rm(}see also Corollaries {\rm\ref{corA-1}\,(2)} and {\rm \ref{propA-3})}, and $\delta$ denotes the connecting map of continuous Galois cohomology associated with \eqref{eqA-2} with $b=i$. The left vertical arrow is obtained from Proposition {\rm\ref{propA-1}\,(2)}.
\end{acor}
\begin{pf}
The second assertion follows from simple computations on the boundary maps of cohomology groups arising from the distinguished triangle
\begin{align*}
R\ol\vG\big(\bE_{\logBb/\Wb}\big) \to R\ol\vG\big(\bE_{\logBb/\logCb}\big) \os{\nu}\to R\ol\vG\big(\bE_{\logBb/\logCb}\big) \os{+1}\to
\end{align*}
in $D^+(\Mod(\gkbZpd))$ (cf.\ Proposition \ref{propA-20}).
The sign $(-1)^{i+1}$ in the diagram arises from the orientation of the distinguished triangle in Proposition \ref{propA-20} and the fact that the construction of connecting morphisms (in the derived category) associated with short exact sequences of complexes commutes with the shift functor $[i]$ up to the sign $(-1)^i$. The details are straight-forward and left to the reader.
\end{pf}

\subsection{Syntomic complexes}\label{sectA-5pre}
We construct the following objects for $r \ge 0$:
\[ \cSbt(r)_{\logA} \in D^+(Y_\et,\bZpd) \quad \hbox{and} \quad \cSbt(r)_{\logB} \in D^+(\ol Y_\et,\gkbZpd), \]
and the following objects for $r$ with $0 \le r \le p-1$
\[ \cSb(r)_{\logA} \in D^+(Y_\et,\bZpd) \quad \hbox{and} \quad \cSb(r)_{\logB} \in D^+(\ol Y_\et,\gkbZpd). \]
\begin{adefn}
Let $\logT$ be a log scheme over $\zp$. A Frobenius endomorphism $\varphi : \logT \to \logT$ is a morphism over $\bZ_p$ such that $\varphi \otimes \bZ/p : \logTo \to \logTo$ is the absolute Frobenius endomorphism in the sense of \cite{K3} Definition {\rm 4.7}.
\end{adefn}
To construct $\cSbt(r)_{\logA}$ and $\cSb(r)_{\logA}$, we fix an \'etale hypercovering $\logAs \to \logA$ and a closed immersion $\logAs \hra \logZs$ of simplicial fine log schemes over $W$ such that $\logZi$ is smooth over $W$ and has a Frobenius endomorphism for each $i \in \bN$. Let $n \ge 1$ be an integer, and let $\logPns$ be the PD-envelope of $\logAns$ in $\logZns$ with respect to the canonical PD-structure on $(p) \subset \Wn$. For $i \ge 1$, let $\cJ^{[i]} \subset \cO_{\cD_n^\star}$ be the $i$-th divided power of the ideal $\cJ=\ker(\cO_{\cD_n^\star} \to \cO_{X_n^\star})$. For $i \le 0$, put $\cJ^{[i]}:=\cO_{\cD_n^\star}$.
Let $\bJ^{[r]}_{\logAns/\Wn}$ be the complex of sheaves on $Y^\star_\et$ \[ \cJ^{[r]} \os{d}\to \cJ^{[r-1]}\otimes_{\cO_{Z_n^\star}} \omega^1_{\logZns/\Wn} \os{d}\to \dotsb \os{d}\to \cJ^{[r-q]}\otimes_{\cO_{Z_n^\star}} \omega^q_{\logZns/\Wn} \os{d}\to \dotsb, \] where $\cJ^{[r]}$ is placed in degree $0$. See \cite{T2} Corollary 1.10 for $d$. The complex $\bE_{\logAns/\Wn}$ we considered in \S\ref{sectA-1} agrees with $\bJ^{[0]}_{\logAns/\Wn}$. We define a complex $\cSnt(r)_{\logAs}$ on $Y^\star_\et$ as the mapping fiber of the homomorphism
\[ p^r-\varphi_n^* : \bJ^{[r]}_{\logAns} \to \bE_{\logAns} \] For $0 \le r \le p-1$, the Frobenius endomorphism on $\logZnrs$ induces a homomorphism of complexes
\[ f_r:=\ol{p^{-r} \cdot \varphi_{n+r}^*} :  \bJ^{[r]}_{\logAns/\Wn} \lra \bE_{\logAns/\Wn} \]
(cf.\ \cite{T2} p.\ 540). We define a complex $\cSn(r)_{\logAs}$ on $Y^\star_\et$ as the mapping fiber of the homomorphism
\[ 1-f_r:  \bJ^{[r]}_{\logAns} \to \bE_{\logAns} \]
Regarding $\cSnt(r)_{\logAs}$ and $\cSn(r)_{\logAs}$ as complexes of projective systems (on $n \ge 1$) of sheaves on $Y^\star_\et$, we define
\begin{align*}
 \cSbt(r)_{\logA}& := R\theta_*(\cSbt(r)_{\logAs}) \qquad \hbox{($r \ge 0$)}, \\
 \cSb(r)_{\logA}& := R\theta_*(\cSb(r)_{\logAs}) \qquad \hbox{($0 \le r \le p-1$),}
\end{align*}
where $\theta : \Shv(Y^\star_\et,\bZpd) \to \Shv(Y_\et,\bZpd)$ denotes the natural morphism of topoi. The resulting objects are independent of the pair $(\logAs,\logZs)$ (cf.\ \cite{KaV} p.\ 212).
\par
We construct $\cSbt(r)_{\logB}$ for $r \ge 0$ and $\cSb(r)_{\logB}$ for $0 \le r \le p-1$ as follows. Fix an \'etale hypercovering $\logAs \to \logA$ and for each finite extension $L/K$ contained in $\ol K$, fix a closed immersion $\logALs \hra \logZLs$ ($\logALs:=\logAs \times_{\logE} \logEL$) of simplicial fine log schemes over $W$ such that $\logZLi$ is smooth over $W$ for each $i \in \bN$ and $L/K$, such that $Z_L^i$ has a Frobenius endomorphism for each $i \in \bN$ and $L/K$, and such that for finite extensions $L'/L$ there are morphisms $\tau_{L'/L} : \logZLps \to \logZLs$ which satisfy transitivity and compatibility with Frobenius morphisms. For a finite extension $L/K$, we define $Y_L^\star$ in a similar way as in \S\ref{sectA-1}. We define complexes $\cSbt(r)_{\logAL}$ and $\cSb(r)_{\logAL}$ on $Y_L^\star$ applying the same constructions as for the complexes $\cSbt(r)_{\logA}$ and $\cSb(r)_{\logA}$, respectively, to the embedding $\logALs \hra \logZLs$, whose inverse images onto $\ol Y{}^\star:=X^\star \otimes_{O_K} \ol k$ yield inductive systems of complexes of objects of $\Shv(\ol Y{}^\star_\et,\gkbZpd)$ with respect to finite extensions $L/K$. We define
\begin{align*}
 \cSbt(r)_{\logBbs} &:= \varinjlim_{K \subset L \subset \ol K} \ \cSbt(r)_{\logALbs}\big|{}_{\ol Y{}^\star}, \\
 \cSbt(r)_{\logBb} &:= R\ol \theta_*(\cSbt(r)_{\logBbs}) \in D^+(\ol Y_\et,\gkbZpd),
\end{align*}
where $\ol \theta : \Shv(\ol Y{}^\star_\et,\gkbZpd) \to \Shv(\ol Y_\et,\gkbZpd)$ denotes the natural morphism of topoi. We construct $\cSb(r)_{\logBb}$ from $\cSb(r)_{\logALbs}$'s in a similar way. When $0 \le r \le p-1$, there is a canonical projection
\[ \psi_\bullet^r : \cSbt(r)_{\logB} \lra \cSb(r)_{\logB} \]
induced by the identity map of $\bE_{\logBb/\Wb}$. There are canonical morphisms
\addtocounter{equation}{1}
\begin{align}
\label{eqA-cant}
\wt{c}_\bullet^r &: \cSbt(r)_{\logB} \lra \bJ^{[r]}_{\logBb/\Wb} \os{1}\lra \bE_{\logBb/\Wb} \quad \hbox{($r \ge 0$)} \\
\label{eqA-can}
c_\bullet^r &: \cSb(r)_{\logB} \lra \bJ^{[r]}_{\logBb/\Wb} \os{1}\lra \bE_{\logBb/\Wb} \quad \hbox{($0 \le r \le p-1$),}
\end{align}
which satisfy
\begin{equation}\label{eqA-trans}
\qquad\qquad\qquad p^r \cdot \wt{c}_\bullet^r =  c_\bullet^r \circ \psi_\bullet^r \qquad \hbox{(when $0 \le r \le p-1$).}
\end{equation}

For $r \ge 0$, we define \[ \bZpd(r)' := (p^a a!)^{-1} \zp(r) \otimes \bZpd \in \Shv(\Xggf,\gkbZpd), \] where $a$ denotes the maximal integer $\le r/(p-1)$. We have $\bZpd(r)' = \mu_{p^\bullet}^{\otimes r}$ canonically when $r \le p-2$. Similarly for $r < 0$, we define
\begin{equation}\label{eq8-t}
 \bZpd(r)' := p^a a!\,\zp(r) \otimes \bZpd \in \Shv(\Xggf,\gkbZpd)
\end{equation}
with $a$ the maximal integer $\le -r/(p-1)$, which will be used in \eqref{eqA-3} below. Let $\ol j : \Xggf \hra \ol X$ and $\ol \iota: \ol Y \hra \ol X$ be the natural immersions. Let $R\ol j_*$ and $\ol \iota^*$ be the following functors, respectively:
\begin{align*}
R\ol j_* & : D^+((\Xggf)_\et,\gkbZpd) \lra D^+(\ol X_\et,\gkbZpd), \\
\ol \iota^* & : D^+(\ol X_\et,\gkbZpd) \lra D^+(\ol Y_\et,\gkbZpd).
\end{align*}
\addtocounter{athm}{4}
\stepcounter{equation}
\begin{athm}[{\bf\cite{T1} \S3.1, \cite{K4} Theorem 5.4, cf.\ \cite{T2} Theorem 5.1}]\label{thmA-3}
For $r \ge 0$, there exists a canonical morphism in $D^+(\ol Y_\et,\gkbZpd)$ compatible with product structures
\[ \eta_\bullet^r : \cSbt(r)_{\logB} \lra \ol \iota^* R\ol j_* \bZpd(r)'. \]
If $r \le p-2$, then $\eta_\bullet^r$ factors through an isomorphism
\begin{equation}\label{eq8-00}
 \cSb(r)_{\logB} \isom \tau_{\le r} \, \ol \iota^* R\ol j_* \mu_{p^\bullet}^{\otimes r}.
\end{equation}
\end{athm}
\begin{pf}
We define $\eta_\bullet^r$ applying the arguments in \cite{T1} \S3.1 in the category of $\bZpd$-sheaves. If $r \le p-2$, then we have $\bZpd(r)'=\mu_{p^\bullet}^{\otimes r}$, and $\eta_\bullet^r$ factors through a morphism $\ol \eta_\bullet^r : \cSb(r)_{\logB} \to \ol \iota^* R\ol j_* \mu_{p^\bullet}^{\otimes r}$ by the construction of $\eta_\bullet^r$ (cf.\ \cite{T1} (3.1.11)). The morphism $\ol \eta_\bullet^r$ induces an isomorphism as claimed, because $\cSb(r)_{\logB}$ is concentrated in $[0,r]$ and $\ol \eta_\bullet^r$ induces isomorphisms on the $q$-th cohomology objects with $0 \le q \le r$ by \cite{K4} 5.4.
\end{pf}
We define
\begin{align*}
H^i_\syn(\logA,\cSb(r)) &:= R^i\vG(\cSb(r)_{\logA}),\\
H^i_\syn(\logA,\cS_\zp(r)) &:=\varprojlim \ H^i_\syn(\logA,\cS_\bullet(r)),\\
H^i_\syn(\logA,\cS_\qp(r)) &:=\qp \otimes_\zp H^i_\syn(\logA,\cS_\zp(r)),\\
H^i_\syn(\logB,\cS_\qp(r)) &:=\qp \otimes_\zp \varprojlim \ R^i\ol\vG\big(\cSb(r)_{\logB}\big).
\end{align*}
where $\vG$ and $\ol \vG$ are as in \eqref{eqA-021} and \eqref{eqA-022}, respectively.
\subsection{continuous(-Galois) syntomic cohomology}\label{sectA-5}
We assume $r \le p-2$ in what follows. For $i \ge 0$, we define the continuous syntomic cohomology as follows:
\[ H^i_\cosyn(\logA,\cS_\zp(r)):=R^i\big(\varprojlim \vG\big)\big(\cS_\bullet(r)_{\logA}\big). \]
Similarly, we define the continuous-Galois syntomic cohomology as follows:
\[ H^i_\dcsyn(\logA,\cS_\zp(r)):=R^i\big(\varprojlim \vG_\Gal \ol\vG\big)\big(\cS_\bullet(r)_{\logB}\big). \]
We put
\begin{align*}
H^i_\cosyn(\logA,\cS_\qp(r)) &:=\qp \otimes_\zp H^i_\cosyn(\logA,\cS_\zp(r)).
\end{align*}
\begin{aprop}\label{propA-5}
Let $i \ge 0$ be an integer.
\begin{enumerate}
\item[(1)]
Let $\eta: H^i_\syn(\logA,\cS_\qp(r)) \to H^i(\Xggf,\qp(r))$ be as in Lemma {\rm\ref{lem8-1}}. Then the kernel of the composite map
\[ H^i_\cosyn(\logA,\cS_\qp(r)) \lra H^i_\syn(\logA,\cS_\qp(r)) \os{\eta}\lra H^i(\Xggf,\qp(r)) \]
agrees with that of the composite map
\[ H^i_\cosyn(\logA,\cS_\qp(r)) \lra H^i_\syn(\logA,\cS_\qp(r)) \lra H^i_\crys(\logB/W)_\qp, \]
which we denote by $H^i_\cosyn(\logA,\cS_\qp(r))^0$, in what follows.
\item[(2)]
If $K$ is a $p$-adic local field {\rm(}i.e., $k$ is finite{\rm)}, then we have
\[ H^i_\cosyn(\logA,\cS_\zp(r)) \isom H^i_\syn(\logA,\cS_\zp(r)). \]
In particular, we have the following canonical map in this case{\rm:}
\[ H^i_\syn(\logA,\cS_\zp(r)) \lra H^i_\cocrys(\logA/W). \]
\item[(3)]
If $r \ge d:=\dim(\Xgf)$, then we have
\[ H^i_\dcsyn(\logA,\cS_\zp(r)) \isom H^i_\cont(\Xgf,\zp(r)),\]
where $H^*_\cont(\Xgf,\zp(r))$ denotes the continuous \'etale cohomology of Jannsen \cite{J}.
\end{enumerate}
\end{aprop}
\begin{pf}
(1) The assertion immediately follows from Lemma \ref{lem8-1}.
\par
(2) There is a short exact sequence analogous to \eqref{eqA-0}
\begin{align*}
0 \lra \Ri 1 \ H^{i-1}_\syn(\logA,\cS_\bullet(r)) & \lra H^i_\cosyn(\logA,\cS_\zp(r)) \\
\notag & \lra H^i_\syn(\logA,\cS_\zp(r)) \lra 0.
\end{align*}
The group $H^{i-1}_\syn(\logA,\cS_n(r))$ is finite for any $i,n \ge 1$ by the properness of $X$ and the finiteness of $k$ (use \eqref{eq8-0} and the arguments in \cite{S2} (4.3.2), \S10.3 to reduced the problem to the case $n=1$, and then use Theorem \ref{thm2-1} of the main body). The assertion follows from these facts.
\par
(3) Since $r \ge d$ by assumption, \eqref{eq8-00} implies $\cSb(r)_{\logB} \simeq \ol \iota^*R\ol j_*\mu_{p^\bullet}^{\otimes r}$. Hence the assertion follows from the isomorphisms in $D^+(\Mod(\gkbZpd))$
\stepcounter{equation}
\begin{equation}\label{eqA-4} R\vG\big(\Xggf,  \mu_{p^\bullet}^{\otimes r}\big)=R\vG\big(\ol X, R\ol j_* \mu_{p^\bullet}^{\otimes r}\big) \isom R\vG\big(\ol Y, \ol \iota^* R\ol j_* \mu_{p^\bullet}^{\otimes r}\big), \end{equation}
where the last isomorphism is a consequence of the proper base-change theorem for the usual \'etale cohomology.
\end{pf}
For $i \ge 0$, put
\begin{align*}
H^i(\Xgf,\qp(r))^0&:=\ker(H^i(\Xgf,\qp(r)) \to H^i(\Xggf,\qp(r))), \\
H^i_\syn(\logA,\cS_\qp(r))^0&:=\ker(\eta : H^i_\syn(\logA,\cS_\qp(r)) \to H^i(\Xggf,\qp(r))).
\end{align*}
We have $H^i(\Xgf,\qp(r))=H^i_\cont(\Xgf,\qp(r))$ when $K$ is a $p$-adic local field. The following corollary is a consequence of Proposition \ref{propA-5}\,(1), (2) and the covariant functoriality of Hochschild-Serre spectral sequences in coefficients (see also the diagram in the proof of Theorem \ref{propA-4} below).
\stepcounter{athm}
\begin{acor}\label{thmA-2}
Assume that $K$ is a $p$-adic local field, and let $e$ be the composite map
\[ e : H^{i+1}_\syn(\logA,\cS_\qp(r))^0 \lra H^{i+1}(\Xgf,\qp(r))^0 \lra H^1(K,H^i(\Xggf,\qp(r))), \]
where the last map is an edge homomorphism of the Hochschild-Serre spectral sequence \eqref{eq8-hs}. If $r \ge d$, then $e$ fits into a commutative diagram
\[\xymatrix{ H^{i+1}_\syn(\logA,\cS_\qp(r))^0 \ar[r] \ar[d]_e & H^{i+1}_\cocrys(\logA/W)_\qp^0 \ar[d] \\
H^1(K,H^i(\Xggf,\qp(r))) \ar[r]^-{\beta^{i,r}} & H^1(K,H^i_\crys(\logB/W)_\qp), }\]
where the top arrow is obtained from Proposition {\rm\ref{propA-5}\,(1)} and {\rm(2)}, the right vertical arrow is the map in Corollary {\rm\ref{corA-2}} and $\beta^{i,r}$ is as in the proof of Theorem {\rm\ref{prop8-1}}.
\end{acor}
\par
We next construct a commutative diagram in $\qt D^+(\Mod(\gkbZpd))$ assuming $r < d$
\addtocounter{equation}{1}
\begin{equation}\label{eqA-3}
\xymatrix{\qt R\vG\big(\ol Y, \cSb(r)_{\logB}\big) \ar[d]_{f^r} \ar[rd]^{h^r} \\ 
\qt R\vG\big(\Xggf, \mu_{p^\bullet}^{\otimes r}\big) \ar[r]^-{g^r} & \qt \big(R\vG\big(\ol Y,\bE_{\logBb/\Wb}\big)\otimes \bZpd(r-d)'\big), }
\end{equation}
which is a key ingredient of the commutative diagram (A) of \S\ref{sect8} for the case $r<d$. See \eqref{eq8-t} for $\bZpd(r-d)'$. We define $h^r$ as the composite morphism
{\small
\begin{align*}
\qt R\vG\big(\ol Y,\cSb(r)_{\logB}\big) & \os{\eqref{eqA-can}}\lra \qt  R\vG\big(\ol Y,\bE_{\logBb/\Wb} \big) \\
& \,\isom \, \qt \big\{R\vG\big(\ol Y, \bE_{\logBb/\Wb}\big)\otimes \bZpd(d-r)' \otimes \bZpd(r-d)'\big\} \\
& \,\lra \, \qt \big\{ R\vG\big(\ol Y,\bE_{\logBb/\Wb}\big)\otimes \bZpd(r-d)'\big\}\,,
\end{align*}
}where the last arrow is induced by the $p^{d-r}$-times of the composite map
\begin{align}
\label{eqA-incl}
\qt \bZ/p^\bullet(d-r)'
& \os{\text{(1)}}\simeq  \qt H^0_\syn(\logBb,\cSbt(d-r)) \\
\notag & \os{\eqref{eqA-cant}}\hra \qt H^0_\crys(\logBb/\Wb)
\end{align}
(see \eqref{eqA-6} below for the isomorphism (1)) and the product of crystalline complexes. We define $f^r$ as the morphism induced by \eqref{eq8-00} and the isomorphisms in \eqref{eqA-4}.
To define $g^r$, we need an isomorphism
\begin{equation}\label{eqA-6} \wt{f}^d : \qt R\vG\big(\ol Y, \cSbt(d)_{\logB}\big) \simeq \qt R\vG(\Xggf, \bZpd(d)') \end{equation}
induced by $\wt{\eta}_\bullet^d$ in Theorem \ref{thmA-3} (cf.\ \cite{T1} Theorem 3.3.2\,(1)). We define
\[\wt{h}^d : \qt R\vG\big(\ol Y, \cSbt(d)_{\logB}\big) \lra \qt R\vG\big(\ol Y,\bE_{\logBb/\Wb} \big) \] in the same way as $h^r$ (using \eqref{eqA-cant} instead of \eqref{eqA-can}) and define $g^d := p^d \cdot (\wt{h}^d \circ (\wt{f}^d)^{-1})$. Finally we define $g^r$ as the composite of the natural map
\[ R\vG\big(\Xggf, \mu_{p^\bullet}^{\otimes r}\big) \lra R\vG\big(\Xggf, \bZpd(d)'\big) \otimes \bZpd(r-d)' \]
and $g^d \otimes \id$. The above diagram is commutative by the definition of $g^r$ and the compatibility of $\eta_\bullet^r$ with products (cf.\ \cite{T1} \S3.1). Now we prove
\addtocounter{athm}{3}
\begin{athm}\label{propA-4}
Assume that $K$ is a $p$-adic local field {\rm(}i.e., $k$ is finite{\rm)}. Then the diagram \eqref{eq8-4} commutes for $r<d$.
\end{athm}
\begin{pf}
We first note the isomorphisms
\begin{align*}
H^{i+1}_\cosyn(\logA,\cS_\qp(r)) & \isom H^{i+1}_\syn(\logA,\cS_\qp(r)), \\
H^{i+1}_\cont(\Xgf,\qp(r)) & \isom H^{i+1}(\Xgf,\qp(r))
\end{align*}
by the assumption that $k$ is finite (cf.\ Proposition \ref{propA-5}\,(2)). For integers $m,s \ge 0$, put
\begin{align*}
 H^m_\dccrys(\logA/W;s)&:= R^m\big(\varprojlim \vG_\Gal\big) \big\{R\vG\big(\ol Y,\bE_{\logBb/\Wb}\big)\otimes\bZ/p^\bullet(s)'\big\},\\
 H^m_\dccrys(\logA/W;s)_\qp^0 &:=\ker(H^m_\dccrys(\logA/W;s)_\qp \to H^m_\crys(\logB/W)_\qp(s)).
\end{align*}
By \eqref{eqA-3}, there is a commutative diagram
\[ \xymatrix{
H^{i+1}_\syn(\logB,\cS_\qp(r)) \ar[d]_{f^r} \ar[rd]^{h^r} \\
H^{i+1}(\Xggf, \qp(r))_\qp \ar[r]^-{g^r} & H^{i+1}_\crys(\logB/W)_\qp (r-d),} \]
where the bottom arrow is the same as $\beta^{i+1,r}$. By this diagram we obtain the arrows $(g^r)^0$ and $(h^r)^0$ in the following commutative diagram:
\[\xymatrix{
H^{i+1}_\syn(\logA,\cS_\qp(r))^0 \ar[d]_{(f^r)^0} \ar[rd]^{(h^r)^0} \\
H^{i+1}_\cont(X, \qp(r))^0 \ar[r]^-{(g^r)^0} \ar[d]_{\text{edge}} & H^{i+1}_\dccrys(\logA/W;r-d)_\qp^0 \ar[d]^{\text{edge}} \\
H^1(K,H^i(\Xggf, \mu_{p^\bullet}^{\otimes r}))_\qp \ar[r]^-{g^r} & H^1(K,H^i_\crys(\logBb/\Wb) \otimes \bZpd (r-d))_\qp \\
H^1(K,H^i(\Xggf, \qp(r))) \ar[u]_{\hspace{-2pt}\wr}^c \ar[r]^-{\beta^{i,r}} & H^1(K,H^i_\crys(\logB/W)_\qp (r-d)), \ar[u]^{\wr\hspace{-2pt}}_{c'} }\]
where the top triangle is induced by \eqref{eqA-3} and the central square commutes by the functoriality of Hochschild-Serre spectral sequences. The arrows $c$ and $c'$ are canonical maps, and the bottom square commutes by the definitions of $\beta^{i,r}$ and $g^r$. The map $c$ is bijective by \cite{J} Theorem 2.2 and the finiteness of $H^i(\Xggf, \mu_{p^n}^{\otimes r})$ for $n \ge 1$. See the proof of Theorem \ref{thmA-1} for the bijectivity of $c'$. Moreover it is easy to see that the composite of the left column agrees with the left vertical arrow of \eqref{eq8-4}, and that the composite of $h^r$ and the right column agrees with $\gamma^{i,r}$ in \eqref{eq8-4}. The commutativity in question follows from these facts.
\end{pf}


\end{document}